\documentclass[10pt]{article}
\RequirePackage{sty_packages}

\begin{document}

	\title{Saecular persistence}
	\author{Robert Ghrist, Gregory Henselman-Petrusek}
	\maketitle	


\let\amsamp=&



\makeatletter 
\newtheorem*{rep@theorem}{\rep@title}
\newcommand{\newreptheorem}[2]{%
\newenvironment{rep#1}[1]{%
 \def\rep@title{#2 \ref{##1}}%
 \begin{rep@theorem}}%
 {\end{rep@theorem}}}
\makeatother


\newtheorem{theorem}{Theorem}
\newreptheorem{theorem}{Theorem}
\newtheorem{proposition}[theorem]{Proposition}
\newtheorem{corollary}[theorem]{Corollary}
\newtheorem{lemma}[theorem]{Lemma}
\newtheorem{definition}{Definition}
\newtheorem{example}{Example}
\newtheorem{counterexample}{Counterexample}
\newtheorem{remark}{Remark}


\newcommand{\GHP}[1]{\textcolor{red}{#1}}


\newcommand{\R}{\mathbf{R}}
\newcommand{\Z}{\mathbf{Z}}
\newcommand{\C}{\mathbf{C}}
\newcommand{\N}{\mathbf{N}}
\newcommand{\Q}{\mathbf{Q}}
\newcommand{\E}{\mathbf{E}}
\newcommand{\F}{\mathbf{F}}
\newcommand{\A}{\mathbf{A}}

\newcommand{\e}{=}
\newcommand{\p}{+}
\newcommand{\m}{-}

\newcommand{\be}{\begin{enumerate}}
\newcommand{\bi}{\begin{itemize*}}
\newcommand{\ei}{\end{itemize*}}
\newcommand{\dd}[2]{\frac{d {#1}}{d {#2}}}
\newcommand{\dn}[2] {\left. \frac{\partial}{\partial {#1}} \right |_{#2}}
\renewcommand{\d}[2]{\frac{\partial {#1}}{\partial {#2}}}
\newcommand{\da}[3]{\frac{\partial {#1}}{\partial {#2}}  \bigg|_{#3}   } 
\newcommand{\dy}[1]{\frac{\partial }{\partial {#1}}}
\newcommand{\la}[1]{\displaystyle{\cal{L}}\left\{{#1}\right\}}
\newcommand{\li}[1]{\displaystyle{\cal{L}}^{-1}\left\{{#1}\right\}}
\newcommand{\ls}{\left \{} 
\newcommand{\rs}{\right \}}
\newcommand{\lb}{\left[}
\newcommand{\rb}{\right]}
\newcommand{\lcc}{\left (}
\newcommand{\rcc}{\right)}
\newcommand{\ba}{\begin{array}}
\newcommand{\ea}{\end{array}}
\newcommand{\bca}{\begin{cases}}
\newcommand{\eca}{\end{cases}}
\newcommand{\bea}{\begin{eqnarray}}
\newcommand{\beaa}{\begin{eqnarray*}}
\newcommand{\eea}{\end{eqnarray}}
\newcommand{\eeaa}{\end{eqnarray*}}
\newcommand{\bmat}{\left[\begin{array}}
\newcommand{\emat}{\end{array}\right]}
\newcommand{\bt}{\begin{thm}}
\newcommand{\et}{\end{thm}}
\newcommand{\bp}{\begin{proof}}
\newcommand{\ep}{\end{proof}}
\newcommand{\bprop}{\begin{prop}}
\newcommand{\eprop}{\end{prop}}
\newcommand{\bl}{\begin{lemma}}
\newcommand{\el}{\end{lemma}}
\newcommand{\bc}{\begin{cor}}
\newcommand{\ec}{\end{cor}}
\newcommand{\bd}{\begin{definition}}
\newcommand{\ed}{\end{definition}}
\newcommand{\ov}{\overline}
\newcommand{\ul}{\underline}
\newcommand{\limm}[2]{\displaystyle \lim_{#1 \to #2}}
\newcommand{\fr}[2]{\displaystyle \frac{#1}{#2}}
\newcommand{\series}[3]{ \displaystyle \sum_{{#1}={#2}}^{#3} }
\newcommand{\intt}[2]{\displaystyle \int_{#1}^{#2}}
\newcommand{\imt}[3]{\int_{#1} #2 \hspace{.1cm} \textmd{d} #3}
\newcommand{\union}[3]{\displaystyle \cup_{{#1}={#2}}^{#3} }
\newcommand{\spaceit}{\hspace{0.5cm}}
\newcommand{\smace}{\hspace{0.25cm}}
\newcommand{\iss}{\hspace{0.15cm} \exists \hspace{0.15cm}}
\newcommand{\foral}{\hspace{0.15cm} \forall \hspace{0.15cm}}
\newcommand{\is}{\hspace{0.15cm} \exists \hspace{0.15cm}}
\newcommand{\teorem}[1]{\fbox{T} \spaceit \textbf{#1} \spaceit}
\newcommand{\cor}[1]{\fbox{C} \spaceit \textbf{#1} \spaceit}
\newcommand{\deff}[1]{\fbox{D} \spaceit \textbf{#1} \spaceit}
\newcommand{\ex}[3]{\textbf {Exercise #1} \hspace {.5 cm}  #2 \\ \textit{ Solution} \\  \hspace{.5cm} \\  #3 }
\newcommand{\sa}[1]{\indent \hspace{0.#1cm}}
\newcommand{\spa}[1]{\indent \hspace{#1cm}}
\newcommand{\dm}{\hspace{.1cm} d \mu}
\newcommand{\inn}{\hspace{.1cm} \in \hspace{.1cm}}
\newcommand{\bma}{\left( \begin{array}}
\newcommand{\ema}{ \end{array} \right )}

\newcommand{\y}[1] {\widetilde{#1}}
\newcommand{\h}[1] {\widehat{#1} }

\newcommand{\tcb}[1]{\textcolor{blue}{#1}}
\newcommand{\tcr}[1]{\textcolor{red}{#1}}
\newcommand{\tcy}[1]{\textcolor{yellow}{#1}}
\newcommand{\tcg}[1]{\textcolor{green}{#1}} 

\newcommand{\f}{\mathfrak}
\newcommand{\tj}[1] {\mathscr{#1}}
\newcommand{\ti}[1]{\textit{#1}}
\newcommand{\mt}[1]{\mathtt{#1}} 
\newcommand{\mf}[1]{\mathfrak{#1}} 

\newcommand{\into}{\hookrightarrow}
\newcommand{\onto}{\twoheadrightarrow}

\renewcommand{\b}{&}
\newcommand{\bij}{\leftrightarrow}
\newcommand{\pd}[1] {^{\mathrm{{#1}}}}
\newcommand{\com}{\circ}
\newcommand{\sln}{\ti{Solution} \sa 2}
\newcommand{\ci}[2]{[{#1}, {#2}] }
\newcommand{\oi}[2]{({#1}, {#2}) }
\newcommand{\ho}[2]{[{#1}, {#2})}
\newcommand{\oh}[2]{({#1}, {#2}] } 
\newcommand{\ergo}{\hspace{.4cm} \therefore \hspace{.4cm}}
\newcommand{\der}[1] {^{(#1)}}
\newcommand{\rt}[1]{\sqrt{#1}} 
\newcommand{\nsg}{\trianglelefteq}
\newcommand{\w}{\wedge} 
\newcommand{\gm}{{\bigwedge}}
\newcommand{\lr}[1] {\overset{#1} {\longrightarrow}}
\newcommand{\lra}{\longrightarrow}
\newcommand{\pp}{\oplus}
\newcommand{\opp}{\oplus \cdots \oplus}
\newcommand{\te}[1]{\otimes_{#1}}
\newcommand{\tn}{\otimes}
\newcommand{\de}[1] {\indent \hspace {.5cm} \tb{(#1)} }
\newcommand{\cd}{\cdots}
\newcommand{\ld}{\ldots}
\newcommand{\lsu}[3]{{#1}^{#2}, \ld, {#1}^{#3}} 
\newcommand{\lst}[2]{#1_1, \ldots, #1_{#2}}
\newcommand{\lsn}[3]{#1_{#2}, \ldots, #1_{#3}}
\newcommand{\lsp}[3]{#1_{#2} + \cdots + #1_{#3}}
\newcommand{\lsm}[4] {#1_{#2}#1_{#3} \cdots #1_{#4}}
\newcommand{\isom}{\hspace{.2 cm} \tilde \rightarrow \hspace{.2cm}}
\newcommand{\st}{\text{ such that }}
\newcommand{\nn}[1]{\Vert #1 \Vert} 
\newcommand{\nr}[1] {\langle {#1} \rangle }
\newcommand{\lcm}{\mathrm{lcm}}
\newcommand{\bs}{\hspace{0.05cm}\backslash\hspace{0.05cm}}
\newcommand{\ot}{\rightharpoonup}
\newcommand{\ee}[1]{ \backepsilon #1 }
\newcommand{\cha}[1]{ \chi{_{_{#1}}}}
\newcommand{\challenge}[2]{\setcounter{equation}{0}   \hspace{2.7cm} --------------------------------- \sa 5   \tb{Exercise {#1}} \sa 5 ---------------------------------   \vspace{.25cm}  \hfill \\  #2  \vspace{.25 cm} \\ . \hspace{3.7cm} ------------------------ \sa 5   Solution \sa 5 ------------------------   \vspace{.25cm}  \hfill .\\ }
\newcommand{\challeng}[2]{\setcounter{equation}{0} \vspace{.3cm}  \textbf{#1} \hspace{.1cm}  #2 $\Diamond$ \\   }
\newcommand{\si}[2] {_{#1_{#2}}}
\newcommand{\ii}[2]{_{{#1}{#2}}}
\newcommand{\fm}[3]{_{#1 = #2}^{#3}} 
\newcommand{\inv}{^{-1}}
\newcommand{\sgn}{\mathtt{sgn}}
\newcommand{\prf}{\hspace{.2cm} \ti{Proof.} \hspace{.2cm}}
\newcommand{\barr}[1] {\left(\begin{array}{#1}}
\newcommand{\earr}{\end{array}\right)}

\newcommand{\supp}{\mathrm{supp}}
\newcommand{\ahom}{\mathrm{AffHom}}
\newcommand{\Fr}{\mathrm{Fr}}
\newcommand{\Rad}{\mathrm{Rad}}
\newcommand{\res}[2] {\mathrm{res}^{#1}_{#2}}
\newcommand{\coind}[2] {\mathrm{coind}^{#1}_{#2}}
\newcommand{\ind}[2] {\mathrm{ind}^{#1}_{#2}}
\newcommand{\Sym}{\mathrm{Sym}}
\newcommand{\spann}{\mathrm{span}}
\newcommand{\Com}{\mathrm{Com}}
\newcommand{\diag}{\mathrm{diag}}
\newcommand{\Ann}{\text{Ann }}
\newcommand{\tr}{\mathrm{tr }}
\newcommand{\Fix}{\mathrm{Fix}}
\newcommand{\Gal}{\mathrm{Gal}} 
\newcommand{\irr}{\mathrm{irr}}
\newcommand{\xar}{\mathrm{char}}
\newcommand{\Frac}{\mathrm{Frac}}
\newcommand{\Ob}{\mathrm{Ob }}
\renewcommand{\hom}[1] {\mathrm{{Hom$_{#1}$}}} 
\newcommand{\Hom}{\mathrm{Hom}}
\newcommand{\End}{\mathrm{End}}
\newcommand{\im}{\mathrm{im }}
\newcommand{\Ab}{\mt{Ab}}
\newcommand{\FinVect}{\mt{FinVect}}
\newcommand{\Vect}{\tt{\mathrm{-}Vect}}
\newcommand{\Groups}{\mt{Groups}}
\newcommand{\Mod}{\text{-}\mt{Mod}}
\newcommand{\Rings}{\mt{Rings}}
\newcommand{\Sets}{\mt{Sets}}
\newcommand{\PO}{\mt{PO}}
\newcommand{\Top}{\mt{Top}}
\newcommand{\HoTop}{\mt{HoTop}}
\newcommand{\Rep}{\mt{Rep}} 
\newcommand{\MatRep}{\mt{MatRep}}
\newcommand{\Alg}{\mathrm{-}\mt{Alg}}
\newcommand{\Fields}{\mt{Fields}}
\newcommand{\mmod}{\mathrm{-}\mt{mod}}

\newcommand{\sym}{\mathrm{Sym}}
\newcommand{\Aut}{\mathrm{Aut}}
\newcommand{\id}{\mathrm{id}}
\newcommand{\idd}[1]{\mathrm{id}_{#1}}
\newcommand{\Inn}{\mathrm{Inn }}
\newcommand{\Out}{\mathrm{Out}}
\newcommand{\Fin}{\mathrm{Fin}}
\newcommand{\Syl}{\mathrm{Syl}}
\newcommand{\ISO}{\mathrm{ISO}}


\renewcommand{\H}{\mathbb {H}}
\renewcommand{\ss}[2]{^{1}{2}} 
\newcommand{\cp}[1]{\mathbb {C}P^{#1}}
\newcommand{\rp}[1]{\mathbb {R}P^{#1}} 
\renewcommand{\c}[1]{\mathbb{C}^{#1}} 
\renewcommand{\r}[1] {\mathbb{R}^{#1}}


\newcommand{\holdsfor}{\mathrm{Y}} 

\renewcommand{\c}[1] {\mathcal{#1}}
\newcommand{\mm}{\mathcal{M}}
\newcommand{\mc}{\mathcal{C}}
\newcommand{\ma}{\mathcal{A}}
\newcommand{\mb}{\mathcal{B}}
\newcommand{\mg}{\mathcal{G}}

\newcommand{\page}[4]{ \ovalbox{#3} \hspace{.5cm} \textbf{#4} \marginpar{\textbf{ #1}.\textit{#2}}}

\newcommand{\se}{\simeq}
%

\newcommand{\su}{\subseteq}
\newcommand{\sn}{\subsetneq}
\newcommand{\ct}{\supseteq}
\newcommand{\qs}{ \subset^ {\varpropto} }  
\newcommand{\cs}{\hspace{.1cm}  \ov {\subset} \hspace{.1cm}  }  
\newcommand{\os}{\hspace{.1cm}  \mathring{\subset}  \hspace{.1cm}}  
\newcommand{\ox}{^{\circ}}
\newcommand{\du}{\coprod}
\newcommand{\oss}[1]{^{ \circ #1 }}
\newcommand{\css}[1]{^{- #1}}

\newcommand{\lp}[2]{\lVert #2 \lVert_{#1}}

\newcommand{\real}[1]{|{#1}_\bullet|}

\newcommand{\ak}{\alpha}
\newcommand{\bk}{\beta}
\newcommand{\ck}{\gamma}
\newcommand{\dk}{\delta}
\newcommand{\ek}{\varepsilon}
\newcommand{\ke}{\epsilon}
\newcommand{\fk}{\varphi}
\newcommand{\kf}{\phi}
\newcommand{\hk}{\eta}
\newcommand{\ik}{\iota}
\newcommand{\kk}{\kappa}
\newcommand{\kl}{\varkappa}
\newcommand{\lk}{\lambda}
\newcommand{\mk}{\mu}
\newcommand{\nk}{\nu}
\newcommand{\pk}{\pi}
\newcommand{\kp}{\varpi}
\newcommand{\qk}{\theta}
\newcommand{\kq}{\vartheta}
\newcommand{\rk}{\rho}
\newcommand{\kr}{\varrho}
\newcommand{\sk}{\sigma}
\newcommand{\ks}{\varsigma}
\newcommand{\uk}{\upsilon}
\newcommand{\tk}{\tau}
\newcommand{\wk}{\omega}
\newcommand{\xk}{\xi}
\newcommand{\zk}{\zeta}
\newcommand{\gk}{\chi}


\newcommand{\ag}{\Gamma}
\newcommand{\ad}{\Delta}
\newcommand{\aq}{\Theta}
\newcommand{\al}{\Lambda}
\newcommand{\ax}{\Xi}
\newcommand{\ap}{\Pi}
\newcommand{\as}{\Sigma}
\newcommand{\au}{\Upsilon}
\newcommand{\af}{\Phi}
\newcommand{\ay}{\Psi}
\newcommand{\aw}{\Omega}


\newcommand{\pf}{\sa 2 ti{Proof.} \sa 2}
\newcommand{\soln}{\sa 2 \ti{Solution.} \sa 2}



\newcommand{\iii}{\begin{enumerate*}}
\newcommand{\OOO}{\item}
\newcommand{\jjj}{\end{enumerate*}}


\newcommand{\tm}[1]{\text{#1}}
\newcommand{\tb}[1]{\textbf{#1}}
\def\MLine#1{\par\hspace*{-\leftmargin}\parbox{\textwidth}{\[#1\]}}


\newcommand{\defont}[1]{\emph{#1}}



\newcommand{\sbt}{\,\begin{picture}(1,1)(-1,-2)\circle*{2}\end{picture}\ }
\renewcommand{\b}{\bullet}
\newcommand{\di}{_{\sbt}}
\renewcommand{\ii}{^{\sbt}}

\newcommand{\sslash}{\mathbin{/\mkern-6mu/}}
\newcommand{\sslashh}{\mathbin{/\mkern-6mu/}^*}

\newcommand{\relativeindex}{relative index }
\newcommand{\module}{functor }
\newcommand{\persistencemodule}{chain functor }

\newcommand*{\dt}[1]{%
  \accentset{\mbox{\bfseries .}}{#1}}

\newcommand{\dz}[1]{{#1}_{>0}}

\newcommand{\statement}{\mathrm{N}}

\newcommand{\rind}{\rho} 
\newcommand{\psrel}[1]{\mathrm{Rel}(#1)}
\newcommand{\dcon}{\prec} 
\newcommand{\dconn}{\succ } 
\newcommand{\Dcon}{\triangleleft} 
\newcommand{\Dconn}{\triangleright} 
\newcommand{\rela}{r} 	

\newcommand{\clh}{{\mathrm{cl}}} 		
\newcommand{\pola}{\xi} 				
\newcommand{\polaa}{\dop \xi} 				

\newcommand{\unbounded}{\mathring}
\newcommand{\groundset}{\mathrm{Grs}}
\renewcommand{\lneq}{<}
\renewcommand{\gneq}{>}
\newcommand{\rel}[1]{\mathrm{Rel}(#1)}
\newcommand{\sgt}{>}
\newcommand{\slt}{>}
\newcommand{\lina}{\sigma} 		
\newcommand{\linb}{\varsigma} 		
\newcommand{\linela}{l} 				
\newcommand{\pred}{\mathrm{pred}} 		
\newcommand{\suq}{\mathrm{succ}}		
\newcommand{\Suq}{\mathrm{Succ}}		
\newcommand{\cov}[1]{{\mathrm{Cov}}(#1)}
\newcommand{\itvseta}{\mathcal{I}} 		
\newcommand{\dphi}{\Phi} 				
\newcommand{\pseta}{\mathrm{P}} 		
\newcommand{\psela}{p} 					
\newcommand{\pselb}{q} 					
\newcommand{\pselc}{r} 					
\newcommand{\psetb}{\mathrm{Q}} 		
\newcommand{\psetc}{\mathrm{R}} 		
\newcommand{\psma}{f} 					
\newcommand{\psmb}{g} 					
\newcommand{\psmc}{h} 					
\newcommand{\psama}{g} 					
\newcommand{\psema}{\epsilon} 			
\newcommand{\psqua}{\eta} 				
\newcommand{\filta}{\mathrm{F}} 
\newcommand{\filtb}{\mathrm{G}} 
\newcommand{\filtc}{\mathrm{L}} 
\newcommand{\ds}[1]{\downarrow (#1)} 	
\newcommand{\dsh}{ \; {\downarrow} \;} 			
\newcommand{\pds}[1]{\; {\mathring \downarrow}\; (#1)} 
\newcommand{\pdsh}{  \; {\mathring \downarrow}\;} 
\newcommand{\pdshp}[1]{  \;\mathring {\downarrow}_{#1} \;} 
\newcommand{\us}[1]{\uparrow(#1)} 		
\newcommand{\ush}{\uparrow} 			
\newcommand{\pus}[1]{\mathring {\uparrow}(#1)} 
\newcommand{\push}{\mathring {\uparrow}}
\newcommand{\dsp}{\downarrow}   
\newcommand{\dso}{\centernot \uparrow} 	

\newcommand{\dseta}{D} 	
\newcommand{\dsma}{a} 		
\newcommand{\dsia}{d} 		

\newcommand{\mdom}{\Phi_\mathsf{M}} 
\newcommand{\mdomG}{\mdom^{\mathsf{Gp}}} 
\newcommand{\foralc}{\forall} 		
\newcommand{\foralcc}{\unbounded \forall} 	
\newcommand{\existsc}{\exists} 		
\newcommand{\existscc}{\unbounded \exists} 
\newcommand{\meet}{\wedge}		
\newcommand{\cdflin}{\Sigma}  	
\newcommand{\cdflinb}{\Sigma'}  	
\newcommand{\cdfm}{\mk} 		
\newcommand{\cdf}{\Lambda} 		
\newcommand{\cdfi}{\Lambda'} 	
\newcommand{\cdfb}{\Gamma} 		

\newcommand{\SBDF}{\mathsf{SBDF}}
\newcommand{\SBDIF}{\mathsf{SBDIF}}
\newcommand{\SCDF}{\mathsf{SCDF}}
\newcommand{\SCDIF}{\mathsf{SCDIF}}

\newcommand{\SBD}{\mathsf{SBD}}
\newcommand{\SBDI}{\mathsf{SBDI}}
\newcommand{\SCD}{\mathsf{SCD}}
\newcommand{\SCDI}{\mathsf{SCDI}}
\newcommand{\FBD}{\mathsf{BD}} 
\newcommand{\IFBD}{\mathsf{BDI}} 
\newcommand{\FCD}{\mathsf{CD}} 
\newcommand{\IFCD}{\mathsf{CDI}} 
\newcommand{\sringa}{R} 		
\newcommand{\xxl}[1]{\mathbb{X}(#1)}
\newcommand{\xxlh}{\mathbb{X}}
\newcommand{\lata}{\mathbf{L}} 	
\newcommand{\latb}{\mathbf{Q}} 	
\newcommand{\latc}{\mathbf{R}} 	
\newcommand{\latx}{\mathbf{X}} 	
\newcommand{\laty}{\mathbf{Y}} 	
\newcommand{\lara}{\lambda} 	
\newcommand{\dcla}{D} 
\newcommand{\ucla}{U} 
\newcommand{\uclsa}{\mathbf U} 
\newcommand{\toda}{F} 
\newcommand{\todb}{G} 
\newcommand{\ffsp}[1]{#1^{\boxplus}}   
\newcommand{\lela}{a} 
\newcommand{\lelb}{b} 
\newcommand{\lelc}{c} 
\newcommand{\leld}{d} 
\newcommand{\lele}{e} 
\newcommand{\pela}{x} 
\newcommand{\pelb}{y} 
\newcommand{\pelc}{z} 
\newcommand{\lssa}{A} 
\newcommand{\fmma}{e}	
\newcommand{\fmmae}{K}	
\newcommand{\fmmaee}{L}	
\newcommand{\xhom}{\mathrm{X}} 	
\newcommand{\bfl}{\dphi(\tob^\pm, \tob^\pm)} 	
\newcommand{\bfls}{\stl{\dz \tob \times \dz \toc}} 
\newcommand{\pto}{\mathrm{H}} 
\newcommand{\pmbf}{\mathbf{H}} 
\newcommand{\lxfa}[1]{\xxl{#1}} 
\newcommand{\lxfah}{\Omega} 
\newcommand{\tobcc}{\dop {\mathbf A}} 	
\newcommand{\tobc}{\mathbf A} 		
\newcommand{\toa}{\mathbf T} 		
\newcommand{\tob}{\mathbf I} 		
\newcommand{\toc}{\mathbf J}	 	
\newcommand{\tod}{\mathbf K} 		
\newcommand{\toe}{\mathbf L}	 	
\newcommand{\todla}{k} 				
\newcommand{\todlb}{l} 				
\newcommand{\toela}{\ell} 			
\newcommand{\tela}{i} 				
\newcommand{\telb}{j} 				
\newcommand{\telc}{k} 				
\newcommand{\teld}{l} 				
\newcommand{\prodatod}{\toa^{2 \downarrow}}
\newcommand{\tofa}{\mathbf{I}} 		
\newcommand{\tofb}{\mathbf{J}} 		
\newcommand{\fala}{i} 				
\newcommand{\falb}{i'} 		 		
\newcommand{\fbla}{j} 				
\newcommand{\fblb}{j'} 		 		
\newcommand{\dla}{\mathbf{M}} 		
\newcommand{\dlala}{m} 				
\newcommand{\cdla}{\mathbf{L}} 		
\newcommand{\cdlala}{u} 
\newcommand{\cdlb}{\mathbf{N}} 		
\newcommand{\cji}[1]{\cjih(#1)}
\newcommand{\cjih}{\mathrm{J}} 		
\newcommand{\kji}[1]{\mathbb{L}(#1)}
\newcommand{\cjii}[1]{\cjih({#1}^*)} 
\newcommand{\kjii}[1]{\mathbb{L}^*(#1)} 
\newcommand{\jila}{j} 				
\newcommand{\jilaa}{j^*} 			
\newcommand{\jilb}{k} 				
\newcommand{\jilbb}{k^*} 			
\newcommand{\jilc}{l} 				
\newcommand{\jilcc}{l^*} 			
\newcommand{\cjie}{\xi} 			
\newcommand{\basea}{B} 				
\newcommand{\basla}{\beta} 			
\newcommand{\cjiso}{{\phi}}			
\newcommand{\cjiiso}{{\phi^*}}		
\newcommand{\tw}{{\textstyle \bigwedge}} 	
\newcommand{\tv}{{\textstyle \bigvee}}		
\newcommand{\laema}{\epsilon} 		

\newcommand{\x}{\sigma}
\newcommand{\xx}{{\sigma^*}}
\newcommand{\LHS}{X}
\newcommand{\RHS}{Y}
\newcommand{\dummysec}{\mathrm{X}} 
\newcommand{\dumva}{\vd{\ss , \s}}
\newcommand{\dumvaa}{\vu{\ss, \s}}
\newcommand{\Xpp}{\vd{\cc_+ , \c_+}}
\newcommand{\Xpps}{\vu{\cc_-, \c_-}}
\newcommand{\fmeet}[2]{\vd{#1,#2}}
\newcommand{\fjoin}[2]{\vu{#1,#2}}
\newcommand{\X}{X_f^\wedge}
\newcommand{\XX}{X_f^\vee}
\newcommand{\vd}[1]{\lfloor #1 \rfloor}
\newcommand{\vu}[1]{\lceil #1 \rceil}
\newcommand{\ev}[1]{|#1|}

\newcommand{\chaina}{C}
\newcommand{\chainb}{D}
\newcommand{\chainc}{E}
\newcommand{\achaina}{A}
\newcommand{\achainb}{B}
\newcommand{\achainaa}{A^*}
\newcommand{\achainbb}{B^*}
\newcommand{\toz}[1]{Z[#1]}
\newcommand{\ttoz}[1]{\mathbf{Z}[#1]}
\newcommand{\tozz}[1]{Z[#1]^*}
\newcommand{\ttozz}[1]{\mathbf{Z}[#1]^*}
\newcommand{\toga}{\ag} 			
\newcommand{\ccompa}{k} 			
\newcommand{\ccompb}{l} 			
\newcommand{\complength}[1]{\mathrm{length}(#1)}

\newcommand{\cpp}[1]{\lambda^{#1}} 
\newcommand{\cpph}{\lambda} 
\newcommand{\fem}{\mk} 
\newcommand{\femm}{\dop \mk} 
\newcommand{\fjseta}{\cji{\mathbf{A}}}
\newcommand{\fjjseta}{\cjii{{\mathbf{A}}}}

\newcommand{\cut}{\mathrm{cut}} 		
\newcommand{\cuta}{c} 				
\newcommand{\cutad}{c^*} 			
\newcommand{\cutb}{d} 				
\newcommand{\cutc}{e} 				
\newcommand{\cutd}{k} 				
\newcommand{\cutdown}[1]{{#1}\di} 	
\newcommand{\cutup}[1]{{#1}\ii} 	
\newcommand{\cutset}[1]{\mathrm{Cut}({#1})} 
\newcommand{\cutsetpm}[2][]{#2^{\pm #1}}
\newcommand{\fpm}[1]{{#1}^{\pm}} 	
\newcommand{\sgna}{\sigma} 			
\newcommand{\sgnb}{\tau} 			
\newcommand{\sgnc}{\upsilon} 		
\newcommand{\signof}[1]{\mathrm{sgn}(#1)}

\newcommand{\iid}{\iota} 
\newcommand{\iclosure}{(\iid \di \iid \ii )}
\newcommand{\allifun}{{\mathrm{Int}}} 
\newcommand{\allitv}[1]{\allitvh(#1)} 
\newcommand{\allitvh}{{\mathbb{I}}} 
\newcommand{\alljtv}[1]{\alljtvh(#1)} 
\newcommand{\alljtvh}{{\mathbb{J}}} 
\newcommand{\ifun}{\mathbb{I}} 				
\newcommand{\sarel}{\mathcal{S}} 				
\newcommand{\ifama}{\mathcal{I}} 				
\newcommand{\ifamb}{\mathcal{J}} 				
\newcommand{\iifama}{I} 				
\newcommand{\iifamb}{J} 				
\newcommand{\iimfa}{b} 					
\newcommand{\iimfb}{c} 					
\newcommand{\halfu}[1]{\phi(#1)} 		
\newcommand{\supu}[1]{\psi(#1)}  		
\newcommand{\itva}{{I}} 			
\newcommand{\itvb}{{J}} 			
\newcommand{\itpa}{\omega} 				
\newcommand{\itpb}{\upsilon} 			

\newcommand{\inta}{m}				
\newcommand{\intb}{n}				
\newcommand{\intc}{l}				
\newcommand{\intmaxa}{N} 			
\newcommand{\intmaxb}{M} 			
\newcommand{\rab}[2]{\{#1, \ld, #2\}} 

\newcommand{\ldiff}{d}
\newcommand{\lsum}{s}

\newcommand{\eqkk}{\eqref{eq_kerrpush} - \eqref{eq_kerrpull} }
\newcommand{\eqk}{\eqref{eq_kerpull} - \eqref{eq_kerpush} }

\newcommand{\genseta}{G}
\newcommand{\dgen}{\downarrow_\genseta} 	
\newcommand{\genla}{g} 				

\newcommand{\lcl}{\mathbb{W}} 		
\newcommand{\lclc}{\mathrm{Lcl}} 	
\newcommand{\efua}{F} 				
\newcommand{\Prp}[1]{\mathrm{Prp}(#1)}
\newcommand{\Nul}[1]{\mathrm{Nul}(#1)}
\newcommand{\ecata}{\mathsf{E}} 				
\newcommand{\ecatb}{\mathsf{F}} 				
\newcommand{\cecata}{\textsf{E}_0}	
\newcommand{\eumoda}{R_0} 			
\newcommand{\setofgraphmorphisms}{T} 

\renewcommand{\Im}{\mathrm{Im}} 		
\newcommand{\Imm}{\mathrm{Im}^*} 		
\renewcommand{\im}{\mathrm{im}} 		
\newcommand{\imm}{\mathrm{im}^*} 		
\newcommand{\Sp}{\mathrm{Sp}} 	 		
\newcommand{\Sq}{\mathrm{Sq}} 			
\newcommand{\Sub}{\mathrm{Sub}} 		
\newcommand{\Quo}{\mathrm{Quo}} 		
\newcommand{\SubQuo}{\mathrm{SubQuo}} 	
\newcommand{\Mlc}{\mathrm{Mlc}}			
\newcommand{\cmc}{\mathrm{CMlc}} 		
\newcommand{\euclidrep}{h}
\newcommand{\Ks}{\mathrm{K}}
\newcommand{\KKs}{\dop {\mathrm{K}}}
\renewcommand{\L}{L}  
\newcommand{\Ker}{\mathrm{Ker}}
\newcommand{\KKer}{\mathrm{Ker}^*}
\renewcommand{\ker}{\mathrm{ker}}
\newcommand{\Cst}{\mathrm{Cst}}
\newcommand{\Cok}{\mathrm{Cok}}
\newcommand{\coker}{\mathrm{coker}}
\newcommand{\kker}{\mathrm{ker}^*}
\newcommand{\Ksym}{\mathrm{K}}
\newcommand{\KKsym}{\mathrm{K}^*}
\newcommand{\ksym}{\mathrm{k}}
\newcommand{\kksym}{\mathrm{k}^*}

\newcommand{\catcx}{\mathsf{Ch}}
\newcommand{\catcd}{\mathsf{CD}}
\newcommand{\catps}{\mathsf{P}}
\newcommand{\catgroup}{\mathsf{Group}}
\newcommand{\catmod}{\mathsf{Mod}}
\newcommand{\cattop}{\mathsf{Top}}
\newcommand{\catgroupab}{\mathsf{AbGroup}}
\newcommand{\catset}{\mathsf{Set}}
\newcommand{\catsetp}{\mathsf{SetP}}
\newcommand{\catvect}{\mathsf{Vect}}
\newcommand{\catvectf}{\mathsf{VectFin}}
\newcommand{\catgsetpt}{\mathsf{T}}

\newcommand{\funca}{\mathsf{F}} 
\newcommand{\funh}{\mathsf{H}} 
\newcommand{\funz}{\mathsf{Z}} 
\newcommand{\funbd}{\mathsf{B}} 
\newcommand{\fdiagram}{\mathsf{Z}} 
\newcommand{\fdiagramI}{\fdiagram_{\mathbb{I}}} 
\newcommand{\funcforget}{\mathsf{Forget}} 
\newcommand{\xfrac}[2]{ \frac{|#1|}{|#2|}}
\newcommand{\sfun}{{\Omega}}
\newcommand{\sfuni}{{\Omega}_{\mathbb{J}}}
\newcommand{\op}{{^{\mathrm{ op}}}} 
\newcommand{\cata}{\mathsf C}				
\newcommand{\catb}{\mathsf D} 				
\newcommand{\objectclass}[1]{\mathrm{IsoCl}(#1)}	
\newcommand{\homset}[1]{#1_{1}} 	
\newcommand{\functa}{F}				
\newcommand{\functb}{G}				
\newcommand{\natransa}{\hk}			
\newcommand{\natransb}{\lambda} 	
\newcommand{\isoa}{\phi}  			
\newcommand{\isob}{\psi}  			
\newcommand{\isoc}{\theta} 			
\newcommand{\autoa}{\phi}  			
\newcommand{\mora}{f} 				
\newcommand{\morb}{g}  				
\newcommand{\morc}{h}  				
\newcommand{\subma}{\iota} 			
\newcommand{\oba}{A}				
\newcommand{\obb}{B}				
\newcommand{\obc}{C}				
\newcommand{\obd}{D}				
\newcommand{\soa}{N} 				
\newcommand{\sob}{M} 				
\newcommand{\quoma}{q} 				
\newcommand{\quomb}{p} 				
\newcommand{\icata}{\mathsf{I}} 				
\newcommand{\subquo}{\mathrm{Subquo}} 	
\newcommand{\dop}[1]{\hat {#1}}
\newcommand{\dople}{ \preceq }
\newcommand{\doplneq}{ \prec }
\newcommand{\dopiso}{\pi} 
\newcommand{\dopisoo}{\dop \pi}
\newcommand{\isoclasses}[1]{\mathrm{IsoClass}}

\newcommand{\lcf}{\Omega} 
\newcommand{\loclo}[1]{\Lambda(#1)} 	
\newcommand{\semitopa}{\Sigma}
\newcommand{\semitopb}{\Upsilon}
\newcommand{\semitop}{\Sigma}
\newcommand{\semitopof}[1]{\semitop(#1)}
\newcommand{\psst}[1]{\Gamma(#1)}
\newcommand{\psstr}[2]{\Gamma^{#2}(#1)}
\newcommand{\finprodsemtop}[1]{\mathrm{Fin^{\downarrow}(#1)}}
\newcommand{\toprodsemtop}[1]{\mathrm{H}_{#1}}
\newcommand{\ivl}[1]{{\mathbb{I}(#1)}} 	
\newcommand{\psl}[1]{{\mathbb P}(#1)}
\newcommand{\pslh}{{\mathbb P}}  
\newcommand{\stl}[1]{{\mathbb{S}(#1)}} 	
\newcommand{\stlh}{{\mathbb{S}}} 		
\newcommand{\dxl}[1]{{\dxlh(#1)}} 	
\newcommand{\dxlh}{\mathbb{D}} 			
\newcommand{\axl}[1]{{\axlh(#1)}} 	
\newcommand{\axlh}{\mathbb{A}} 			
\newcommand{\axela}{c}					
\newcommand{\axelaa}{c^*}				
\newcommand{\axelb}{Y} 					
\newcommand{\bxl}[1]{\bxlh(#1)} 	
\newcommand{\bxlh}{\mathbb{\unbounded A}} 			
\newcommand{\acl}[1]{ \aclh ({#1})} 		
\newcommand{\aclh}{{\mathbb{V}}} 	
\newcommand{\phl}[1]{\Omega(#1)} 		
\newcommand{\phlh}{\Omega} 			
\newcommand{\fpl}[1]{{\mathbb{G}(#1)}}
\newcommand{\fplh}[1]{{\mathbb{G}}}
\newcommand{\findn}[1]{\mathrm{Fin}^\downarrow(#1)}
\newcommand{\finup}[1]{\mathrm{Fin}^\uparrow(#1)}
\newcommand{\antic}[1]{\mathrm{Anti}(#1)}
\newcommand{\cla}{A}  					
\newcommand{\clb}{B}
\newcommand{\clc}{C}
\newcommand{\mseta}{E}  				

\newcommand{\fcxa}{\mathrm{A}} 
\newcommand{\cxa}{A}
\newcommand{\bda}{\delta}
\newcommand{\bdaa}{\partial}
\newcommand{\chainC}{{C}}
\newcommand{\chainF}{{\fl}}
\newcommand{\icb}[1]{{\mathbf B}_{#1}}
\newcommand{\icbname}{bifilter category  }
\newcommand{\icr}[1]{{\mathbf R}_{#1}}
\newcommand{\icrname}{reflective bifilter category }
\newcommand{\ics}[1]{{\mathbf S}_{#1}}
\newcommand{\ils}{\mathbf S} 
\newcommand{\icsname}{sequential category }
\newcommand{\icc}[1]{{\mathbf C}_{#1}}
\newcommand{\iccname}{conjugate category }
\newcommand{\ice}[1]{{\Xi}_{#1}^*}
\newcommand{\icee}[1]{{\Xi}_{#1}^\uparrow}
\newcommand{\icename}{universal category }
\newcommand{\uxtoa}[1]{\mathbf{C}^*(\mathbf{D}_{#1})}  		
\newcommand{\uxtoaa}[1]{\mathbf{C}^\uparrow(#1)}  		
\newcommand{\uxtofa}[2]{\mathbf{C}_{#2}^\downarrow(#1)}  		
\newcommand{\uxoa}{\cxa_{\ice{}}}  			
\newcommand{\uxoaa}{\cxa^{\icee{}}}  		
\newcommand{\uxba}{\delta_{\uxtoa{\tob}}} 		
\newcommand{\uxbaa}{\partial^{\icee{}}}
\newcommand{\uxfa}{\filta_{\ice{}}} 			

\newcommand{\tcxa}{h} 						
\newcommand{\nbecxa}{g} 					

\newcommand{\chainsof}[1]{\mathbf{Ch}(#1)}
\newcommand{\chainsoff}[1]{\mathbf{Ch}^*(#1)}
\newcommand{\chainfunctor}[1]{\mathbf{Ch}({#1})}
\newcommand{\chainfunctorr}[1]{\mathbf{Ch}^*({#1})}

\newcommand{\cyclesa}{Z} 					

\DeclareRobustCommand{\coprod}{\mathop{\text{\fakecoprod}}}
\newcommand{\fakecoprod}{%
  \sbox0{$\prod$}%
  \smash{\raisebox{\dimexpr.9625\depth-\dp0}{\scalebox{1}[-1]{$\prod$}}}%
  \vphantom{$\prod$}%
}

\newcommand{\cdmapa}{\phi} 					
\newcommand{\cdmw}{\Lambda} 				
\newcommand{\cdow}{\Omega} 					
\newcommand{\cccd}[1]{(#1)} 				
\newcommand{\saq}{D} 						

\newcommand{\iseta}{A} 	
\newcommand{\inda}{\alpha} 		
\newcommand{\isetb}{B} 	
\newcommand{\indb}{\beta} 		
\newcommand{\indc}{\gamma} 		
\newcommand{\indd}{\ell} 	

\newcommand{\ela}{s} 	
\newcommand{\elb}{t} 	
\newcommand{\elc}{u} 	
\newcommand{\eld}{v} 	
\newcommand{\seta}{S} 
\newcommand{\setaa}{S^*} 
\newcommand{\setb}{T} 
\newcommand{\setc}{U} 
\newcommand{\setd}{V} 
\newcommand{\sete}{W} 
\newcommand{\setx}{X}
\newcommand{\sety}{Y}
\newcommand{\finss}[1]{{\mathcal F}(#1)} 
\newcommand{\sfuna}{f} 	
\newcommand{\partfa}{p} 	
\newcommand{\partfb}{q} 	
\newcommand{\eval}[1]{|#1|}  

\newcommand{\dgB}{G} 	
\newcommand{\dga}{f}  	
\newcommand{\dgb}{g} 	
\newcommand{\dgc}{h} 	
\newcommand{\dgd}{e} 	
\newcommand{\dgq}{k} 	
\newcommand{\dgca}{x}  	
\newcommand{\dgcb}{y}  	
\newcommand{\otype}[1]{\otypeh(#1)} 		
\newcommand{\otypeh}{\mathrm{Obj}} 		
\newcommand{\odual}{^*} 	

\newcommand{\spd}{\cpd^{\mathrm{sae}}} 		
\newcommand{\spdSET}{\spd_{\mathsf{cst}}} 		
\newcommand{\spdNCL}{\spd_{\mathsf{ncl}}} 		
\newcommand{\cpd}{\mathrm{PD}} 			
\newcommand{\gpda}{\mathrm{PD}^{\mathsf A}} 
\newcommand{\gpdb}{\mathrm{PD}^{\mathsf B}} 
\newcommand{\barcode}[1]{\mathrm{BC}_{#1}}			
\newcommand{\iindex}[1]{\Omega_{#1}} 		
\newcommand{\pda}[2]{P_{#1}(#2)} 			
\newcommand{\rpda}{V} 						
\newcommand{\rpdb}{W} 						
\newcommand{\prfa}[1]{\mathrm{R}_{#1}^\amita} 	
\newcommand{\prfb}[1]{\mathrm{R}_{#1}^\amitb} 	
\newcommand{\pdta}[1]{\gpda_{#1}}				
\newcommand{\pdtb}[1]{\gpdb_{#1}}				
\newcommand{\pitvc}{\mathrm{Dgm}}				
\newcommand{\poga}{\mathrm{G}} 					
\newcommand{\birth}[1]{\mathrm{birth}(#1)} 		
\newcommand{\death}[1]{\mathrm{death}(#1)} 		
\renewcommand{\t}{i}
\renewcommand{\tt}{{\dop i}}
\renewcommand{\a}{a}
\renewcommand{\aa}{{\dop a}}
\newcommand{\z}{s}
\newcommand{\zz}{\dop s}
\newcommand{\s}{s}
\renewcommand{\ss}{\dop s}
\renewcommand{\f}{f}
\newcommand{\fop}{f^{op}}
\newcommand{\ff}{f^*}
\newcommand{\g}{g}
\renewcommand{\gg}{\dop g}
\renewcommand{\c}{c}
\newcommand{\cc}{\dop c}
\renewcommand{\dd}{\dop d}
\renewcommand{\d}{d}
\renewcommand{\k}{k}
\renewcommand{\kk}{\dop k}

\newcommand{\iia}{\f(\aa \dople \a)\ii}
\newcommand{\dia}{\f(\aa \dople \a)\di}

\newcommand{\K}{\mathrm{K}_f}
\newcommand{\KK}{ \dop {\mathrm K}_f}
\newcommand{\KKd}{\mathrm{K}_{\dop f}}
\newcommand{\Lf}{{\mathrm{L}_f}}

\newcommand{\dimia}{\phi} 		
\newcommand{\dimib}{\psi} 		

\newcommand{\cutma}{L} 			
\newcommand{\isq}[1]{\beta(#1)} 	
\newcommand{\isqh}{\beta} 		

\newcommand{\ringa}{R} 						
\newcommand{\fielda}{{\mathds k }} 			
\newcommand{\groundfield}{{\mathbbm{k}}}	

\newcommand{\crmod}[1]{{#1}\text{-}\mathtt{Mod}}
\newcommand{\rmoda}{V} 			
\newcommand{\rmodb}{U} 			
\newcommand{\rmodsa}{H}			
\newcommand{\rmodsb}{K}  		
\newcommand{\rmodsc}{L} 		
\newcommand{\rmodla}{u} 		
\newcommand{\rmodlb}{v} 		
\newcommand{\jhv}[1]{\jhvh(#1)} 	
\newcommand{\jhvh}{\mathrm{JH}} 	
\newcommand{\smic}{\amits(\ecata)} 	
\newcommand{\abga}{G} 			
\newcommand{\absga}{H}			
\newcommand{\absgb}{K}  		
\newcommand{\absgc}{L} 			
\newcommand{\crup}[1]{C_{#1}} 	

\newcommand{\amita}{\mathsf{A}}
\newcommand{\amitb}{\mathsf{B}}
\newcommand{\amitj}{\mathsf{I}}
\newcommand{\amits}{\mathsf{S}}

\newcommand{\posts}{\seta}
\newcommand{\post}{\ela}
\newcommand{\qt}{h} 
\newcommand{\generator}{\zeta} 
\newcommand{\conset}{W}
\newcommand{\consela}{w}
\newcommand{\cutembed}[1]{\epsilon^{#1}}


\begin{abstract}

A persistence module is a functor $\dga: \tob \to \ecata$, where $\tob$ is the poset category of a totally ordered set. This work introduces \emph{saecular decomposition}: a categorically natural method to decompose $\dga$ into simple parts, called interval modules. Saecular decomposition exists under generic conditions, e.g., when $\tob$ is well ordered and $\ecata$ is a category of modules or groups.   This represents a substantial generalization of existing factorizations of 1-parameter persistence modules, leading to, among other things, persistence diagrams not only in homology, but in homotopy.

Applications of saecular decomposition include inverse and extension problems involving filtered topological spaces, the 1-parameter generalized persistence diagram, and the Leray-Serre spectral sequence.  Several examples -- including cycle representatives for generalized barcodes -- hold special significance for scientific applications.

The key tools in this approach are modular and distributive order lattices, combined with Puppe exact categories.

\begin{quote}
\end{quote}
\end{abstract}

\section{Framework}

\begin{quote}
    \tb{{s{\ae}culum}} \emph{(s{\ae}cul\=\i, n.)} : generation; a period of long duration; the longest fixed interval.
\end{quote}

Advances in topological data analysis have generated growing attention to functors of the form $\mora: \tob \to \ecata$, where $\tob$ is the poset category of a totally ordered set.  As the order-theoretic term for a totally ordered set is a \emph{chain}, we will call such objects \emph{chain diagrams} or \emph{chain functors}. For $\ecata$ the category of vector spaces and linear transformations, chain functors are called \emph{persistence modules}. 

Topological persistence began with a set of independent works getting at parametrized topological features \cite{frosini92,frosini99,RobinsTowards99,RobinsComputational02} and crystallized in the emergence of a persistence diagram (PD) for the homologies of filtered simplicial complexes \cite{ELZ2000,ELZTopological02}. This is not a functor, but a function from the plane to the non-negative integers encoding birth and death of homology classes. 

The first revolution in persistence came from a change of perspective, that the values of the PD correspond to multiplicities of indecomposable factors in a Krull-Schmidt decomposition of a persistence module \cite{CZC+barcode04,CZC+Persistence05}. From this foundation, richer aspects of persistence were appended, including the representation-theoretic interpretation, with relations to quivers and Gabriel's Theorem \cite{CSMZigzag09} providing both clarity and extension. 

The persistence module perspective was transformative for several reasons.  For scientific and industrial applications, it solved the fundamental problem of mapping persistent features back to data, thus making individual points in the PD and bars in the barcode interpretable.  For theorists, the module perspective provided the formal basis on which duality, stability, and other concepts later grew or expanded. In the process, the module perspective also expanded the domain where persistence was defined, and made the persistence diagram itself more readily computed. Thus our understanding of persistence advanced in two stages: first with the definition of the persistence diagram, and second with the decomposition of an underlying algebraic object.

In 2016, Patel launched a new iteration of this two-stage framework for persistence. In \cite{patel2018generalized}, he introduced the notion of a \emph{generalized} persistence diagram, which is well-defined for tame chain functors in any abelian or essentially small symmetric monoidal category. 
As the name suggests, the generalized PD represents a substantial advance over the classic persistence diagram of \cite{ELZTopological02}, which requires field coefficients. However, this new formulation of the generalized persistence diagram also provides no decomposition of the underlying functor itself.  

The purpose of the present work is to complement the cycle launched by Patel by providing the desired decomposition of chain functors. We do so for chain functors valued in abelian categories -- as well as Puppe exact categories and $\catgroup$, the category of groups and group homomorphisms.  In so doing, we provide new tools for mapping features back to data; for structural investigation of the persistence module; and for algebraic connections with mathematics more broadly.  
In particular, we provide the first notion of a cycle representative for a generalized persistent feature, and the first definition of a barcode for persistent homotopy.  We also formulate an enumeration theorem that relates the generalized persistence diagram to the Leray-Serre spectral sequence. A second enumeration theorem relates the saecular decomposition to the generalized PD.  This relationship is less rigid than the one that relates Krull-Schmidt decomposition of linear-coefficient persistence modules with classical PD's, but similarly informative.

The essential ingredients in this approach are Puppe exact (p-exact) categories and order lattices.  P-exact categories give a lightweight axiomatic framework for homological algebra which brings much of the essential structure into  focus.  Many of the ideas in this work came directly or indirectly from existing work on p-exact categories, particularly \cite{grandis12}.  Order lattices are deeply entwined in theory of exact categories, and provide much of the technical heavy lifting in the present work.  At the same time, they help to bring the main ideas into sharp relief.  Order lattices were also essential to extending the basic ideas of this work from p-exact categories to the (not  p-exact) category of \emph{groups}. The details of this extension are surprising, and relate to a remarkable structure theorem concerning half-normal bifiltrations.

\section{Results}

\newcommand{\At}{A}
\newcommand{\Bt}{B}
\newcommand{\Zt}{L}
\newcommand{\zt}{\lambda}
\newcommand{\lhomfamily}{\mathcal H}
\newcommand{\factormap}{c}
\newcommand{\factorsuperset}{\mathcal C}
\newcommand{\apatha}{\sigma} 
\newcommand{\apathb}{\tau} 
\newcommand{\apathc}{\upsilon} 
\newcommand{\sss}{\Sigma}  
\newcommand{\factor}{\mathrm{Fact}} 
\newcommand{\ncl}{\mathrm{ncl}} 

This section summarizes several of the main ideas introduced throughout this work.  All theorems are modifications or special cases of results introduced later; appendix \ref{sec_intro_proofs} contains complete proofs.  

To begin, fix a functor 
	$$
	\dga: \tob \to \ecata
	$$ 
where $\ecata$ is a category of (left or right) $R$-modules, and  $\tob$ is the posetal category of a well ordered set. We will consider more general types of source and target categories throughout; however, restricting to well-ordered sets and $R$-modules will simplify the language needed for a general summary.

Recall that an \emph{interval} in a poset $\pseta$ is a subset $\itva \su \pseta$ such that $b \in \itva$ whenever $\itva$ contains an upper bound of $b$ and a lower bound of $b$.  This notion coincides with the usual notion of interval, when $\tob$ is the real line.  A functor $\dgb: \tob \to \ecata$ is an \emph{interval diagram of support type $\itva$} if  $\dga_\a = 0$ for $\a \notin \itva$, and the restricted functor $\dga|_{\itva}$ carries each arrow to an isomorphism.  The family of \emph{nonempty} intervals in $\tob$ is denoted $\allitv{\tob}$.  The family of interval diagrams with support type $\itva$ is denoted $\allifun(\itva)$.

\subsection{Saecular homomorphisms and interval decomposition}  

\noindent Our first contribution is a formal definition of what it means to decompose the functor $\dga$ into a family of interval diagrams.  We do this with the aid of a uniquely defined lattice homomorphism $\sfun$, which we call the \emph{saecular homomorphism}, as per the following.

	\begin{itemize}
	\item \emph{Domain:}  Write  $\axl{\pseta}$ for the family of down-closed subsets of any poset $\pseta$. We regard $\axl{\pseta}$ as an order lattice, under inclusion.  
Given $\itva \in \allitv{\tob}$ and $\seta, \seta', \setb, \setb' \in \axl{\tob}$, let us define $\itva \sim (\seta, \setb) \iff \itva = \setb - \seta$, and define $(\seta, \setb) \le (\seta', \setb') \iff \seta \su \seta'$ and $\setb \su \setb'$.  We can then endow $\allitv{\tob}$ with a unique partial order such that $ \itva \le \itvb \iff \itva \sim (\seta, \setb) \le (\seta', \setb') \sim \itvb$ for some $\seta, \seta', \setb, \setb' \in \axl{ \tob}$.  The domain of $\sfun$ is $\axlh{\allitv{\tob}}: = \axl{\allitv{\tob}}$.
	\item \emph{Codomain:} Recall that the category $[\tob, \ecata]$ of functors $\tob \to \ecata$ is abelian. The functor $\dga$ is an object in $[\tob, \ecata]$, and the family of subobjects of $\dga$, denoted $\Sub(\dga)$ or $\Sub_\dga$, is a complete order lattice under inclusion.  In particular, every subset of $\Sub_\dga$ has a meet and a join.  	The codomain of $\sfun$ is $\Sub(\dga)$.
	\end{itemize}

\noindent Recall that a homomorphism of order lattices is \emph{complete} if it preserves arbitrary meets and joins, and, in addition, its domain and codomain are complete lattices. To streamline the definition of $\sfun$, we denote downsets as $\dsh \itva = \{ \itvb: \itvb \le \itva\}$ and $\pdsh \itva = \{ \itvb: \itvb < \itva\}$ (for nonstrict and strict downsets respectively).

\begin{theorem}
\label{thm_sample_hom}
There exists a unique complete lattice homomorphism $\sfun: \axlh{\allitv{\tob}}\to \Sub(\dga)$ such that
	\begin{align}
	\frac{ \sfun  \dsh \itva \;  } { \sfun \pdsh \itva \;  }
	\in \allifun( \itva )
	\label{eq_sample_hom}
	\end{align}
for each $\itva \in \allitv{\tob}$.  
\end{theorem}

\begin{definition}
The \emph{saecular CDI homomorphism}  is the  complete lattice homomorphism defined in Theorem \ref{thm_sample_hom}.  We say that $\sfun$ decomposes $\dga$ into the indexed family of interval functors defined in Equation \eqref{eq_sample_hom}.
\end{definition}

The abbreviation \emph{CDI} stands for \emph{completely distributive interval} and stems from $\sfun$ being a free completely distributive homomorphism whose domain is an Alexandrov topology on a poset of intervals. The need for delineating various types of saecular lattice homomorphisms will arise later in \S\ref{sec_lattic_free_hom}.


\subsection{Saecular functors}  

\noindent The saecular CDI homomorphism can  be extended to a more informative object, called the \emph{saecular CDI functor}. To construct this functor, first define the  \emph{locally closed} category,  $\lclc(\allitv{\tob})$, as follows.
 Objects are sets of form $\seta - \setb$, where $\seta, \setb \in \axlh{\allitv{\tob}}$; equivalently, objects are locally closed sets, with respect to the Alexandrov topology on the poset $\allitv{\tob}$.  An arrow $X \to Y$ is a subset $S \su X \cap Y$ such that $S$ is up-closed in $X$ and down-closed in $Y$.  Composition of arrows $X \xrightarrow{S} Y \xrightarrow{T} Z$ is  intersection: $T \com S = T \cap S$.    
It can be shown that $\lclc(\allitv{\tob})$  has a zero object, kernels, and cokernels.
A functor $\mathsf{F}$ from $\lclc(\allitv{\tob})$ to any other category with zeros, kernels, and cokernels is  \emph{exact} it if preserves these universal objects.

\begin{theorem}
\label{thm_main_result_finite_fun}
There exists an exact functor  $\sfun: \lclc (  \allitv{\tob}) \to [\tob, \ecata]$  such that 
	$$
	\sfun \{ \itva \} \in \allifun( \itva)
	$$ 
for each nonempty interval $\itva$.  
This functor is unique up to unique isomorphism, if we require  $\sfun$ to restrict to a complete lattice homomorphism  $\axlh{\allitv{\tob}} \to \Sub(\dga)$.  In this case $\sfun|_{\axlh{\allitv{\tob}}}$ is the saecular CDI homomorphism.
\end{theorem}

\begin{definition}
The \emph{saecular CDI functor} of $\dga$ is the unique (up to unique isomorphism) functor defined in Theorem \ref{thm_main_result_finite_fun}.  
\end{definition}

\begin{remark}
\label{rmk_sample_cd_finite_existence}
When $\tob$ is finite, we can replace $\ecata$ with any Puppe exact category in Theorems \ref{thm_sample_hom} and \ref{thm_main_result_finite_fun}.  Puppe exact categories include all abelian categories, as well as many non-additive categories.  If we require  $\Sub(\dga_\a)$ to have finite height for each $\a \in \tob$, then $\tob$ can be any totally ordered set. Several other relaxations are available, in addition.
\end{remark}

\subsection{Saecular functors for nonabelian groups} 

\noindent If $\ecata$ were the category $\catgroup$, rather than a category of $R$ modules, then there would be at least two reasonable ways to interpret a fraction $\soa/\sob$ for any nested pair of subdiagrams $\sob \su \soa \su \dga$.  

In the first interpretation, $\soa/\sob$ is the \emph{cokernel} of the inclusion $\sob \su \soa$.  This is the sense in which we write the fraction $\xfrac{ \sfun \dsh \itva } { \sfun \pdsh \itva }$ in Equation \eqref{eq_sample_hom}.  The fraction can also be interpreted as the quotient of $\soa$ by $\sob$.  For concreteness, we sometimes write this as $\Cok (\soa/\sob)$.  
In the second interpretation, $(\soa/\sob)_\a$ is the family of cosets of $\sob_\a$ in $\soa_\a$, for each $\a \in \tob$; then $\soa/\sob$ is a diagram of pointed sets, where the distinguished point is the 0-coset.  For concreteness, we sometimes write this as $\Cst (\soa/\sob)$.  Note that $\Cst (\soa/\sob)$ inherits group structure from $\soa$ via the usual construction of quotient groups iff $\sob_\a \unlhd \soa_\a$ for all $\a$, in which case $\Cst (\soa/\sob) = \Cok (\soa/\sob)$.
Remarkably,  Theorem \ref{thm_sample_hom} remains true under the coset interpretation, and true (up to loss of uniqueness) under the cokernel interpretation.

\begin{theorem}
\label{thm_sample_hom_group}
Theorem \ref{thm_sample_hom} remains true when $\ecata= \catgroup$, if we interpret $\frac{ \sfun  \dsh  \itva \;  } { \sfun \pdsh \itva \;  }$ as $\Cst \left (\frac{ \sfun  \dsh  \itva \;  } { \sfun \pdsh \itva \;  } \right)$.  
Moreover, if $\sfun$ is the unique homomorphism that satisfies Theorem \ref{thm_sample_hom} under this interpretation, then 
	\begin{align}
	\Cok \left ( \frac{\sfun \dsh \itva}{\sfun \pdsh \itva} \right ) \in \allifun(\itva)
	\label{eq_intro_cok_factor}
	\end{align}
for each interval $\itva$.
\end{theorem}

\begin{definition}	
We call $\Cok \left ( \frac{\sfun \dsh \itva}{\sfun \pdsh \itva} \right )$  and $\Cst \left ( \frac{\sfun \dsh \itva}{\sfun \pdsh \itva} \right )$ the \emph{$\itva^{th}$ saecular factor} of $\dga$ and the $\itva^{th}$ saecular \emph{coset} factor of $\dga$, respectively.
\end{definition}

\begin{remark}
Theorem \ref{thm_main_result_finite_fun} also has an analogue for groups, but the statement is  more involved.
\end{remark}

\section{Applications}

This section summarizes several applications to existing problems and theories in related fields.  Appendix \ref{sec_intro_proofs} contains complete proofs.

\subsection{Barcodes, persistence diagrams, and inverse problems}

The \emph{barcode}  of a chain diagram $\dgb: \R \to \ecata$ is a multiset of intervals, i.e.\  elements of $\allitv{\R}$.  
It is well defined whenever $\ecata = \fielda \catvectf$, the category of finite dimensional $\fielda$-linear vector spaces, for some field $\fielda$.   We denote the barcode  $\barcode{ \dgb }$. 

The associated \emph{persistence diagram} (PD) \cite{ELZTopological02} is the function $\cpd_\dgb: \R \times \R \to \Z$ sending $(p,q)$ to the multiplicity of interval $[p,q)$ in $\barcode { \dgb }$.  Barcodes and persistence diagrams carry identical information when $\barcode{ \dgb }$ contains only left-closed, right-open intervals, but each introduces minor impositions in formal arguments; multisets entail a variety of semantic subtleties, and the notion of a persistence diagram fails to capture information about intervals of form $(a,b), (a,b]$, and $[a,b]$.
Rather than introduce a third term, we break from convention, and use both persistence diagram and barcode to refer to the function $\allitv{\R} \to \Z$ sending $\itva$ to the number of copies of $\itva$ in $\barcode{ \dgb }$.  It has been shown \cite{gabriel1972unzerlegbare, ZCComputing05, botnan2020decomposition} that the barcode uniquely determines the isomorphism class of $\dgb$; indeed, one can  generalize this result to chain functors with any totally ordered index category.  

The notion of a \emph{generalized persistence diagram} \cite{patel2018generalized} is a natural abstraction of the persistence diagram, which is suitable for a wider range of target categories $\ecata$.  Generalized PD's come in two varieties.  A type-$\amita$ generalized PD is a function $\gpda_\dgb: \R^2 \to \Z^{W}$, where $W$ is the set of object isomorphism classes in $\ecata$. Type-$\amita$ PD's are well defined for any tame (meaning constructible, \emph{cf.} \cite{patel2018generalized}) functor $\dgb$ from $\R$ to an essentially small symmetric monoidal category with images.  A type-$\amitb$ generalized PD is a function $\gpdb_\dgb: \R^2 \to \Z^{V}_{\ge 0}$, where $V$ is the set of isomorphism classes of simple objects  in $\ecata$ (meaning objects with no proper nonzero subobjects), and $\Z_{\ge 0}$ is the set of nonnegative integers.  Type-$\amitb$ PD's are well defined for any constructible functor from $\R$ to an abelian category.

The \emph{saecular persistence diagram}, or the \emph{saecular barcode}, is a fourth notion, distinct from $\cpd, \gpda$, and $\gpdb$.  It is well defined for any chain functor $\dga: \tob \to \ecata$ where $\tob$ is well ordered and $\ecata \in \{ R \catmod, \catmod R, \catgroup\}$; in fact, it can be defined much more generally.

\begin{definition}
The \emph{saecular barcode} (or \emph{saecular persistence diagram}) of $\dga$, is the function $\spd_\dga$  that sends an interval $\itva$ to the $\itva^{th}$ saecular factor of $\dga$.  The saecular \emph{coset} barcode is defined similarly, when $\ecata = \catgroup$.
\end{definition}

The saecular barcode $\spd_\dgb$ agrees with the classical  $\cpd_\dgb$ in the following sense:

\begin{theorem}
\label{thm_sample_cpd_agrees}
Let $\toc$ be totally ordered set and $\dgb: \toc \to \fielda \catvectf$ be a functor. Then  $\cpd_\dgb$ and $\spd_\dgb$ are well defined, and
	\begin{align*}
	\cpd_\dgb(\itva) = \dim( \spd_\dgb(\itva))
	\end{align*}
for each $\itva \in \allitv{\toc}$, where $\dim( \dgc ): = \max_{\a \in \tob} \dim(\dgc_\a)$.
\end{theorem}

The saecular barcode $\spd_\dgb$ also agrees with $\gpdb_\dgb$   in a natural sense (Theorem \ref{thm_sample_gpd_agrees}).  
To streamline the statement of this result, recall that the \emph{length} of a module -- or, more generally, the length of an object $\oba$ in an abelian category -- is the supremum of all $n$ such that $\oba$ admits a strictly increasing sequence of subobjects $0<\soa_1< \ldots < \soa_n = \oba$.  
When $\oba$ has finite length, we may define the Jordan-H\"older vector of $\oba$ to be unique element $\jhv{ \oba } \in \Z^V$ such that $\jhv{ \oba }_C = \#\{ i: \soa_{i}/\soa_{i-1} \cong C \}$ for each $C \in V$, where $\soa$ is any maximal-with-respect-to-inclusion chain of subobjects.  It can be shown that $\jhv{ \oba }$ does not depend on the choice of $\soa$ \cite[Theorem 1.3.5]{grandis12}.  

Given an interval functor $\dgb$ with support type $\itva$, we write $\jhv{ \dgb}$ for
the Jordan-H\"older vector of any object $\dgb_\a$, where $\a \in \itva$.  This vector does not depend on the choice of $\a$.

\begin{theorem}
\label{thm_sample_gpd_agrees}
Let  $\dgb: \R \to \ecata$ be a constructible functor in the sense of \cite{patel2018generalized}. If $\ecata$ is abelian and $\Sub(\dga_\a)$ has finite length for each $\a \in \R$, then  $\gpdb_\dga$ and $\spd_\dga$ are well defined, and
	\begin{align*}
	\gpdb_\dga(\itva) = \jhv{ \spd_\dga(\itva) }
	\end{align*}
for each $\itva \in \allitv{\tob}$.
\end{theorem}

Theorems \ref{thm_sample_cpd_agrees} and \ref{thm_sample_gpd_agrees} illustrate the overlap that exists in information captured by $\cpd, \gpda, \gpdb$, and $\spd$.  However, the saecular approach complements the preceding notions in several important ways:
\\

\noindent \tb{Inverse problems and cycle representatives.}  The classical persistence diagram, PD, plays a preeminent role in modern topological data analysis.  It is used when data takes the form of a filtered topological space $X_1 \su \cdots \su X_n$, from which one can derive a chain functor $\dga: [n] = \{1, \ldots, n\} \to \fielda \catvectf, \; p \mapsto \funh_m(X_p ; \fielda)$.  A \emph{basis of cycle representatives} for the functor $\dga$ consists of two parts: (i) an internal direct sum decomposition $\dga = \bigoplus_k \dga^k$, where each $\dga^k$ is a dimension-1 interval functor, and (ii) for each $k$, a cycle  $z^k \in \funz_m(X_n; \fielda)$ that generates that summand $\dga^k$.  

Cycle representatives have played a formative role in the development of homological persistence, because they allow one to map bars in the barcode back to data; concretely, one can localize the $k$th bar at the support of $z^k$, $\supp(z^k) \su X$.  However,  the notion of a persistent cycle representative is only defined for chain functors $\dga: \tob \to \fielda\catvectf$, at present, because chain functors valued in other categories generally cannot be expressed as direct sums of interval diagrams.

The saecular framework bridges this gap.  In particular, we can define a \emph{generating set} for interval $\itva$ to be a family of cycles $\seta \su \funz_m(X_n)$ that collectively generate the interval functor $\spd_\dga(\itva)$, in the natural sense.  Indeed, we can also define generators for barcodes in persistent \emph{homotopy}, \emph{cf.} \S\ref{sec_homotopy}.  This is one of several concrete benefits of functoriality.
\\

\noindent \tb{Vanishing of free components.}  
In generalized persistence, the torsion-free components of abelian groups either vanish (type-$\amitb$ diagrams) or require the introduction of \emph{sign conventions} (type-$\amita$ diagrams).  The saecular barcode naturally accommodates the torsion-free component of abelian groups without the need for sign conventions.
\\

\noindent \tb{Constraints on $\dga$ and $\ecata$.} The generalized persistence diagrams  $\gpdb$ are defined only for chain functors  $\dga: \R \to \ecata$ that meet certain criteria.  For example, $\ecata$ must be abelian and $\dga$ must be \emph{constructible}.   Constructible functors  resemble finite collections of constant functors, in the sense that every constructible functor $\dga$ admits a finite collection of disjoint intervals $\itva_1, \ldots, \itva_m$ such that $\R = \itva_1 \cup \cdots \cup \itva_m$ and $\dga(a \le b)$ is an isomorphism whenever there exists an interval $\itva_k$ such that  $a,b \in \itva_k$.
The saecular persistence diagram exists under substantially looser conditions.
\\

\noindent \tb{Extension problems.} As formalized in Theorem \ref{thm_sample_gpd_agrees}, the saecular PD solves a basic extension problem posed by the type-$\amitb$ generalized PD.  In particular, if $\ecata \in \{ R \catmod , \catmod R\}$, the functor $\dga$ is constructible, and each $\dga_\a$ has finite length as a module, then  $\gpdb_\dga(\itva)$ represents the family of composition factors of the nonzero objects $\spd_\dga(\itva)$.
\\

\noindent \tb{Functoriality and uniqueness.} The saecular PD sits within a lager categorical framework.  In particular, Theorem \ref{thm_main_result_finite_fun} implies that $\spd_\dga$ is essentially the restriction of the saecular functor $\sfun$ to the (full) subcategory of simple objects (i.e.\ singletons) in $\lclc ( \allitv{\tob})$.  This fact has diverse structural implications, many of which will be expanded below.  Moreover, while a given chain functor $\dga: \tob \to \ecata$ could admit many Krull-Schmidt decompositions in general, it admits at most \emph{one} saecular homomorphism, and, up to canonical isomorphism, only one saecular functor and saecular PD.
\\

There exists, moreover, a productive interaction between the formalisms of generalized and saecular persistence.  
The authors of the present work might never have considered persistence outside $\fielda \catvectf$ were it not for the introduction of generalized persistence modules by \cite{patel2018generalized} -- and had not the author of that work posed the problem of how to compute $\gpdb$, algorithmically, as an open problem.  Conversely, the machinery of saecular persistence  served as inspiration for important theoretical results concerning generalized persistence, e.g.\ \cite{mccleary2020edit}.

\subsection{Series}

The notion of a series has many useful realizations in abstract algebra -- subnormal series, central series, composition series, etc.  Saecular decomposition can also be repackaged as a series.

\begin{theorem} 
\label{thm_cseries_sample}
Suppose that $\tob$ is well ordered and $\ecata \in \{R\catmod, \catmod R, \catgroup \}$.  For each linearization $\lina$ of   $\allitv{\tob}$, there exists a unique $\forall$-complete lattice homomorphism $\cdflin: \axl{\lina} \to \Sub_\f$ such that
	\begin{align}
		\frac
		{ \cdflin \dsh_\lina   (\itva) }
		{ \cdflin \pdshp\lina (\itva) }		
	\in
	\allifun(\itva )
	\label{eq_intervalconstraints}
	\end{align}
for each $\itva \in \allitv{\tob}$.  To wit, 
	$$
	\cdflin = \sfun|_{\axl{\lina}}
	$$
where $\sfun$ is the saecular CDI homomorphism.
\end{theorem}

\begin{definition}
We call $\cdflin$ the \emph{subsaecular series} of $\dga$ subordinate to $\lina$.  The fraction 
	$
	\frac
		{ \cdflin \dsh_\lina   (\itva) }
		{ \cdflin \pdshp\lina (\itva) }	
	$
is the \emph{$\itva^{th}$ 	subsaecular factor} of $\lina$.  When $\ecata = \catgroup$, the $\itva^{th}$ subsaecular \emph{cokernel} factor is $\Cok \left ( \frac{\cdflin \dsh \itva}{\cdflin \pdsh \itva} \right )$.
\end{definition}

\emph{Mirabile dictu}, the associated factors are uniquely determined, up to isomorphism.

\begin{theorem}
\label{thm_sample_saecular_factors_unique}
The $\itva^{th}$ subsaecular factor of $\lina$ is canonically isomorphic to $\sfun (\dsh \itva)/\sfun(\pdsh \itva)$, for all $\itva$, where $\sfun$ is the saecular homomorphism.  In particular, the isomorphism type of the $\itva^{th}$ subsaecular factor is independent of $\lina$.
\end{theorem}

\begin{remark}
The canonical isomorphism in Theorem \ref{thm_sample_saecular_factors_unique}  is described concretely in \S\ref{sec_exact_fun_from_lattice_hom} and \S\ref{sec_coset_functors}.  The construction requires no special machinery, only a pair of Noether isomorphisms.
If $\ecata = \catgroup$, then the canonical isomorphism is an arrow in $[\tob,\catsetp]$; it is an  arrow in $[\tob, \catgroup]$ iff $(\cdflin \pdsh \itva)_\a \unlhd (\cdflin \dsh \itva)_\a$ for all $\a \in \tob$.
\end{remark}

If $\tob = [n] = \{1, \ldots, n\}$, equipped with the usual order, then we can represent the subsaecular series $\cdflin$ visually, via a  schematic of form \eqref{eq_subsaecular_absract_schematic}.  The construction proceeds as follows.  Poset $\lina$ is a sequence of intervals $([a_1, b_1), \ldots, [a_m, b_m))$.  For convenience, write $\cdflin_p$ for $\cdflin \dsh_\lina [a_p, b_p)$.  Let $\Zt = (\Zt_0, \Zt_1, \ldots)$ be the sequence obtained  from $(0 = \sss_0,  \ldots ,\sss_m = \dga)$ by deleting each $\sss_p$ such that  $\sss_p = \sss_{p-1}$.  Define $\zt_k = \Zt_k/\Zt_{k-1}$, and define $(p_\bullet, q_\bullet)$ such that $\Zt_k = \cdflin \dsh_\lina {[p_k, q_k)}$ for each $k$.  Thus each $\zt_k$ is a nonzero interval module of support type $[p_k,q_k)$, and the sequence $\zt = (\zt_1, \zt_2, \ldots)$ runs over all nonzero interval factors of $\dga$.  We call $\zt$ the \emph{reduced series} of $\apatha$.  Schematic  \eqref{eq_subsaecular_absract_schematic} illustrates this data.

\begin{equation}
\begin{tikzpicture}

\draw[gray, thick] (1,4) -- (1.75, 3.25);
\draw[gray, thick] (2.25, 2.75) -- (4,1);

\foreach \i in {2, 3}
    {
     \pgfmathtruncatemacro{\x}{( 4 - \i )};
    \draw[gray, thick] (1+\i, 4-\i) -- (\i, 4-\i-1);
    \filldraw[black] (1+\i, 4-\i) circle (2pt) node[anchor=south west] {$ \Zt_{\x} $};	
    \filldraw[black] (1+\i, 4-\i) circle (2pt) node[anchor=east] {$ [p_{\x}, q_{\x} ) \:\: $};	
    \filldraw[black] (\i, 4-\i-1) circle (2pt) node[anchor=north east] {$  \zt_{\x}  $} ;	
    }

\draw[gray, thick] (1, 4) -- (0, 3);
\filldraw[black] (1, 4) circle (2pt) node[anchor=south west] {$ \Zt_{\ell} $};	
\filldraw[black] (1, 4) circle (2pt) node[anchor=east] {$ [p_{\ell}, q_{\ell} ) \:\: $};	
\filldraw[black] (0, 3) circle (2pt) node[anchor=north east] {$  \zt_{\ell}  $} ;	    

\foreach \j in {0.85, 1, 1.15}
{
\filldraw[black] (1+\j, 4-\j) circle (1pt) node {};	
}

\end{tikzpicture}
\label{eq_subsaecular_absract_schematic}
\end{equation}

\begin{example}[Cyclic groups]
\label{ex:subsaecular_cyclic_new}

Consider the diagram of cyclic groups
\begin{equation}
    \begin{tikzcd}
		C_9  \arrow[r, "\times 2"] &      	
		C_6  \arrow[r, "\times 2"] &  
    	C_4
	\end{tikzcd}
    \label{eq_noetheriso}
\end{equation}
where $C_p$ is formally regarded as $\Z / p \Z$ and each map sends the coset $x + a\Z$ to the coset  $2x + b\Z$. Let us choose a linear order $\lina$ such that $[a,b) \le [c,d)$ when either $a < c$ or $a = c$ and $b \le d$.    The corresponding subsaecular series can then be expressed as follows.  For convenience, we write $C^x_y$ for the order-$y$ subgroup of  $C_x$.

\begin{equation}
\label{eq_compseriesexample}
\begin{tikzpicture}
\draw[gray, thick] (1,4) -- (4,1);
\foreach \i in {0,...,3}
{
 \pgfmathtruncatemacro{\x}{( 4 - \i )};
\draw[gray, thick] (1+\i, 4-\i) -- (\i, 4-\i-1);
}

\filldraw[black] (1, 4) circle (2pt) node[anchor=east] {$[3, 4 ) \;\; $} ;	
\filldraw[black] (2, 3) circle (2pt) node[anchor=east] {$[2,4) \;\; $} ;	
\filldraw[black] (3, 2) circle (2pt) node[anchor=east] {$[1,3) \;\; $} ;	
\filldraw[black] (4, 1) circle (2pt) node[anchor=east] {$[1,2) \;\; $} ;	

\filldraw[black] (1, 4) circle (2pt) node[anchor=south west] {$  C_9 \to C_6   \to C_4    $} ;	
\filldraw[black] (2, 3) circle (2pt) node[anchor=south west] {$  C_9 \to C_6   \to C^4_2  $} ;	
\filldraw[black] (3, 2) circle (2pt) node[anchor=south west] {$  C_9 \to C^6_3 \to 0  $} ;	
\filldraw[black] (4, 1) circle (2pt) node[anchor=south west] {$  C^9_3 \to 0 \to 0  $} ;

\filldraw[black] (0, 3) circle (2pt) node[anchor=north east] {$  0 \to 0   \to  \frac{C_4}{C^4_2}    $} ;	
\filldraw[black] (1, 2) circle (2pt) node[anchor=north east] {$  0 \to \frac{C_6}{C^6_3}   \to C^4_2  $} ;	
\filldraw[black] (2, 1) circle (2pt) node[anchor=north east] {$  \frac{C_9 }{C_3^9 }\to C^6_3 \to 0  $} ;	
\filldraw[black] (3, 0) circle (2pt) node[anchor=north east] {$	 C_3^9 \to 0 \to 0  $} ;	

\end{tikzpicture}
\end{equation}

\noindent Observe, in particular, that each factor $\zt_k$ is a type-$[p_k, q_k)$ interval functor.
\end{example}

\subsection{Cumulative density functions}

How can one calculate subsaecular series, in practice?  Paradoxically, the only reliable strategy we know is to solve the (seemingly harder) problem of evaluating  $\sfun$.  The work required to do so is greatly reduced by  Theorem \ref{thm_build_omega}.  To state this result, for each $\itva \in \allitv{\tob}$ and accompanying nested pair $(\seta, \setb) \sim \itva$, define 
	\begin{itemize*}
	\item $\Ks(\seta)$ to be the maximum-with-respect-to-inclusion subdiagram of $\dga$ such that $\Ks(\seta)_\a = 0$ for $\a \notin \seta$
	\item $\KKs(\seta)$ to be the minimum-with-respect-to-inclusion subdiagram of $\dga$ such that $\KKs(\seta)_\a = \dga_\a$ for $\a \in \seta$	
	\item $\At_{\seta, \setb}$ to be the intersection (i.e. the meet) $\KKs(\seta) \wedge \Ks(\setb)$
	\end{itemize*}
We call the function $\At: \allitv{\tob} \to \Sub_\dga, \; \itva \mapsto \At_\itva$ the \emph{saecular joint cumulative distribution function} of $\dga$.

\begin{theorem}
\label{thm_build_omega}
If $\sfun$ is the saecular homomorphism, then for each $\setc \in \axlh \allitv{\tob}$
	$$
	\sfun \setc = \bigvee_{ (\seta, \setb) \sim \itva \in \setc} \At_{\seta, \setb}
	$$
\end{theorem}

\begin{remark}
\label{rmk_generalization_of_interval_join_formula}
In later sections we will define the saecular CDI functor much more generally than it appears here. Theorem \ref{thm_build_omega} holds for this broader class of homomorphisms, as well.
\end{remark}

\begin{example}[Cyclic groups, continued]
\label{ex_cyclic_continued}
Let us compute the reduced series of the diagram $\dga$ shown in Example \ref{ex:subsaecular_cyclic_new}.  For convenience, let us take advantage of existing notation by writing $\Zt_1 \le \cdots \le \Zt_4$ for the reduced series shown in \eqref{eq_compseriesexample}.  To reduce notation further, we will write $\At_{pq}$ for $\At_{[1,p),[1,q)}$, where $\At$ is  the saecular joint cumulative distribution function. If we set $\Bt_{pq} = \frac{\At_{pq}}{\At_{p,q-1} \vee \At_{p-1,q}} = \frac{\sfun \; \dsh \; [p,q)}{\sfun \pdsh [p,q)}$, then $\Bt_{p_k q_k} \cong \zt_k$ for each $k$, by Theorem \ref{thm_sample_saecular_factors_unique}.  One may thus verify directly that 
\begin{align*}
	A
	=
    \begin{array}{r|ccccc|}
    \cline{2-6}  
     4 &  0 & \Zt_2 & \Zt_3 & \Zt_4 & \Zt_4 \\
     3 &  0 & \Zt_2 & \Zt_2 & \Zt_2 & \Zt_2\\ 
     2 &  0 & \Zt_1 & \Zt_1 & \Zt_1 & \Zt_1\\  
     1 &  0 & 0 & 0  & 0 & 0 \\  
     0 & 0 & 0  & 0 & 0 & 0\\ \cline{2-6}     
    \multicolumn{1}{c}{}& \multicolumn{1}{c}{0} & \multicolumn{1}{c}{1} & \multicolumn{1}{c}{2}& \multicolumn{1}{c}{3} & \multicolumn{1}{c}{4} \\     
    \end{array}     
	&&
	B
	=
    \begin{array}{r|ccccc|}
    \cline{2-6}  
     4 &   &  & \zt_3 & \zt_4 &  \\
     3 		&   & \zt_2 &  &  & \\ 
     2      &   & \zt_1 &  &  & \\  
     1 		&   &  &   &  &  \\  
     0 &  &  &  &  & \\ \cline{2-6}     
    \multicolumn{1}{c}{}& \multicolumn{1}{c}{0} & \multicolumn{1}{c}{1} & \multicolumn{1}{c}{2}& \multicolumn{1}{c}{3} & \multicolumn{1}{c}{4} \\     
    \end{array}     	    
\end{align*}
The reduced series is easily obtained from $\At$, via Theorem \ref{thm_build_omega}, and the reduced factors are obtained directly from $\Bt$.
\end{example}

\subsection{Homotopy}
\label{sec_homotopy}

Let $X$ be a filtered topological space, i.e. a functor $\tob \to \cattop$ such that $X(i \le j)$ is  inclusion for $i \le j$.  Composing with the homotopy functor $\pi_n$  yields a chain diagram $\dga = \pi_n \circ X$ valued in $\catgroup$ (for $n>0$; to avoid pathology, assume that $\tob$ has a minimum element and $X_i$ is path connected for all $i$). 

It is no surprise that the saecular barcode of $\pi_n \circ X$ captures information that the homological barcode of $X$ does not.  Indeed, if $\tob$ is a singleton, this reduces to the observation that $\funh_n(X)$ and $\pi_n(X)$ capture different information.  Example \ref{ex_saecular_homotopy_discriminates} gives a more interesting illustration.

\begin{example}
\label{ex_saecular_homotopy_discriminates}
Filtered spaces with isomorphic persistent homology, pointwise-isomorphic homotopy, and distinct saecular homotopy factors.

Figure \ref{fig_homotopy} presents a pair of filtered CW complexes, $X$ and $Y$.  The space $X_3$ has one vertex and three edges, denoted $a, b, c$, respectively. Filtrations $X$ and $Y$ agree up to index 3, i.e.  $X_i = Y_i$  for $i \le 3$.   Spaces $X_4$ and $Y_4$ each obtain by attaching a 2-cell to $X_3$; in each case, the resulting space is isomorphic to a wedge sum of $S^1$ with $T^2$. 

Equation \eqref{eq_homotopy_cdf} shows the simple cumulative distribution functions of $\pi_1 \circ X$ and $\pi_1 \circ Y$, denoted $\At^X$ and $\At^Y$, respectively. As in Example \ref{ex_cyclic_continued}, we write $\At_{pq}$ for $\At_{[1,p),[1,q)}$.  These functions may be calculated directly from the definition.  The corresponding saecular cokernel factors may be obtained from $\At^X$ and $\At^Y$ by taking quotients.  For economy, we adopt the following  conventions:
	\begin{itemize*}
	\item $\nr{ r_1, \ldots, r_n }$ is the free group on symbols $r_1, \ldots, r_n$
	\item $\xi( q | r_1, \ldots, r_n)$ is the cokernel of $\nr{r_1, \ldots, r_n} \su \nr{ q,  r_1, \ldots, r_n}$ in $\catgroup$; it  is isomorphic to $\Z$
	\item $\gamma(p,q | r_1, \ldots, r_n)$ is the minimum normal subgroup of $\nr{r_1, \ldots, r_n}$ that contains the commutator $pq p^{-1}q^{-1}$
	\item $F = \gamma(a,c | a,b,c)$, and $G = \gamma(a,b | a,b)$.
	\end{itemize*}

\end{example}

				\begin{figure}[h]
                  \centering
                    \includegraphics[width=0.9\textwidth]{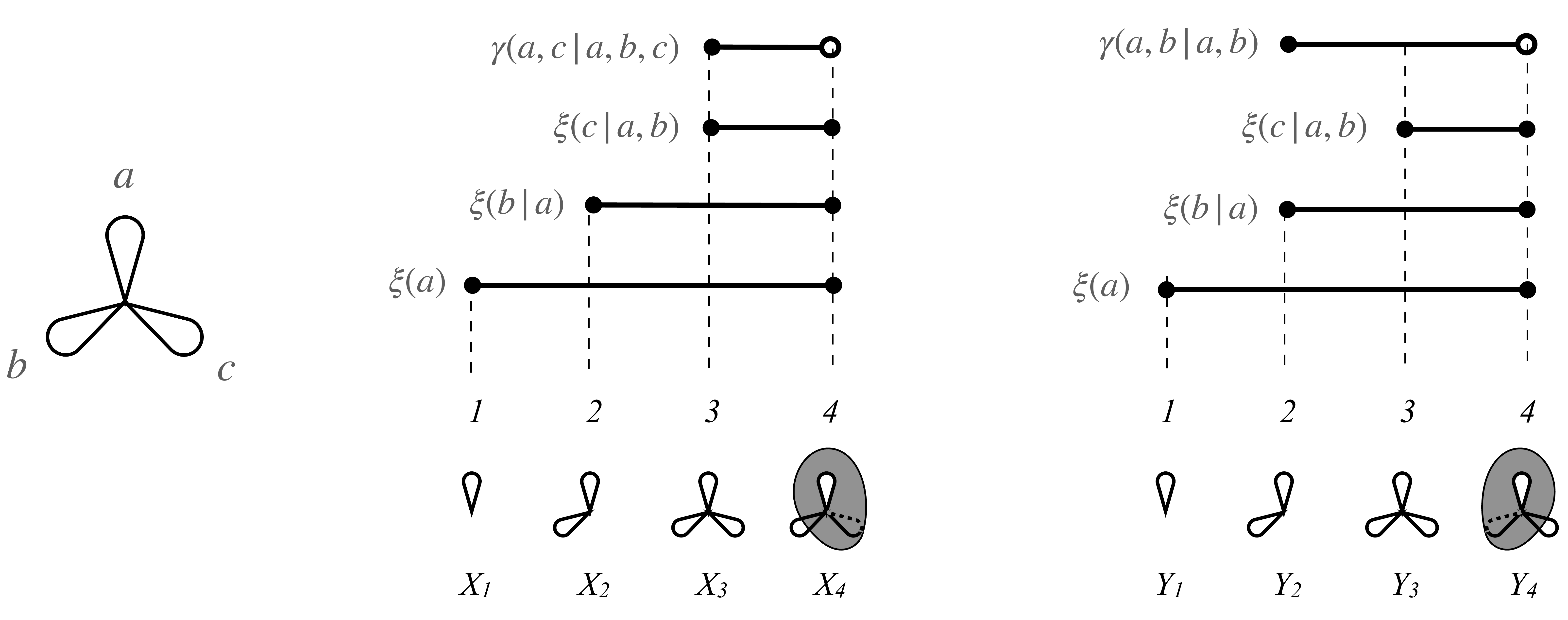}
                \caption{                	
					The saecular cokernel factors of two different filtered CW complexes, $X$ and $Y$.  Each nonzero factor appears as a closed or half-open interval.  
                }
                \label{fig_homotopy}
                \end{figure}

\begin{align}
	A^X
	&=
    \begin{array}{r|ccccc|}
    \cline{2-6}  
     5 &  0 & \nr{a} & \nr{a,b} & \nr{a,b,c} & \nr{a,b,c} \\
     4 &  0 & 0 & 0 & F & F \\
     3 &  0 & 0 & 0 & 0 & 0\\ 
     2 &  0 & 0 & 0 & 0 & 0\\  
     1 &  0 & 0 & 0 & 0 & 0 \\  
     0 &  0 & 0 & 0 & 0 & 0\\ \cline{2-6}     
    \multicolumn{1}{c}{}& \multicolumn{1}{c}{0} & \multicolumn{1}{c}{1} & \multicolumn{1}{c}{2}& \multicolumn{1}{c}{3} & \multicolumn{1}{c}{4} \\     
    \end{array}        
    &
	A^Y
	&=
    \begin{array}{r|ccccc|}
    \cline{2-6}  
     5 &  0 & \nr{a} & \nr{a,b} & \nr{a,b,c} & \nr{a,b,c} \\
     4 &  0 & 0 & G & G & G \\
     3 &  0 & 0 & 0 & 0 & 0\\ 
     2 &  0 & 0 & 0 & 0 & 0\\  
     1 &  0 & 0 & 0 & 0 & 0 \\  
     0 &  0 & 0 & 0 & 0 & 0\\ \cline{2-6}     
    \multicolumn{1}{c}{}& \multicolumn{1}{c}{0} & \multicolumn{1}{c}{1} & \multicolumn{1}{c}{2}& \multicolumn{1}{c}{3} & \multicolumn{1}{c}{4} \\     
    \end{array}      
    \label{eq_homotopy_cdf}  
\end{align}

\section{Organization}

In \S\ref{sec_literature} and \S\ref{chp:definitions} we review existing literature and introduce several notational conventions.  In \S\ref{sec_pec} and \S\ref{sec_lattic_free_hom} we recall the definitions of a Puppe exact category and a free lattice homomorphism, respectively; we also prove  structural properties needed for subsequent constructions.  In \S\ref{sec_saecular} we introduce saecular decomposition for chain diagrams valued in Puppe exact categories.    In \S\ref{sec_saecular_group} we adapt and extend this framework to diagrams valued in $\catgroup$.  In \S\ref{chp:hompersistence}, \S\ref{chp_pepersistence}, and \S\ref{sec_lsss} we relate saecular barcodes with (classical, $\fielda$-linear) persistence diagrams, generalized persistence diagrams, and Leray-Serre spectral sequences, respectively.  In \S\ref{sec_modified_birkhoff} we review a key technical result concerning free modular lattices.  In \S\ref{sec_wo_cd_extension} we prove a related result, which is used to establish existence of saecular CD homomorphisms in general, for functors from a well ordered set to $\ecata \in \{R \catmod, \catmod R, \catgroup\}$.  In \S\ref{sec_intro_proofs} we supply proofs for Theorems \ref{thm_sample_hom} - \ref{thm_build_omega}.

\section{Literature}
\label{sec_literature}

\emph{Functors valued in categories other than $\fielda\catvectf$.} Persistence modules valued in a category other than finite-dimensional $\fielda$-vector spaces find diverse interesting examples in data science.  The earliest examples of persistent topological structure in data include spaces with torsion, e.g.\ the  Klein bottle \cite{CIS+Local08}.  The method of circular coordinates developed in \cite{MSVPersistent11}, for example, relies explicitly on the use of integer coefficients.  Persistence for circle-valued maps has also been considered by Burghelea and Dey  in \cite{burghelea2013topological} and by Burghelea and Haller in \cite{burghelea2017topology}.

Mendez and Sanchez-Garcia have recently explored notions of persistence for homology with semiring coefficients, motivated by applications in biology, neuroscience, and network science \cite{mendez2021directed}. 

Work in Floer homology has recently prompted active exploration of the infinite-dimensional case.  In \cite{usher2016persistent}, Usher and Zhang introduces persistence for Floer homology via a non-Archimedean singular value decomposition of the boundary operator of the chain complex.  This work has generated a rapidly growing body of literature \cite{polterovich2017persistence}, including works on autonomous Hamiltonian flows \cite{polterovich2016autonomous} and rational curves of smooth surfaces \cite{browning2017rational}.

The decomposition of integer homological persistence modules (i.e., modules obtained by applying the $n$th homology functor with $\Z$ coefficients to a filtered sequence of topological spaces) has been studied by Romero et al.\ in \cite{romero2014defining}.  This treatment is closely entwined with spectral sequences, and has led to interesting directions in multi-persistence and spectral systems \cite{romero2018computational}.  This approach also has several interesting algorithmic aspects.  Filtrations and  bifiltrations play a prominent role.
\\

\noindent \emph{Generalized persistence.} The generalized persistence diagram was introduced for $\R$-parametrized constructible modules valued in Krull-Schmidt categories by Patel  in  \cite{patel2018generalized}.  Patel and McClearly have recently adapted the principles introduced in this work to achieve novel stability results in multiparameter persistence \cite{mccleary2020edit}.
\\

\noindent \emph{Chain functors, lattices, and exact categories.} The theory of Puppe exact categories has an expansive treatment in \cite{grandis12} and other works by M.\ Grandis. Of particular relevance, \cite{grandis12} concretely describes the universal RE model and p-exact classifying category for an $\tob$-indexed diagram in $E$, where $\tob$ is any totally ordered set and $\ecata$ is any exact category \cite{grandis12}. One of only a few technical challenges to relating this work with interval decomposition of persistence modules can be summarized as extending this construction from free lattice homomorphisms to free \emph{complete} lattice homomorphisms; the present work expand and addresses this challenge in detail.
The Puppe exact category of partial matchings has also appeared in the literature of persistence, specifically with regard to stability \emph{cf.}\ \cite{bauer2020persistence}.  
\\

\noindent \emph{Order lattices.}  The lattice-theoretic approach to persistence presented in this work was previewed by an analogous matroid theoretic treatment in \cite{henselman2016matroid} and \cite{henselman2017}. Lattice theoretic structure in persistence modules has also been explored in \cite{costa2014variable} and \cite{costaskraba13}. 
More broadly, the description of persistence in terms of the fundamental subspaces (kernel and image) has been a recurring theme since the emergence of the field.  The basics of this analysis appear in the seminal works of Robins \cite{RobinsComputational02}, Edelsbrunner, Letscher, and Zomorodian \cite{ELZTopological02}; as well as in Carlsson and Zomorodian \cite{ZCComputing05} which reformulates the persistence diagram in terms of graded modules. Bifiltrations play a central role in the reformulation of Carlsson and de Silva \cite{CSZigzag10} and in diverse other works.
\\

\noindent \emph{Persistent homotopy.}  The notion of a persistent homotopy group was introduced by Letscher \cite{letscher2012persistent}, as a tool to study knotted complexes which become unknotted under inclusion into larger spaces.  The definition of a persistent homotopy group, in this work, has a natural analog to that of a persistent homology group.  The family of all such groups corresponds conceptually to what Patel \cite{patel2018generalized} terms the \emph{rank function}, rather than a persistence diagram or barcode.  Letscher's treatment is formally disjoint from that of Patel, for several reasons; first, it includes homotopy in dimension 1, which does not fit into an abelian category; second, persistent homotopy groups are subquotients of loop spaces in which class representatives are well-defined -- this level of details is abstracted away in generalized PD's; third and finally, Letscher does not define a persistence diagram or barcode to complement the notion of a persistent homotopy group.  More recently,  M\'{e}moli and Zhou  \cite{memoli2019persistent} have explored several possible notions of persistent homotopy for metric spaces.      

In \cite{blumberg1705universality}, Blumberg and Lesnick show that interleavings form a universal pseudometric on persistent homotopy groups (unlike \cite{letscher2012persistent}, these authors use the term persistent homotopy group to refer to any functor obtained by composing a functor $\R \to \cattop$ with a homotopy functor $\pi_n$).  Several useful and fundamental properties of persistent homotopy are explored in \cite{batan2019persistent}.  Further properties of persistent homotopy have been explored in \cite{jardine2020persistent}.
\\

\noindent \emph{Spectral sequences}  Connections between persistence and spectral sequences have been suggested essentially since the formalization of persistent homology as a concept \cite{EZComputing01, EHComputational10}.  A  formal relation between the two was introduced in \cite{basu2017spectral}.  Revisions and extensions of this relationship to integral homology are discussed in \cite{romero2014defining}, though we do not know a refereed reference for this discussion as of yet.

\section{Definitions}
\label{chp:definitions}

We assume working familiarity with the definitions of poset, order lattice, and lattice homomorphism.  Definitions for complete, completely distributive (CD), modular, and algebraic lattices can be found in standard references.  As certain terms in order theory have distinct but related definitions elsewhere, however, there are several points with real risk of confusion.  We endeavor to indicate these cases as they arise, and, where necessary, to disambiguate  by introducing nonstandard variants of the standard terms.\\

\noindent
\tb{Partially ordered sets.} For economy of notation, $\pseta$ will denote both a poset and the associated  posetal category. The dual poset $\pseta^*$ represents the dual poset and the  opposite category $\pseta\op$.  The  partial order on $\pseta$ is denoted $\rel{\pseta}$, and the set of covering relations is  $\cov{\pseta} \su \rel{\pseta}$.  
Given $\lela \in \pseta$, we write
    	\begin{align*}
		\downarrow_\pseta \lela 
		&= 
		\{\lelb \in \pseta: \lelb \le \lela\} & \mathring \downarrow_\pseta \lela 
		= 
		\{\lelb \in \pseta: \lelb \lneq \lela \}
		\\
		\uparrow_\pseta \lela 
		&= 
		\{\lelb \in \pseta: \lelb \ge \lela\}  & \mathring \uparrow_{\pseta} \lela 
		= 
		\{\lelb \in \pseta: \lelb \gneq \lela \}
		\end{align*}
We omit subscripts on arrows where context leaves no room for confusion.

The maximum of $\pseta$, should such an element exist, is denoted  $1_\pseta$ or  $1$.  The minimum element is  denoted $0_\pseta$ or $0$.  A poset with both minimum and maximum is \emph{bounded}. A poset homomorphism $\mora: \pseta \to \psetb$ preserves \emph{existing} maxima   if $\pseta = \downarrow \psela$ implies $\psetb = \downarrow \mora_\psela$.  Preservation of existing minima (respectively, existing bounds) is defined similarly.  A bound-preserving poset homomorphism with a bounded domain and codomain is called \emph{bounded}.

Unless otherwise indicated, $\pseta \times \psetb$ denotes the canonical realization of the \emph{product} of $\pseta$ and $\psetb$ in the category of posets and order-preserving maps.  Concretely, if $\groundset(\pseta)$ denotes the ground set of $\pseta$, then 
\begin{align*}
\groundset(\pseta \times \psetb) = \groundset(\pseta) \times \groundset(\psetb)
\end{align*}
 and  $(\lela', \lelb') \le (\lela, \lelb)$ iff $\lela' \le \lela$ and $\lelb' \le \lelb$.    The $\inta$-fold product of $\pseta$ with itself is $\pseta^\inta$.
An \emph{interval} in $\pseta$ is a subset $\itva \su \pseta$ such that $\lela, \lelb \in \itva$ and $\lela \le x \le \lelb$ implies $x \in \itva$.  The family of intervals is denoted  $\allitv \pseta $.  Given $\lela, \lelb \in \pseta$, we write $[\lela, \lelb]$ for the interval $\{ \psela \in \pseta : \lela \le \psela \le \lelb\}$.

We will often wish to work with both a poset $\pseta$ and a copy of $\pseta$, whose elements are disjoint from, but nevertheless comparable to, the elements of $\pseta$.  We achieve this by defining $\dop \pseta = \pseta \times \{ 0\}$, the product of $\pseta$ with a 1-element poset.  We have $\pseta \cap \dop \pseta = \emptyset$, but there exists a canonical isomorphism $\pi: \dop \pseta \to \pseta, \; (\psela, 0) \mapsto \psela$.  Moreover, we can compare elements of $\dop \pseta$ to those of $\pseta$ by defining a binary relation (not a partial order) $\dople$ such that $\dop \psela \dople \psela$ iff $\pi(\dop \psela) \le \psela$, for each $\dop \psela \in \dop \pseta$ and $\psela \in \pseta$.
\\

\noindent\tb{Relations.} Given sets $\seta$, $\setb$ and a binary relation $\sim$ contained in $\seta \times \setb$, we write
	\begin{align*}
	\seta_{\sim \elb} := \{ \ela \in \seta : \ela \sim \elb \}
	\end{align*}
for each $\elb \in \setb$.  Similarly, given any function $\mora: \setc \to \seta$, we write
		$
		\mora_{\sim \elb}:   = \{ \elc \in \setc : \mora(\elc) \sim \elb\}.
		$
so that, for example, 
	\begin{align*}
	\N = \uparrow_\Z 0 =  \Z_{\ge 0}.
	\end{align*}

\noindent \tb{Complete lattices and homomorphisms.} 
A set function $\mora: \lata \to \latb$ \emph{preserves existing (nonempty) suprema} if $\bigvee \mora(\seta)$ exists and satisfies $ \bigvee \mora(\seta) = \mora(\bigvee \seta)$ for each (nonempty)  $\seta \su \lata$ such that $\bigvee \seta$ exists.  Preservation of existing (nonempty) infima is defined similarly.

A lattice $\lata$ is \emph{complete}\footnote{ \emph{Caveat lector}, the term complete assumes distinct meanings in the context of lattices generally and that of totally ordered sets specifically.  That of lattices follows the definition above, and is used exclusively throughout this work.  As an added safeguard, we will sometimes use the term \emph{totally ordered lattice} in place of \emph{totally ordered set}.}  if $\bigvee \seta$ and $\bigwedge \seta$ exist for each $\seta \su \lata$.   A lattice homomorphism $\mora$ is 
	\begin{itemize}
	\item \emph{$\existsc$-complete} if $\mora$ preserves existing suprema and infima
	\item \emph{$\existscc$-complete} if $\mora$ preserves existing nonempty suprema and infima
	\item \emph{$\foralc$-complete} if $\mora$ preserves existing suprema and infima, and has a complete domain
	\item \emph{$\foralcc$-complete} if $\mora$ preserves existing nonempty suprema and infima, and has a complete domain	
	\item \emph{complete } if $\mora$ preserves existing suprema and infima, and has a complete domain and codomain
	\end{itemize}
A sublattice $\latb \su \lata$ is  \emph{$\existsc$-complete}  if the inclusion $\latb \su \lata$ is $\existsc$-complete.  The same convention applies to all other notions of completeness listed above. \\

\noindent \tb{Upper and lower continuity.} 
 A lattice is \emph{upper continuous} or \emph{meet continuous} if for any $\lela \in \lata$ the poset homomorphism $x \mapsto x \wedge \lela$ preserves upward directed suprema \cite{maeda2012theory}.  It is \emph{lower continuous} or \emph{join continuous} if $\lata^*$ is meet continuous; that is, if  for any $\lela \in \lata$ the poset homomorphism $x \mapsto x \vee \lela$ preserves downward directed infima.
 \\

\noindent \tb{Set rings.} 
It will be convenient to assign a fixed notation to several set rings associated to a partially ordered set $\pseta$.  In particular, we define the following.
	\begin{description}
	\item[$\psl{\pseta}$] the power set lattice on the ground set of $\pseta$		
	\item[$\axl{\pseta}$] the lattice of decreasing subsets of $\pseta$. Equivalently, the complete set ring of {Alexandrov closed} subsets of $\pseta$.
	\item[$\bxl{\pseta}$] the lattice of proper nonempty decreasing subsets of $\pseta$
	\item[$\stl{\pseta}$] the minimum sublattice of $\axl{\pseta}$ containing $\{\dsh \psela : \psela \in \pseta\}$.   Unlike the preceding examples, this object is in general only a \emph{semitopology}.\footnote{By definition, a semitopology on a set $\seta$ is a bounded sublattice of $\psl{\seta}$.}
	\end{description}
These operations may be composed.  For example, the complete lattice $\axlh^2(\pseta)=\axlh(\axlh(\pseta))$ of  decreasing sets of decreasing sets will play a central role in the following story. \\

\noindent \tb{Free embeddings.} Given $\xxlh \in \{ \axlh, \bxlh\}$,   the \emph{free embedding}  $ \fem: \pseta \to \axlh \xxlh(\pseta)$ is the map such that
	\begin{align*}
	\fem(\psela)= \{ \setb \in \xxl{\pseta} : \psela \notin \setb\}
	\end{align*} 
for each $\psela \in \pseta$.   The term \emph{free} is justified in this context by Theorem \ref{thm_tunconstruction} (Tunnicliffe \cite{tunnicliffe1985free}). Some useful facts about free embeddings follow directly  from the definition: 
	\begin{enumerate*}
	\item If $\xxlh = \axlh$, then $\fem(\psela) = \axl{\pseta - \uparrow \psela}$.
	\item The free embedding fails to preserve suprema and infima, in general.\footnote{Indeed, failure to preserve limits its essential to the role of the free embedding.}
	\item The free embedding preserves existing maxima (respectively, minima) iff $\xxlh = \bxlh$.	
	\end{enumerate*}

\noindent \tb{Order chains.}  The term \emph{chain} is synonymous with \emph{totally ordered set} in classical order theory.\footnote{It has come to take a more restrictive meaning in other branches of mathematics, namely that of a totally ordered \emph{subset}.} A \emph{chain in $X$} means a totally ordered subset of $X$.  The notions of  boundedness and completeness which we have introduced for posets and lattices are also suitable for every chain, when regarded as an order lattice.  Thus a \emph{$\mathring \forall$-complete chain in $\lata$} is a  $\mathring \forall$-complete sublattice of $\lata$ that is totally ordered as a poset, a \emph{bounded chain} in $\pseta$ is a totally ordered subset containing maximum and minimum, etc. \\

\noindent \tb{Chain functors.}
  A \emph{chain functor} in a category $\ecata$ is a functor of form $\dga: \tob \to \ecata$ for some chain $\tob$.
When $\ecata$ has a zero object, the \emph{support} of $\mora$ can be defined
	\begin{align*}
	\supp(\dga): = \{ \tela \in \tob : \dga_\tela \neq 0\}.
	\end{align*}  
A functor has \emph{support type $\seta$} if $\supp(\dga) \in \{\seta, \emptyset\}$.  In particular,  the 0 functor has every support type.
	
A chain functor is \emph{interval} if 
	\begin{align*}
	\mora(\tela \le \telb) \quad \tm{is} \quad
		\begin{cases}
		\tm{invertible} & \tela, \telb \in \supp(\mora) \\
		\tm{zero} & else.
		\end{cases}
	\end{align*}
The support of every interval chain functor is an interval of $\tob$.  The \emph{object type} of a nonzero interval functor is the isomorphism class $\dga_\tela$, where $\tela$ is any element of $\supp(\dga)$.  The object type of the zero functor is the isomorphism class of the zero object. \\

\noindent \tb{Categories.} We write $\catcd$ for the category of completely distributive lattices and complete lattice homomorphisms, and $\catps$ for the category of partially ordered sets and order-preserving functions.  The category of sets and functions is $\catset$, and the category of pointed sets is $\catsetp$.  The category of groups and group homomorphisms is $\catgroup$.  The category of vector spaces and finite-dimensional vector spaces over $\fielda$ are denoted $\fielda\catvect$ and $\fielda\catvectf$, respectively.

\section{P-exact categories}
\label{sec_pec}

Puppe exact categories are generalizations of abelian categories that may or may not be additive.  Formally, a well-powered category $\cata$ is  \emph{Puppe exact} (or p-exact) if (i) it has a zero object, kernels, and cokernels, (ii) every mono is a kernel and every epi is a cokernel, (iii) every morphism has an epi-mono factorization.  Here we review the definitions and prove some basic structural results used in later sections. The reader is referred to \cite{grandis12} for further details.

\subsection{First principles}

Let  $\ecata$ be a p-exact category.  The direct and inverse image operators of an arrow $\mora: \oba \to \obb$ are denoted
	\begin{align*}
	\mora\di: \Sub(\oba) \to \Sub(\obb)
	&&
	\mora\ii: \Sub(\obb) \to \Sub(\oba)
	\end{align*}
respectively.

\begin{proposition}  
[Grandis, {\cite[p. 48-52]{grandis12}}]
\label{prop_pexact_modular_principles}
Let $\mora: \oba \to \obb$ be an arbitrary arrow in $\ecata$.  
	\begin{enumerate*}
	\item The poset $\Sub(\oba)$ is a modular lattice.
	\item The pair $(\mora\di, \mora\ii)$ is a \emph{modular connection}.  Concretely, this means that $\mora\di$ and $\mora\ii$ preserve order, and
		\begin{align*}
		\mora\ii \mora\di (\lela) = \lela \vee \mora\ii (0)
		&&
		\mora\di \mora\ii( \lelb) = \lelb \wedge \mora\di(1)
		\end{align*}
	for all $\lela \in \Sub(\oba)$ and all $\lelb \in \Sub(\obb)$.  In particular, $(\mora\di, \mora\ii)$ is a Galois connection.
	\item Consequently, if $\seta \su \Sub(\oba)$ and $\bigvee \seta$ exists, then $\bigvee(\mora\di \seta)$ exists, and $\mora\di(\bigvee \seta) = \bigvee \mora\di(\seta)$.  Dually $\mora\ii(\bigwedge \seta) = \bigwedge \mora\ii(\seta)$ for any $\seta$ such that $\bigwedge \seta$ exists.	
	\end{enumerate*}
\end{proposition}

The family of bounded modular lattices and modular connections forms a p-exact category under the composition rule
	$
	(\morb\di, \morb\ii) \com (\mora\di, \mora\ii) : = (\morb\di \mora\di, \mora\ii \morb\ii).
	$
We denote this category $\Mlc$, and denote the full subcategory of \emph{complete} modular lattices by $\cmc$.  If one defines $\Sub(\mora)$ to be the direct/inverse image pair $(\mora\di, \mora\ii)$ for each arrow $\mora$ in a Puppe exact category $\ecata$, then $\Sub$ is an (exact)  functor  $\ecata \to \Mlc$.  For details see \cite{grandis12}.

Each arrow $\mora: \oba \to \obb$ in $\Mlc$ admits a  unique epi-mono factorization of form
\begin{center}
    \begin{tikzcd}
    	\downarrow  \mora\ii 0  \arrow[r,"m", hook]   &\oba \arrow[rr,"\f"]  \arrow[d,two heads, "q"] && \obb \arrow[r,two heads, "p"] & \uparrow \mora\di 1 \\
	& \uparrow \mora\ii 0 \arrow[rr,"g"] && \downarrow \mora\di 1 \arrow[u,hook, "n"]
    \end{tikzcd}
    \label{eq_signedextensioncom}
\end{center}   
where
	\begin{align*}
	q\di(x) & = x \vee \mora\ii0  && n\di(y) = y && g\di(x) = \mora\di(x) \\
	q\ii(x) & = x  && n\ii(y) = y \wedge \mora\di 1 && g\ii(y) = \mora\ii(y)
	\end{align*}
Each subobject of $\lata$, regarded as subobject in either  $\Mlc$ of $\cmc$, engenders a unique mono of form
	\begin{align*}
	(\mora\di, \mora\ii): [0, \lela]_\lata \to \lata
	&&
	\mora\di \lelb = \lelb
	&&
	\mora\ii \lelb = \lela \wedge \lelb.
	\end{align*}
Thus  subobjects of $\lata$ are in canonical 1-1 correspondence with elements of $\lata$, and, equivalently, with intervals of form $[0, \lela]_\lata$.  A similar correspondence holds between quotients of $\lata$, connections of form $(\mora\di, \mora\ii) : \lata \to [\lela, 1]$, where $\mora\ii$ is inclusion and $\mora\di(\lelb) = \lela \vee \lelb$, elements of $\lata$, and intervals of form $[\lela, 1]_\lata$.

\begin{lemma}
\label{lem_direct_image_complete_on_CD_complete_sublattice} 
The direct image operator of an arrow $\f: \oba \to \obb$ restricts to a complete lattice homomorphism on each $\forall$-complete CD sublattice of $\Sub(\oba)$ containing  the kernel of $\f$.  The dual statement holds for inverse image operator. Similarly, the direct image operator of an arrow $\f: \oba \to \obb$ restricts to a lattice homomorphism on each distributive sublattice of $\Sub(\oba)$ containing  the kernel of $\f$.  The dual statement holds for inverse image operator.
\end{lemma}
\begin{proof}
Let $\lata$ be a $\forall$-complete CD sublattice of $\Sub(\oba)$ that contains  $k := \ker(\f)$, and  fix $\seta \su \lata$.  Complete distributivity then provides the third of the following identities
	\begin{multline*}
	\textstyle
	\f \di (\bigwedge \seta) 
	= 
	\f \di(\f \ii \f \di \bigwedge \seta)
	=
	\f \di (( \bigwedge \seta) \vee k)
	=
	\f \di ( \bigwedge \{ \ela \vee k : \ela \in \seta \})
	=
	\\
	\textstyle	
	\f \di ( \bigwedge \{ \f \ii \f \di \ela : \ela \in \seta \})
	= 
	\f \di( \f \ii \bigwedge \{ \f \di \ela : \ela \in \seta \})
	=
	\bigwedge \f \di (\seta).
	\end{multline*}
This establishes the claim for the direct image operator, since $\f \di$ preserves existing joins.  The claim for inverse image is equivalent, by duality.  If we assume, instead, that $\lata$ is distributive, then similar argument shows the direct and inverse image operators preserve binary meets and joins.
\end{proof}

Given any index category $\icata$, we denote the category of  $\icata$-shaped diagrams in $\ecata$ by	$[\icata, \ecata]$. Diagrams make new p-exact categories from old.

\begin{lemma}[P-exact diagrams]   
\label{lem_subdiagrams}
Let $\icata$ be a small category.
	\begin{enumerate*}
	\item The category of diagrams $[\icata, \ecata]$ is p-exact. 
	\item Kernels and cokernels obtain object-wise.  Concretely,  $0 \rightarrow \oba \xrightarrow{m} \obb \xrightarrow{p} \obc \rightarrow 0$ is exact in $[\icata, \ecata]$ if and only if $0 \rightarrow \oba_\tela \xrightarrow{m_\tela} \obb_\tela \xrightarrow{p_\tela} \obc_\tela \rightarrow 0$ is exact in $\ecata$, for each $\tela$.
	\item An arrow $\oba \xrightarrow{m} \obb$ in $[\icata, \ecata]$ is mono iff $m_\tela$ is mono for each  $\tela$.  Likewise for epis.
	\item If $(\obb^\telb)_{\telb \in J}$ is an  indexed family  of subdiagrams of  $\oba$ and $\bigvee_\telb  \obb^\telb_\tela$ exists for each  $\tela$, then  $\bigvee \obb$ exists  and obtains object-wise.  The dual statement holds for meets.  
	\end{enumerate*}
\end{lemma}
\begin{proof}  Claim 1 is a routine exercise, as observed in \cite[p.\ 48]{grandis12}.  Coincidence of exact and object-wise exact sequences follows from the universal property of (co)kernels. Claim 3 follows Claim 2, since p-exact categories have normal epis and monos \cite[p.\ 46]{grandis12}.  Claim 4 follows Claim 3.
\end{proof}

\subsection{Regular and canonical induction}
\label{sec_induced_maps}

If $\cata$ is any category with subobjects and $\oba$ is an object in $\cata$, then we may define a poset $\pseta(\oba)$ such that $\groundset(\pseta(\oba)) = \{ (\sob, \soa) \in \Sub(\oba)^2 : \sob \le \soa \le \oba\}$, ordered such that $(\sob, \soa) \le (\sob', \soa')$ iff $\soa \le \soa'$ and $\sob \le \sob'$.  We say that $(\sob, \soa) \le (\sob', \soa')$  is a \emph{doubly nested pair}, and call $\pseta(\oba)$ the poset of nested pairs in $\oba$.  We call the corresponding homomorphism $\soa/\sob \to \soa'/\sob'$ the \emph{regularly induced map}, \emph{cf.}\ \cite{grandis12}.  The following lemma can be proved directly using the axioms of kernels and cokernels.

\begin{lemma}
\label{lem_functor_exists_from_nested_pair_poset}
Regularly induced maps are closed under composition.  More precisely, for each object $\oba$ in $\ecata$ 
there exists a functor $\funca: \pseta( \oba) \to \ecata$ sending ${(\sob, \soa)}$ to $\soa / \sob$ and sending each relation $(\sob, \soa) \le (\sob', \soa')$ to the regularly induced map $\soa/\sob \to \soa'/\sob'$.
\end{lemma}

Following \cite{grandis12}, let us define a binary relation $\mdom$ on the family of subquotients of $\oba$ such that 
	\begin{align}
	\soa/ \sob \; \mdom  \; \soa' / \sob'
	&&
	\iff
	&&
	M \vee N' = N \vee M'
	\quad
	\text{and}
	\quad
	M \wedge N' = M' \wedge N.
	\label{eq_defofmutualdomination}
	\end{align}

\begin{remark}
Grandis defines this relation, \emph{mutual domination}, in a rather more general context than that which we require here.  Equivalence between  \eqref{eq_defofmutualdomination} and this more abstract formulation appears in expression (4.3.2) of  \cite{grandis12}, and a more accessible interpretation appears in the proof of Theorem 4.3.5.
\end{remark}

\begin{lemma}
\label{lem_phi_M_relation_for_semitopologies}
Suppose that $\sob \su \soa$ and $\sob' \su \soa'$ are elements of the lattice $\lata$ of closed sets in a semitopology.  Then 
	\begin{align*}
	\soa/ \sob \; \mdom  \; \soa' / \sob' \quad 
	\iff 
	\quad \soa - \sob = \soa' - \sob'
	\end{align*}
\end{lemma}
\begin{proof}
If $ \soa - \sob = \soa' - \sob'$ then $\soa \vee \soa' = \soa \vee \sob' = \sob \vee \soa'$ and $\sob \wedge \sob' = \sob \wedge \soa' = \sob' \wedge \soa$.  The converse holds by a similar set theoretic argument; see the note following Lemma 1.2.6 in \cite[p.\ 23]{grandis12}.
\end{proof}

\begin{theorem}[Induced isomorphisms \cite{grandis12}]  
\label{thm_mutualdomination}
If $\oba$ is an object in $\ecata$ and $\soa, \soa' \in \Sub( \oba)$, then the map $\phi$ defined by commutativity of  
\begin{equation}
    \begin{tikzcd}
    	\soa  \arrow[r, tail] & \soa \vee \soa'  \arrow[r, two heads] & (\soa \vee \soa') / \soa \\
		\soa \wedge \soa' \arrow[u, tail] \arrow[r, tail] & \soa' \arrow[u, tail] \arrow[r, two heads] & \soa' / (\soa \wedge \soa')  \arrow[u, "\phi" ', dashed]
    \end{tikzcd}
    \label{eq_noetheriso}
\end{equation} 
is an isomorphism.   Moreover, if $\soa / \sob \; \mdom \; \soa' / \sob'$ then there exists a commutative diagram of isomorphisms
\begin{equation}
    \begin{tikzcd}
    	& (\soa \vee \soa') / (\sob \vee \sob')    \\
		\soa/\sob   \arrow[ur] \arrow[rr, "\psi", dashed] && \soa'/\sob' \arrow[ul] \\
		& (\soa \wedge \soa') / (\sob \wedge \sob') \arrow[ul] \arrow[ur]
    \end{tikzcd}
    \label{eq_canonicallyinducediso}
\end{equation} 
where each solid arrow is derived as in \eqref{eq_noetheriso}.
\end{theorem}

\begin{proof}  
These results appear in \cite{grandis12}, Lemma 2.2.9 and Proposition 4.3.6.
\end{proof}

\begin{definition}
The maps $\phi$ and $\psi$ defined in \eqref{eq_noetheriso} and \eqref{eq_canonicallyinducediso} are \emph{Noether isomorphisms of the first and second kind}, respectively.  The latter may also be termed the \emph{canonical} or \emph{canonically induced} isomorphism.
\end{definition}

\begin{theorem}[Cohesion of canonically induced isomorphisms \cite{grandis12}]
\label{thm_grandis_distributive_coherence}
Let $\oba$ be an object in $\ecata$ and $\lata$ be a distributive sublattice of $\Sub(\oba)$.  Then the family of canonically induced isomorphisms between subquotients with numerator and denomenator in $\lata$ is closed under composition.
\end{theorem}
\begin{proof}
Grandis  \cite[pp. 23-24, Theorem 1.2.7]{grandis12} proves both the stated result and its converse in the special case where $\ecata$ is the category of abelian groups.  However, the same argument carries over verbatim for any p-exact category, if we use Theorem \ref{thm_mutualdomination} to ``stand in'' as a generalized version of the diagram \cite[p. 20, (1.16)]{grandis12} used in that argument.
\end{proof}

\newcommand{\soaR}{\hat \soa}
\newcommand{\sobR}{\hat \sob}
\newcommand{\soc}{Q}
\newcommand{\sod}{P}

\begin{theorem}[Compatibility of canonical isomorphisms with regular induction {\cite[p180, Theorem 4.3.7]{grandis12}} ]
\label{thm_regular_cohesion}
Let $(\sob, \soa) \le (\sob', \soa')$ and $(\sod, \soc) \le (\sod', \soc')$ be doubly nested pairs.  Let $\oba_{{\soa}/{\sob}}$ and $\oba_{{\soa'}/{\sob'}}$ be objects in $\ecata$ that realize the quotients $\soa/\sob$ and $\soa'/\sob'$, respectively, and define $\obb_{\soc/\sod}$ and $\obb_{\soc'/\sod'}$ similarly.

If $\soa/ \sob \; \mdom  \; \soc / \sod$ and $\soa'/ \sob' \; \mdom  \; \soc' / \sod'$, then the following diagram commutes, where horizontal arrows are regularly induced maps and vertical arrows are canonically induced isomorphisms.
\begin{equation}
    \begin{tikzcd}
    	\oba_{{\soa}/{\sob}} \ar[rr] 
		&& 
		\oba_{{\soa'}/{\sob'}} 
		\\
    	\obb_{{\soc}/{\sod}} \ar[rr, "" ]  \ar[u, dashed, "\cong"]
		&& 
		\obb_{{\soc'}/{\sod'}} \ar[u, dashed, "\cong"']
    \end{tikzcd}
    \label{eq_regular_cohesion}
\end{equation} 
\end{theorem}

\begin{remark}
\label{rmk_conflating_quotient_objects}
Theorem \ref{thm_regular_cohesion} intentionally avoids the conventional abuse of notation that conflates $\soa/\sob$ with a representative of a cokernel in $\ecata$.  Note, in particular, even when $(\sob, \soa) = (\sob', \soa')$ and $(\sod, \soc) = (\sod', \soc')$, it is possible that $\oba_{\soa/\sob} \neq \obb_{\soc/\sod}$.  Thus, in particular, Theorem \ref{thm_regular_cohesion} says something nontrivial about the relationship between regularly induced maps and the isomorphisms that exist (via the universal property of cokernels) between distinct representatives of a cokernel. 
\end{remark}

A \emph{quasi-regular induction} on a  sublattice $\lata \su \Sub_\oba$ is a commutative square of form
\begin{equation}
    \begin{tikzcd}
    	\frac{\soa}{\sob} \ar[rr, "\mora"] 
		&& 
		\frac{\soc}{\sod} 
		\\
    	\frac{\soa'}{\sob'} \ar[rr, "" ]  \ar[u, dashed, "\cong"]
		&& 
		\frac{\soc'}{\sod'} \ar[u, dashed, "\cong"']
    \end{tikzcd}
    \label{eq_define_quasiregularlyin_induced}
\end{equation} 
where the numerator and denominator of each subquotient lie in $\lata$, vertical arrows are canonically induced isomorphisms (in particular  $\soa/\sob \mdom \soa'/ \sob'$ and  $\soc/\sod \mdom \soc'/ \sod'$) and the lower horizontal arrow is regularly induced.  In this case we say that $\mora$ is $\lata$-quasi-regularly induced, and write $f : (\soa, \sob, \soc, \sod) \sim (\soa', \sob', \soc', \sod')$.


\begin{lemma}
\label{lem_quasiregularly_induced_maps_unique}
If  $\lata$ is a distributive sublattice of $\Sub(\oba)$ and $\soa \le \sob$ and $\soc \le \sod$ are nested pairs in $\lata$, then there exists at most one $\lata$-quasi-regularly induced arrow $\soa / \sob \to \soc/\sod$.
\end{lemma}
\begin{proof}
Suppose that $f: (\soa, \sob, \soc, \sod) \sim (\soa', \sob', \soc', \sod')$ and $f: (\soa, \sob, \soc, \sod) \sim (\soa'', \sob'', \soc'', \sod'')$. Then we have a diagram of regularly induced maps (solid) and canonically induced isomorphisms (dashed)
\begin{equation}
    \begin{tikzcd}
    	&
		\frac{\soa}{\sob} \ar[rr, dotted, " \mora"] 
		&& 
		\frac{\soc}{\sod} 
		\\
		&
    	\frac{\soa'}{\sob'} \ar[rr, "" ]  \ar[u, dashed, "\cong"]
		&& 
		\frac{\soc'}{\sod'} \ar[u, dashed, "\cong"']
		\\
    	\frac{\soa''}{\sob''} \ar[rr, "" ]  \ar[ur, dashed, "\cong"] \ar[uur, dashed, "\cong"]		
		&& 
		\frac{\soc''}{\sod''} \ar[ur, dashed, "\cong"']		\ar[uur, dashed, "\cong"]		
    \end{tikzcd}
    \label{eq_quasiregularlyin_induced_unique}
\end{equation}
Each triangle of dashed arrows in \eqref{eq_quasiregularlyin_induced_unique} commutes by distributive cohesion (Theorem \ref{thm_grandis_distributive_coherence}).  The bottom face, composed of two parallel solid arrows and two parallel dashed arrows, commutes by Theorem \ref{thm_regular_cohesion}.  Thus any dotted arrow $\mora$ that makes the back face commute also makes the slanted face commute, and vice versa. 
\end{proof}

\subsection{Exact functors via lattice homomorphisms}
\label{sec_exact_fun_from_lattice_hom}

Let $\lata$ be the lattice of closed sets in a semitopological space $\setx$.  Let $\cdf: \lata \to \Sub_\oba$ be a bound-preserving lattice homomorphism from $\lata$ to the subobject lattice of an object $\oba$ in a p-exact category $\ecata$.

The main result of this section is Theorem \ref{thm_lclcextension}.  This closely resembles (and may be subsumed by) several existing results in \cite{grandis12}. We give a proof which is completely elementary, and which can be modified to work outside the the setting of p-exact categories, which we will do in later sections.

\begin{theorem}
\label{thm_lclcextension}
There exists an exact functor
	\begin{align*}
	\sfun: \lclc(\lata) \to \ecata
	\end{align*}
which agrees with $\cdf$ on $\lata$, in the sense that 
	$
	\sfun|_{\lata} = \cdf.
	$
	
If $\sfun'$ is any other functor satisfying the same criterion, then there exists a unique natural isomorphism $\eta: \sfun \cong \sfun'$ such that $\eta|_{\lata}$ is identity.  Transformation $\eta$ assigns the canonically induced isomorphism $\eta_{\oba/\obb} : \sfun_{\oba/\obb} \cong \sfun'_{\oba/\obb}$ to each nested pair $\obb \le \oba \in \lata$.
\end{theorem}
\begin{proof}
We will construct $\sfun$ directly.  To reduce notation, write $|\seta|$ for $\cdf \seta$.
Invoking the axiom of choice if necessary, for each object $\oba \in \lclc(\lata)$ select a nested pair $\sob_\oba \le \soa_\oba \in \lata$ such that $\oba = \soa_\oba - \sob_\oba $, and define $\sfun \oba = \xfrac{\soa_\oba}{\sob_\oba}$.

For each regularly induced map $\mora: \obb \to \obb'$ in $\lclc(\lata)$, select a doubly nested pair $(\sob, \soa) \le (\sob', \soa')$ such that $\oba = \soa - \sob$ and $\oba' = \soa' - \sob'$. Let $\sfun \mora$ be the unique map such that the following diagram commutes, where vertical arrows are canonically induced isomorphisms; such  isomorphisms exist by Lemma \ref{lem_phi_M_relation_for_semitopologies} and Theorem \ref{thm_mutualdomination}.
\begin{equation}
    \begin{tikzcd}
    	\xfrac{\soa_\oba}{\sob_\oba} \ar[rr, "\sfun \mora"] 
		&& 
		\xfrac{\soa_{\oba'}}{\sob_{\oba'}} 
		\\
    	\xfrac{\soa}{\sob} \ar[rr  ]  \ar[u, dashed, "\cong"]
		&& 
		\xfrac{\soa'}{\sob'} \ar[u, dashed, "\cong"']
    \end{tikzcd}
    \label{eq_define_sfun_on_arrows}
\end{equation} 

\noindent \underline{Claim 1.} Every choice of of doubly nested pair $(\sob, \soa) \le (\sob', \soa')$ yields the same value for $\sfun \mora$.  
\emph{Proof.} Map $\sfun \mora$ is quasi-regularly induced on $\Im(\cdf)$; uniqueness follows from Lemma \ref{lem_quasiregularly_induced_maps_unique}. 
\\

\noindent \underline{Claim 2.} Function $\sfun$ preserves composition of regularly induced maps.  More precisely, if $f, g, h$ are regularly induced maps such that $h = gf$, then $\sfun h = (\sfun g) (\sfun f)$.    
\emph{Proof.} Select doubly nested pairs $(\sob, \soa) \le (\sob', \soa')$ and $(\sob', \soa') \le (\sob'', \soa'')$ for $f$ and $g$.  Then $(\sob, \soa) \le (\sob'', \soa'')$ is a doubly nested pair for $h$. Then we have a diagram of for \eqref{eq_quasiregularlyin_induced_unique}, with regularly induced arrows (solid, lefthand side), arrows determined by $\sfun$ (solid, righthand side) and canonically induced isomorphisms (dashed).

\begin{equation}
    \begin{tikzcd}
    	&
		\xfrac{\soa''}{\sob''} \ar[rr, dashed,""] 
		&& 
		 \sfun \frac{\soa''}{\sob''} 
		\\
		&
    	\xfrac{\soa'}{\sob'} \ar[rr, dashed,"" ]  \ar[u, ""]
		&& 
		 \sfun \frac{\soa'}{\sob'} \ar[u, "\sfun \morb"']
		\\
    	\xfrac{\soa}{\sob} \ar[rr, dashed, "" ]  \ar[ur, ""] \ar[uur, ""]		
		&& 
		 \sfun \frac{\soa}{\sob} \ar[ur, "\sfun \mora"']		\ar[uur, "\sfun \morc" near end]		
    \end{tikzcd}
    \label{eq_quasiregularlyin_induced_unique}
\end{equation}
The triangle on the lefthand side of \eqref{eq_quasiregularlyin_induced_unique} commutes by Lemma \ref{lem_functor_exists_from_nested_pair_poset}.  The three rectangular faces of  \eqref{eq_quasiregularlyin_induced_unique} by Claim 1.  Thus the entire diagram commutes.  The desired conclusion follows.
\\

We have now defined $\sfun$ on regularly induced maps, a family which includes all epis and monos.  To each map $h$ in $\lclc(\lata)$ corresponds a unique epi-mono decomposition $h = gf$, so we may define a new function $\sfun'$ on arbitrary maps via $\sfun' h = (\sfun' g)(\sfun' f)$.  In fact, $\sfun'$ agrees with $\sfun$ on regularly induced maps by Claim 2, so we may simply extend the domain of definition of $\sfun$ to include all morphisms, via $\sfun'$.
\\

\noindent \underline{Claim 3.} Function $\sfun$ is a functor.  
\emph{Proof.} It suffices to show that $\sfun$ preserves composition of morphisms.  Fix arrows
	\begin{align*}
	x  & = ba : M / N \to P / Q \\
	y  & = dc : P / Q \to S / T
	\end{align*}
where $a$, $c$ are epi and $b$, $d$ are mono.  Diagram \eqref{eq_sfunpreservescomp_new} is then uniquely determined by commutativity  and the condition that $fe$ be an epi-mono decomposition.

\begin{equation} 
    \begin{tikzcd}
    M/N 	\arrow[r, two heads, "a"] 	\arrow[ddr, two heads, "g"']
    &
    P'/ Q 	\arrow[r, tail, "b"]	\arrow[dd, two heads, "e"] \arrow[ddr]
    &
    P/Q 	\arrow[dd, two heads, "c"]    
    \\
    \\
    &
    (P' \vee R)/R 	\arrow[ddr, tail, "h"'] 	\arrow[r, tail, "f"]
    &
    P/ R 	\arrow[dd, tail, "d"]    
    \\
    \\
    &&
    S/T 
    \end{tikzcd}
    \label{eq_sfunpreservescomp_new}
\end{equation} 

\noindent Each arrow in diagram \eqref{eq_sfunpreservescomp_new} is regularly induced, so the diagram commutes by Claim 2.  It follows that $\sfun (yx) = (\sfun y) (\sfun x)$, since $hg$ is the epi-mono decomposition of $yx$.
\\

\noindent \underline{Claim 4.}  Functor $\sfun$ is exact. That is, $\sfun$ preserves short exact sequences.  
\emph{Proof.}  Fix a short exact sequence $0 \to \oba \xrightarrow{f} \obb \xrightarrow{g} \obc \to 0$ in $\lclc(\lata)$.  Then there exist nested pairs $(\sob, \soa') \le (\sob, \soa) \le (\soa', \soa)$ such that the following diagram commutes, where vertical arrows are canonically induced isomorphisms
\begin{equation} 
    \begin{tikzcd}
    0 \ar[r]
    &
	\xfrac{\soa'}{\sob} \ar[r]
	&
	\xfrac{\soa}{\sob} \ar[r]	
	&
	\xfrac{\soa}{\soa'} \ar[r]
	&
	0
	\\ 	
	0 \ar[r]
	&
	\sfun \oba \ar[u, dashed] \ar[r, "f"]
	&
	\sfun \obb \ar[u, dashed] \ar[r, "g"]
	&
	\sfun \obc \ar[u, dashed] \ar[r]
	&
	0	
    \end{tikzcd}
    \label{eq_sfun_exactness}
\end{equation}

This concludes the proof that a suitable exact functor $\sfun$ exists.  Now posit a second functor $\sfun'$ with the same properties. We must show that the map $\eta$ defined in the theorem statement is a natural isomorphism $\sfun \cong \sfun'$, and that $\eta = \eta'$ for any natural isomorphism $\eta'$  such that $\eta'_\seta = 1$ for each $\seta \in \lata$.  

To verify that $\eta$ is a natural transformation, fix a doubly nested pair $(\sob, \soa) \le (\sob', \soa')$. Since $\sfun( \soa / \sob) \mdom |\soa|/|\sob|$ and $\sfun'( \soa / \sob) \mdom |\soa|/|\sob|$, and because $\sfun( \soa' / \sob') \mdom |\soa'|/|\sob'|$ and $\sfun'( \soa' / \sob') \mdom |\soa'|/|\sob'|$, Theorem \ref{thm_regular_cohesion} can be invoked to argue that the following diagram commutes, where $\mora$ is the regularly induced map $\soa/\sob \to \soa'/ \sob'$ and vertical arrows are canonically induced isomorphisms (the proof is similar to that of Lemma \ref{lem_quasiregularly_induced_maps_unique}).
\begin{equation}
    \begin{tikzcd}
    	\sfun(\soa /\sob ) \ar[rr, "\sfun \mora"]
		&& 
		\sfun(\soa' /\sob' )
		\\
    	\sfun'(\soa /\sob ) \ar[rr, "\sfun' \mora"]  \ar[u, dashed, "\eta_{\soa /\sob}"]
		&& 
		\sfun'(\soa' /\sob' )  \ar[u, dashed, "\eta_{\soa' /\sob'}"']
    \end{tikzcd}
    \label{eq_sfun_nat_iso_exists}
\end{equation} 

\noindent This proves that $\eta$ is a natural isomorphism, since every map in $\lclc(\lata)$ is a composition of regularly induced morphisms.

For uniqueness, posit a natural isomorphism $\eta': \sfun \to \sfun'$ such that $\eta'_\seta$ is identity for each $\seta' \in \lata$.  Then  for each object $\oba = N/D$ in $\lclc(\lata)$ we have a commutative diagram

\begin{equation} 
    \begin{tikzcd}
    \sfun D \arrow[r, tail] 
    & \sfun N  \arrow[r, two heads] 
    & \sfun(N/D)
    \\
    \sfun' D \arrow[r, tail] \arrow[u, equal]
    & \sfun' N  \arrow[r, two heads] \arrow[u, equal]
    & \sfun'(N/D) \arrow[u, "\cong", "\hk'_{\oba}" ']
    \end{tikzcd}
    \label{eq_inducednaturaliso}
\end{equation} 

\noindent By the uniqueness, $\hk'_{\oba}$ coincides with the canonically induced isomorphism.  Thus $\eta' = \eta$, as desired.
\end{proof}

\begin{remark}
One can apply Theorem \ref{thm_lclcextension} to any bounded distributive lattice $\lata$ by identifying $\lata$ with a semitopological space via Stone duality, which states that every bounded distributive lattice is isomorphic to the lattice of compact open topology of a spectral space, or via Priestley duality, which states that each bounded distributive lattice is isomorphic to the lattice of clopen up-sets in a Priestley space.
\end{remark}

\begin{remark}
\label{rmk_induced_exact_fuctor_for_groups}
Theorem \ref{thm_lclcextension} carries over unchanged if $\ecata$ is not a p-exact category but $\catgroup$, provided that $\cdf$ factors through the lattice $\latb$ of normal subgroups of $\Sub_\oba$.  Indeed, the proofs of Theorems  \ref{thm_mutualdomination} - \ref{thm_regular_cohesion} and Lemma \ref{lem_quasiregularly_induced_maps_unique} carry over unchanged to the setting of groups, so long as we assume that all subobjects listed in the theorem statements are normal; the same holds for the proof of Theorem \ref{thm_lclcextension}.
\end{remark}

\begin{definition}  
We all $\sfun$ the \emph{locally closed functor} of $\cdf$.  

For economy, we often abuse notation by writing $\cdf$ for $\sfun$.  This leaves little room for confusion in practice, since $\cdf$ and $\sfun$ agree on $\lata \su \lclc(\lata)$.
\end{definition}


\section{Lattice homomorphisms}
\label{sec_lattic_free_hom}

Free lattices and free homomorphisms play a fundamental role in our story.  Here we review the essential details.

\subsection{Free maps}

The following two classical results from the theory of free lattices are paramount.

\begin{theorem}[Tunnicliffe 1985, \cite{tunnicliffe1985free}] 
\label{thm_tunconstruction}
  If $ \fem$ is the free embedding $\pseta \to \axlh^2(\pseta)$, then the following are equivalent.
	\begin{enumerate*}
	\item There exists an $\forall$-complete, CD sublattice of $\lata$ containing the image of $\psma$.
	\item There exists an $\forall$-complete lattice homomorphism $\morc$ such that $\mora = \morc \fem$.  
	\end{enumerate*}	
\begin{equation} 
            \begin{tikzcd}
        				&& {\axlh^2(\pseta)} \arrow[dll , "\fem{}" ' , <-] \arrow[drr, "\morc", dashed]	&& \\
        			\pseta \arrow[rrrr,  "\mora"] &&&& \cdla
        	\end{tikzcd}
    \label{eq_freediagram}
\end{equation}  
When it exists, the commuting homomorphism $\morc$ satisfies 
	\begin{align}
	\morc(\seta )
	=
	\bigvee _{X \in M_\seta}   \bigwedge  f (X)
	=
	\bigwedge_{Y \in N_\seta} \bigvee f(Y)
	\label{eq_cdfformula}
	\end{align}
where 
	$
	M_\seta = \{X \su \pseta : \bigwedge \fem(X) \su \seta   \}
	$
and
	$
	N_\seta = \{X \su \pseta : \seta \su \bigvee \fem(X) \}.
	$
In particular, $\morc$ is unique.
\end{theorem}
\begin{proof}
If $\lata$ is itself CD then the desired conclusion follows verbatim from \cite{tunnicliffe1985free}.  Otherwise one obtains a commuting homomorphism $h^\latb:{\axlh \xxlh(\pseta)} \to \latb $ for each $\mathsf{X}$-complete, CD sublattice $\latb \su \lata$.  It is straightforward to show that $h^\latb = h^{\latb'}$ for any two $\mathsf{X}$-complete CD sublattices $\latb$ and $\latb'$, since both must satisfy \eqref{eq_cdfformula}.
\end{proof}

\begin{remark}
\label{rmk_free_cd_formula_alternate_form}
The set $\bigwedge \fem(X)$ contains a maximum element $\xi(X) := \pseta - \bigcup_{x \in X}(\uparrow x)$, since $\fem(x) = \axl{\pseta - \uparrow x }$.  Therefore $M_\seta = \{ X \su \pseta : \xi(X) \in \seta\}$.
\end{remark}

We call  $\axlh^2(\pseta)$  the \emph{free completely distributive (CD) lattice on $\pseta$}.  The corresponding map $\morc$ is the \emph{free completely distributive homomorphism (FCDH) of $\mora$},  denoted   $\FCD(\mora)$.

\begin{theorem}[Birkhoff] 
\label{thm_birkhoff}
Let $\lata$ be a bounded lattice and $\dop \tob, \tob$ be totally ordered sets.  Let $\pseta$ be the coproduct in $\catps$ of $\axlh^2(\dop \tob )$ and $\axlh^2( \tob)$, denoted  $\axlh^2(\dop \tob ) \sqcup  \axlh^2( \tob)$.  Define $\cpph: \axlh^2(\tob) \to \axlh^2(\dop \tob \sqcup \tob)$ by $\cpph(\seta) = \{ \dop \setb \sqcup \setb : \dop \setb \in \axlh(\dop \tob), \; \setb \in \seta \}$.  Define $\dop \cpph: \axlh^2(\dop \tob) \to \axlh^2(\dop \tob \sqcup \tob)$ similarly.  Let $\dxlh \axlh( \dop \tob \sqcup \tob)$ be the sublattice of $\axlh^2( \dop \tob \sqcup \tob)$ generated by $\Im(\dop \cpph) \cup \Im(\cpph)$.  Finally, let $\morb: \pseta \to \lata$ be any order-preserving map that preserves the bounds of $\axlh^2(\dop \tob)$ and $\axlh^2(\tob)$.  
Then the following are equivalent.
	\begin{enumerate*}
	\item There exists a modular sublattice of $\lata$ containing the image of $\morb$.
	\item There exists lattice homomorphism $\morc$ (which happens to preserve bounds) such that $\morb = \morc \circ (\dop \cpph \sqcup \cpph)$.  
	\end{enumerate*}	
\begin{equation} 
            \begin{tikzcd}
        				&& \dxlh \axlh (\pseta) \arrow[dll , "\dop \cpph \sqcup \cpph" ' , <-] \arrow[drr, "\morc", dashed]	&& \\
        			\pseta \arrow[rrrr,  "\morb"] &&&& \cdla
        	\end{tikzcd}
    \label{eq_freediagram_modular_specialized}
\end{equation}  
When a commuting homomorphism $\morc$ exists, it is unique.
\end{theorem}
\begin{proof}
This is a special case of a classic result of Birkhoff \cite[p66, III, Theorem 9]{birkhoff1973lattice}.  See  \S\ref{sec_modified_birkhoff} for further details.
\end{proof}

We call  $\dxlh \axlh(\pseta)$ the \emph{free bounded distributive (BD) lattice on $\pseta$}.  The corresponding map $\morc$ is the \emph{free bounded distributive homomorphism (FBDH) of $\morb$},  denoted $\FBD(\morb)$.  We will be particularly interested in the case where $\morb = \FCD(\mora|_{\dop \tob}) \sqcup \FCD(\mora|_{\tob})$ for some order-preserving map $\mora: \dop \tob \sqcup \tob \to \lata$.  In this case we write $\ov \mora$ for $\morb$, hence $\FBD({\bar \mora}) = \FBD(\morb)$.

\begin{remark}
For economy of notation, we will sometimes use $\psma$ to denote either $\FBD(\psma)$ or $\FCD(\psma)$. 
\end{remark}

\subsection{Free interval maps}

Let $\tob$ be a totally ordered set and $\dop \tob = \tob \times \{0 \}$ be a disjoint copy of $\tob$ equipped with the canonical isomorphism $\dopiso: \dop \tob \to \tob, \; (\a, 0) \mapsto \a$.  Define $\dop \cpph$ and $\cpph$ such that the following diagram commutes, where $\dop \fem, \fem$, and $\nu$ are free embeddings:
\begin{equation} 
    \begin{tikzcd}
	\axlh^2( \dop \tob) \ar[rr, "\dop \cpph"]
	&&
	\axlh^2( \dop \tob \sqcup \tob )
	&&
	\axlh^2( \tob)  \ar[ll, "\cpph" ']
	\\ 	
	\dop \tob \ar[u, "\dop \fem"] \ar[rr]
	&&
	\dop \tob \sqcup \tob \ar[u, "\nu"] 
	&&
	\tob \ar[u, "\fem" '] \ar[ll]
    \end{tikzcd}
    \label{eq_free_embedding_system}
\end{equation} 
Concretely, $\cpph(\seta) = \{ \dop \setb \sqcup \setb : \dop \setb \in \axl{ \dop \tob}, \; \setb \in \seta \}$.  Map $\dop \cpph(\dop \seta)$ follows a similar rule.

Let 	
	$$
	\alljtvh( \dop \tob \sqcup \tob) = \{ \dop \seta \sqcup \seta \in \axl{ \dop \tob \sqcup \tob} : \dopiso( \dop \seta ) \subsetneq \seta \} 
	$$
and define
	\begin{align*}
	\iid \di : \axlh \alljtvh( \dop \tob \sqcup \tob) \leftrightarrows  \axlh \axlh( \dop \tob \sqcup \tob): \iid \ii
	&&
	\iid \di( \seta) = \dsh_{\axlh( \dop \tob \sqcup \tob)} \seta 
	&&
	\iid \ii( \seta) = \seta \cap \alljtvh( \dop \tob \sqcup \tob)
	\end{align*}
Maps $\iid \ii$ and $\iid \di$ are polarities, and the composition $\iclosure$ is a closure operator.  

\begin{remark}[Equivalence of $\alljtvh( \dop \tob \sqcup \tob)$ and $\allitvh(\tob)$]
\label{rmk_equivalence_of_j_and_i}
There is a canonical bijection $\phi: \alljtv{ \dop \tob \sqcup \tob } \to \allitv{ \tob }, \; \dop \seta \sqcup \seta \mapsto \seta - \dopiso(\seta)$.  We endow $\allitvh(\tob)$ with the partial order such that $\phi$ is an isomorphism of posets; this agrees with the partial order on $\allitv{\tob}$ described in the introduction.  
We call $\phi(\dop \setb \sqcup \setb)$ the interval of $\dop \setb \sqcup \setb$, and $\dop \setb \sqcup \setb = \phi^{-1}(\itva)$ the split partition of $\itva$.

  In consideration of this isomorphism, we will abuse notation by passing elements of $\allitv{\tob}$ as arguments to functions which take $\alljtv{\dop \tob \sqcup \tob}$ as inputs.
\end{remark}

As in Theorem \ref{thm_birkhoff}, let $\dxlh \axlh( \dop \tob \sqcup \tob)$ denote the (bounded, incidentally) sublattice of $\axlh^2( \dop \tob \sqcup \tob)$ generated by the union of images $\Im(\dop \cpph) \cup \Im(\cpph)$.  Similarly, let $\dxlh \alljtvh( \dop \tob \sqcup \tob)$ denote the (bounded) sublattice of $\axlh \alljtvh ( \dop \tob \sqcup \tob)$ generated by the union  $\Im( \iid \ii \dop \cpph) \cup \Im( \iid \ii \cpph)$.  

Let $\mora: \dop \tob \sqcup \tob \to \lata$ be any order-preserving function.  Suppose that $\mora(\dop \tob)$ and $\mora(\tob)$ include into (possibly distinct) $\forall$-complete sublattices of $\lata$, and let $\dgb = \FCD( \mora|_{ \dop \tob}) \sqcup  \FCD( \mora|_{ \tob})$ be the copairing in $\catps$ of the free CD homomorphisms generated by  $\mora|_{\dop \tob}$ and $\mora|_{\tob}$.   Then  for each $\xxlh \in \{ \axlh , \dxlh\}$ we obtain a commuting diagram of solid arrows \eqref{eq_free_diagram_xdi}.  We call this the \emph{bifiltered mapping diagram} of $\mora$.

\begin{equation} 
    \begin{tikzcd}
	\xxlh \alljtvh (\dop \tob \sqcup \tob) 
		\ar[rrrr, dashed, "\morc'"]
	&&&&
	\lata
	\\    
	& 
	\xxlh \axlh (\dop \tob \sqcup \tob) 
		\ar[urrr, dotted, "\morc"]
		\ar[ul, "\iid \ii"]
	&&
	\axlh^2(\dop \tob) \sqcup \axlh^2(\tob)
		\ar[ll, "\dop \cpph \sqcup \cpph" ']
		\ar[ur, near start , "\bar \mora"]		
	\\
	&& 
	\dop \tob \sqcup \tob 
		\ar[ul,  "\nu"  ]
		\ar[uurr, bend right,  "\mora" ']
		\ar[ur, "\dop \fem \sqcup \fem"  ]	
    \end{tikzcd}
    \label{eq_free_diagram_xdi}
\end{equation} 

A free  BD \emph{interval} homomorphism (FBDIH) of $\mora$ is a bound-preserving homomorphism $\morc'$ such inserting $\morc'$ into the solid arrow diagram \eqref{eq_free_diagram_xdi} produces a new commuting diagram, where $\xxlh = \dxlh$.  Similarly, a free  CD \emph{interval} homomorphism (FCDIH) is a map $\morc'$ such that inserting $\morc'$ into the solid arrow diagram \eqref{eq_free_diagram_xdi} produces a new commuting diagram, where $\xxlh = \axlh$.  These maps are unique, if the exist, by Theorems \ref{thm_free_distributive_mapping} and \ref{thm_free_cd_mapping}.   We may therefore speak of \emph{the}  free BDIH and free CDIH,  denoted $\IFBD(\bar \mora) $ and  $\IFCD(\mora)$, respectively.

\begin{theorem}
\label{thm_free_distributive_mapping}
Let $\mora$ be an order-preserving function with bifiltered mapping diagram \eqref{eq_free_diagram_xdi}.
	\begin{enumerate*}
	\item If $\bar \mora$ admits a free BDH (respectively, BDIH), then this homomorphism is unique.	
	\item If $\bar \mora$ admits a free BDIH 
	then it admits a free BDH, 
	 and 
		\begin{align}
		\IFBD(\bar \mora) = \FBD(\bar \mora) \iid \di
		&&
		\FBD(\bar \mora) = \IFBD(\bar \mora) \iid \ii				
		\end{align}

	\item Map $\bar \mora$ admits a free BDH iff $\Im(\bar \mora) \su \latb$ for some modular sublattice $\latb \su \lata$.
	\item Map $\bar \mora$ admits a free BDIH iff it admits a free BDH and $\bigwedge \mora( \setb) \le \bigvee \mora( \dop \setb )$ for every split partition $(\dop \setb, \setb)$.
	\end{enumerate*}
\end{theorem}

\begin{theorem}
\label{thm_free_cd_mapping}
Let $\mora$ be an order-preserving function with bifiltered mapping diagram \eqref{eq_free_diagram_xdi}.
	\begin{enumerate*}
	\item If $\mora$ admits a free CDH (respectively, CDIH), then this homomorphism is unique.	
	\item If $\mora$ admits a free CDIH 
	then it admits a free CDH, 
	 and 
	 \begin{align}
		\IFCD(\mora) = \FCD(\mora) \iid \di
		&&
		\FCD(\mora) = \IFCD(\mora) \iid \ii				
		\end{align}

	\item Map $\mora$ admits a free CDH iff $\Im(\morb) \su \latb$ for some $\forall$-complete CD sublattice $\latb \su \lata$.
	\item Map $\mora$ admits a free CDIH iff it admits a free CDH and $\bigwedge \mora( \setb) \le \bigvee \mora( \dop \setb )$ for every split partition  $(\dop \setb, \setb)$.
	\end{enumerate*}
\end{theorem}

\begin{example}
\label{ex_unit_interval_map}
Let $\tob = \lata$ be the real unit interval $[0,1] \su \R$. We will use  Theorem \ref{thm_free_cd_mapping} to show that the map
	\begin{align*}
	\textstyle
	\Psi: \axlh{\allitv{\tob}}\to \tob 
	\quad
	\quad
	\quad
	\quad
	\seta \mapsto \bigvee_{\itva \in \seta} \bigwedge \itva
	\end{align*}
is a complete lattice homomorphism.  Later, we will use this result to construct a pathological counterexample (specifically, Counterexample \ref{cx_wo_needed_for_uniqueness}) in \S\ref{sec_uniqueness}.  Moreover, we claim that 
	\begin{align}
	\Psi \dsh \itva = \Psi \pdsh \itva
	\label{eq_psi_dsh_equa_pdsh}
	\end{align}
for all $\itva \in \allitv{\tob}$.  This fact is most easily proved in terms of an example.  Suppose that $\itva = [a,b]$.  If $0<a$, then $\pdsh [a,b]$ contains every interval of form $[a - \epsilon, b]$ for $\epsilon \in (0,a)$, hence, by completeness $\Psi \pdsh [a,b] = a$.  On the other hand, if $\a = 0$ then $0 \le \Psi \pdsh [a,b] \le \Psi \pdsh [a,b] = 0$.  

To prove that $\Psi$ is complete, let $\mora: \dop \tob \sqcup \tob \to \lata$ be the poset map that sends $t$ and $(t,0)$ to $t$, for each $t \in \tob$.  Then $\mora$ admits a free CD homomorphism $\sfun$ and a  free CDI homomorphism $\sfun' : \alljtv{\dop \tob \sqcup \tob} \to \lata$, by Theorem \ref{thm_free_cd_mapping}.  Let $\xi(X) = (\dop \tob \sqcup \tob) - \bigcup_{x \in X} \uparrow x$,  as in Remark \ref{rmk_free_cd_formula_alternate_form}.  Then completeness of $\sfun'$ implies the first of the following identities, and Remark \ref{rmk_free_cd_formula_alternate_form} implies the third:
	\begin{align}
	\sfun' \seta
	=	
	\bigvee_{\dop \setb \sqcup \setb \in \seta } \sfun' \dsh_{\alljtvh(\dop \tob \sqcup \tob) }(\dop \setb \sqcup \setb)
	=		
	\bigvee_{\dop \setb \sqcup \setb \in \seta } \sfun \dsh_{\axlh(\dop \tob \sqcup \tob) }(\dop \setb \sqcup \setb)
	=		
	\bigvee_{\xi(X)\in \seta} 
	\bigwedge \mora(X)
	\label{eq_interval_formula}
	\end{align}
Every element  $\dop \setb \sqcup \setb \in \alljtv{\dop \tob \sqcup \tob}$ can be expressed in form $\xi(X)$, where $X = ( \dop \tob \sqcup \tob) - (\dop \setb \sqcup \setb)$ and $\inf \mora(X)  = \inf \mora(X \cap \dop \tob) = \inf \mora( \setb - \dopiso(\dop \setb)) = \inf(\setb - \dopiso(\dop \setb))$.  Thus $\Psi = \sfun' \circ \Phi$, where $\Phi: \axl{ \allitv{ \tob } } \to \axl{ \alljtv{ \dop \tob \sqcup \tob } }$ is the isomorphism induced by the isomorphism $\phi$ defined in Remark \ref{rmk_equivalence_of_j_and_i}.  In particular, $\Psi$ is the composition of two complete homomorphisms, hence a complete homomorphism, which was to be shown.
\end{example}

\subsection{Proof of Theorems \ref{thm_free_distributive_mapping} and \ref{thm_free_cd_mapping}}

\begin{proof}[Proof of Theorem \ref{thm_free_distributive_mapping}]

Assertion 3 follows from Theorem \ref{thm_birkhoff}, as does uniqueness of  $\FBD(\bar \mora)$.  If $\bar \mora$ admits a free BDIH $\morc$, then $\morc': = \morc  \iid \di$ is a free BDH.  In particular,  it is the \emph{unique} free BDH, $\FBD(\bar \mora)$.  Since $\iid \ii \iid \di = 1$, it follows that $\morc = \FBD(\bar \mora) \iid \di$.  This proves uniqueness of $\morc$, which we may now call \emph{the} free BDIH.  Assertions 1 and 2 follow.

To prove assertion 4, assume that $\IFBD(\bar \mora)$ exists.  Then $\FBD(\bar \mora)$ exists, and $\FBD(\bar \mora) =  \FBD(\bar \mora) \iclosure$ by assertion 2.  Thus $\bigwedge \mora( \setb) \le \bigvee \mora( \dop \setb )$ for every split partition $(\dop \setb, \setb)$, by Lemma \ref{lem_closure_comp_eq}.  Conversely, suppose that $\FBD(\bar \mora)$ exists, and that $\bigwedge \mora( \setb) \le \bigvee \mora( \dop \setb )$ for every split partition.  Then $\FBD(\bar \mora) =  \FBD(\bar \mora) \iclosure$, by Lemma \ref{lem_closure_comp_eq}.  Thus the homomorphism $\morc': = \FBD(\bar \mora) \iid \di$ satisfies $\morc'  \iid \ii \dop (\cpph \sqcup \cpph) = \FBD(\bar \mora) (\cpph \sqcup \cpph) = \bar \mora$.  Thus $\bar \mora$ admits a free BDIH, namely $\morc'$.
\end{proof}

\begin{proof}[Proof of Theorem \ref{thm_free_cd_mapping}]

Assertion 3 follows from Theorem \ref{thm_tunconstruction}, as does uniqueness of  $\FCD( \mora)$.  If $\mora$ admits a free CDIH $\morc$, then $\morc': = \morc  \iid \di$ is a free CDH.  In particular,  it is the \emph{unique} free CDH, $\FCD(\mora)$.  Since $\iid \ii \iid \di = 1$, it follows that $\morc = \FCD(\mora) \iid \di$.  This proves uniqueness of $\morc$, which we may now call \emph{the} free CDIH.  Assertions 1 and 2 follow.

To prove assertion 4, assume that $\IFCD(\mora)$ exists.  Then $\FCD(\mora)$ exists, and $\FCD(\mora) =  \FCD(\mora) \iclosure$ by assertion 2.  Thus $\bigwedge \mora( \setb) \le \bigvee \mora( \dop \setb )$ for every split partition $(\dop \setb, \setb)$, by Lemma \ref{lem_closure_comp_eq}.  Conversely, suppose that $\FCD(\mora)$ exists, and that $\bigwedge \mora( \setb) \le \bigvee \mora( \dop \setb )$ for every split partition.  Then $\FCD(\mora) =  \FCD(\mora) \iclosure$, by Lemma \ref{lem_closure_comp_eq}.  Thus the homomorphism $\morc': = \FCD(\mora) \iid \di$ satisfies $\morc'  \iid \ii \dop (\cpph \sqcup \cpph) = \FCD(\mora) (\cpph \sqcup \cpph) = \mora$.  Thus $\mora$ admits a free CDIH, namely $\morc'$.
\end{proof}

\newcommand{\predxi}{\xi}
\newcommand{\predxii}{\dop \xi}
Given a subset $\seta \su \tob$, define $\predxi(\seta) = \{ \setc \in \axlh(\tob) : \setc \subsetneq \setb \text{ for some } \setb \in \seta \} \in \axlh^2(\tob)$.  Define $\predxii(\dop \seta) \in \axlh^2(\dop \tob)$ similarly, for $\dop \seta \su \dop \tob$.

\begin{lemma}
\label{lem_deflation}
For each $\dop \seta \in \axlh^2(\dop \tob)$ and each $\seta \in \axlh^2(\tob)$ one has 
	\begin{align}
	\iclosure \cpph (\seta) 
	&= 
	\dop \cpph  \predxii  ( \dopiso \seta) \wedge \cpph(\seta)
	\label{eq_deflation_formula_a}	
	\\
	\iclosure \dop \cpph (\dop \seta) 
	&= 
		\begin{cases}
		\dop \cpph (\dop \seta)  & \seta < 1 \\
		\cpph \predxii (1) & \seta = 1
		\end{cases}
	\label{eq_deflation_formula_aa}
	\end{align}	
As a special case, for every split partition $(\dop \setb, \setb)$ one has
	\begin{align}
	\textstyle
	\iclosure (\bigwedge_{\a \in \setb} \nu_\a ) 
	&= 
	\textstyle
	(\bigvee_{\aa \in \dop \setb} \nu_{\aa} ) \wedge (\bigwedge_{\a \in \setb} \nu_\a )	
	\label{eq_deflation_formula_split_a}	
	\end{align}	
Moreover, for each $\seta \in \axlh^2(\tob)$, exactly one of the following statements holds true: (i) $\cpph (\seta)$ can expressed in form $\bigwedge \nu(\setb) = \bigwedge \{\nu _\a : \a \in \setb\}$ for some split partition  $(\dop \setb, \setb)$, and $\seta$ contains a maximum element $\setc = \max(\seta)$, (ii)   $\seta$ contains no maximum, and $\iclosure \cpph \seta = \dop \cpph (\dopisoo \seta) \wedge \cpph \seta$.   
\end{lemma}
\begin{proof}
One has
	\begin{align}
	(\iid \di \iid \ii) ( \cpph \seta)
	&=
	\{	
		\dop \setc \sqcup \setc \in \axlh^2( \dop \tob \sqcup \tob)
		:
		\dop \setc \sqcup \setc \su  \dop \setx \sqcup \setx
		\text{ for some }
		\dop \setx \sqcup \setx \in \cpph \seta
		\text{ such that }
		\dopiso (\dop \setx ) \subsetneq \setx
	\}
	\\
	&=
	\{	
		\dop \setc \sqcup \setc \in \axlh^2( \dop \tob \sqcup \tob)
		:
		\dop \setc \in \predxii( \dopisoo (\seta)), \; \setc \in \seta 
	\}	
	\\
	& = 
	\dop \cpph  \predxii (\dopisoo \seta) \wedge \cpph(\seta)
	\label{eq_deflation_formula_proof_a}
	\end{align}
This establishes \eqref{eq_deflation_formula_a}.  It is simple to check that $\iclosure \dop \cpph (\dop \seta) = \dop \cpph (\dop \seta)$ when $\dop \seta < 1$, and that $\iclosure \dop \cpph (1) = \dop \cpph \predxii(1)$.  Formula \eqref{eq_deflation_formula_aa} therefore follows from the observation that $\iclosure$ preserves binary meets, since $\dop \seta \wedge 1 = \dop \seta$.  Equation \eqref{eq_deflation_formula_split_a} follows from the case where $\seta = \bigwedge \fem(\setb)$.  
The remainder of the theorem is a routine exercise.
\end{proof}

\begin{lemma}
The closure operator $\iclosure$ is a join-complete lattice homomorphism, and $\iclosure(\dxlh \axlh(\dop \tob \sqcup \tob)) \su \dxlh \axlh(\dop \tob \sqcup \tob)$.  This operator preserves 0 but not 1.
\end{lemma}
\begin{proof}
It is simple to verify that $\iclosure$ preserves binary meets and arbitrary joins.  It maps $\dxlh \axlh(\dop \tob \sqcup \tob)$ into $\dxlh \axlh(\dop \tob \sqcup \tob)$ by \eqref{eq_deflation_formula_a} and \eqref{eq_deflation_formula_aa}.  The set $\iclosure 1$ does not contain $\dop \tob \sqcup \tob$; in particular, $\iclosure 1 \neq 1$.
\end{proof}

\begin{lemma}
\label{lem_closure_comp_eq}
If $\mora$ admits a free BDH (respectively, a free CDH) $\morc$, then $\morc = \iclosure \morc$ iff $\bigwedge \mora( \setb) \le \bigvee \mora( \dop \setb )$ for every split partition $(\dop \setb, \setb)$.  
\end{lemma}
\begin{proof}
Let  $\morc$ be a free distributive or CD mapping.  First suppose that $\morc = \iclosure \morc$, and and fix a split partition $(\dop \setb, \setb)$.  Then equation \eqref{eq_deflation_formula_split_a} implies that $ \bigwedge \mora(\setb) = \bigwedge (\morc \nu[\setb]) =  \morc \bigwedge \nu[ \setb] = \morc \iclosure \bigwedge \nu[\setb] = \morc ( \bigvee \nu[  \dop \setb]) \wedge \morc(\bigwedge \nu[ \setb]) = (\bigvee \mora( \dop \setb)) \wedge (\bigwedge \mora( \setb)) $, hence $\bigvee \mora( \dop \setb) \ge \bigwedge \mora( \setb)$.  This establishes one direction.  

For the converse, suppose  that $\bigwedge \mora( \setb) \le \bigvee \mora( \dop \setb )$ for every split partition $(\dop \setb, \setb)$.  Then one can argue that 
	\begin{align}
	\mora(\a) \le \mora( \dopisoo (\a) )	
	\label{eq_split_inequality}
	\end{align}
for each $\a\in \tob$.

Let $\dop \seta \in \axlh^2(\dop \tob)$ be given.  If $\dop \seta < 1$ then $\iclosure \dop \cpph \seta = \cpph \seta$.  Otherwise $\seta = 1$, and $\morc \iclosure 1 = \morc ( \cpph \predxii (1)) = \morc (\bigvee_{\aa \in \dop \tob} \nu_\aa) =\bigvee \mora(\dop \tob) \ge \bigwedge \mora( \emptyset) = 1 = \morc(1)$ by Equation \eqref{eq_deflation_formula_aa}.  Thus $\morc$ and $\morc \iclosure$ agree on $\Im(\dop \cpph)$.

Let $\dop \seta \in \axlh^2( \tob)$ be given.  If $\cpph (\seta)$ can be expressed in form $\bigwedge \nu[\setb]$ for some split partition $(\dop \setb, \setb)$, then $\morc \iclosure \cpph(\seta) = \morc \iclosure \bigwedge \nu[\setb] = \morc ( \bigvee \nu[  \dop \setb]) \wedge \morc(\bigwedge \nu[ \setb]) = (\bigvee \mora( \dop \setb)) \wedge (\bigwedge \mora( \setb)) = \bigwedge \mora( \setb) = \morc \bigwedge \nu[ \setb] = \morc \cpph(\seta)$.  Otherwise $\iclosure(\cpph \seta) = \dop \cpph (\dopisoo \seta) \wedge \cpph \seta$ by Lemma \ref{lem_deflation}, and in this case  $\morc \iclosure \cpph (\seta) = \morc \cpph (\seta) $ by Equation \eqref{eq_split_inequality}.  Thus $\morc$ and $\morc \iclosure$ agree on $\Im(\cpph)$.

It follows that $\morc \iclosure$ and $\morc$ agree on the image of $\dop \cpph \sqcup \cpph$.  If $\morc$ is a free \emph{distributive} mapping, then, since every element of $\dxlh \axlh(\dop \tob \sqcup \tob)$ can be expressed as a finite join of elements of form $\dop \cpph(\dop \seta) \wedge \cpph (\seta)$, it follows that $\morc \iclosure = \morc$.  If $\morc$ is a  \emph{CD} mapping, then, since every element of $\axlh^2(\dop \tob \sqcup \tob)$ can be expressed as a (possibly infinite) join of elements of form $\dop \cpph(\dop \seta) \wedge \cpph (\seta)$, and both $\morc$ and $\morc \iclosure$ preserve arbitrary joins, it again follows that $\morc = \morc \iclosure$.
\end{proof}

\begin{remark} Map $\iid \ii$ is a complete lattice homomorphism. Map $\iid \di$ is a \emph{join}-complete lattice homomorphism, i.e.\ it preserves binary meets and arbitrary joins.  However, $\iid \di$ may fail to preserve empty or infinite meets.  As an example, consider the case where $\tob$ is the disjoint union of $\toc_-$ and $\toc_+$, where $\toc_-$ is a down-closed subset of $\tob$ isomorphic to $\Z_{\ge 0}$ and $\toc_+$ is a up-closed subset of $\tob$ isomorphic to $\Z_{\le 0 }$. Let $(\dop \setb, \setb)$ be the split partition such that $\dop \setb = \toc_-$ and $\setb= \toc_+$.  Then 
	\begin{align*}
	\textstyle
	(\iid \di \iid \ii) (\bigwedge_{\a \in \setb} \nu_\a ) = (\bigvee_{\aa \in \dop \setb} \nu_{\aa} ) \wedge (\bigwedge_{\a \in \setb} \nu_\a )
	&&
	\textstyle	
	\bigwedge_{\a \in \setb}  \left ( (\iid \di \iid \ii) \nu_\a \right ) = (\bigwedge_{\pi(\aa) \in \setb} \nu_\aa ) \wedge  (\bigwedge_{\a \in \setb} \nu_\a ) 
	\end{align*}
The expression on the right contains $\dop \setb \cup \toc_-$, however the expression on the left does not.
\end{remark}

\section{Saecular persistence}
\label{sec_saecular}

Fix a functor
	\begin{align*}
	\f : \tob \to \ecata
	\end{align*}
from a totally ordered set $\tob$ to a p-exact category $\ecata$.  

\subsection{The saecular filtrations}

Recall that the category of $\tob$-shaped diagrams in $\ecata$ is p-exact (Lemma \ref{lem_subdiagrams}).  Diagram $\f$ is an object in this category, and the poset
$
\Sub_\f
$
is an order lattice under inclusion (also by Lemma \ref{lem_subdiagrams}).  We may therefore define maps $\K : \tob \to \Sub_\f$ and $\KK: \dop \tob \to \Sub_\f$ by 
	\begin{align}
	\KK(\tela) & = \bigwedge  \{ \dgb \in \Sub_\f : \dgb_\lela = 1_{\Sub(\f_\a)} \foral \lela \dople \tela \} 
	\label{eq_kkerdef}
	\\	
	\K(\tela) & = \bigvee \{ \dgb \in \Sub_\f : \dgb_\lela = 0_{\Sub(\f_\a)} \foral \lela \ge \tela \} .
	\label{eq_kerdef}	
	\end{align}
These suprema and infima always exist; for example, $\KK(i)_\a = \dga$ when $\a \le i$ and $\KK(i)\a = \Im( \dga(i \le \a))$ when $i \le \a$.  We call $\KK$ and $\K$ the \emph{saecular filtrations} of $\dga$.  The \emph{saecular  CDH, CDIH, BDH, and BDIH} are 
	\begin{align*}
	\SCD(\mora) &= \FCD(\beta)
	&
	\SBD(\mora) = \FBD(\bar \beta)	
	\\
	\SCDI(\mora) &= \IFCD(\beta)		
	&
	\SBDI(\mora) = \IFBD(\bar \beta)		
	\end{align*}
respectively, where $\beta: = \KK \sqcup \K : \dop \tob \sqcup \tob \to \Sub_\mora$ is the copairing of $\KK$ and $\K$ in category $\catps$.   For economy of notation, we often write $|\seta|$ for the value taken by any one of these homomorphisms on a down-closed set $\seta$.  We call $\beta$ the \emph{saecular bifiltration}.




\begin{theorem}
\label{thm_sacular_existence_iff}
Diagram $\dga$
	\begin{enumerate*}
	\item admits a saecular BDH iff each saecular filtration factors through a complete sublattice of $\Sub_\dga$.
	\item admits a saecular BDH iff it admits a saecular BDIH.
	\item admits a saecular CDH iff it admits a saecular CDIH.	
	\end{enumerate*}
In particular, $\mora$ admits a BDH if $\ecata = \in \{R\catmod, \catmod R\}$, or $\ecata$ is any abelian category with complete subobject lattices.
\end{theorem}
\begin{proof}
Lattice $\Sub_\dga$ is modular, by Proposition \ref{prop_pexact_modular_principles}.  Assertion 1 thus follows from Theorem \ref{thm_free_distributive_mapping}.

Now suppose that each saecular filtration factors through a complete sublattice of $\Sub_\dga$ (otherwise there exists no saecular homomorphism of any kind).  Then meet and join are well-defined for every subset of $\Im(\KK)$ and every subset of $\Im(\K)$.  If $(\dop \setb, \setb)$ is a split partition of $\tob$, then $\bigvee_{\aa \in \dop \seta } \KK(\aa)$ is the minimum subdiagram $\morb$ of $\dga$ such that $\dgb_\aa = \dga_\aa$ for each $\aa \in \dop \seta$, and $\morc = \bigwedge_{\a \in \seta} \K(\a)$ is the maximum subdiagram of $\mora$ such that $\morc_\a = 0$ for all $\a \in \seta$.  Thus $\bigvee \KK[ \dop \seta ] \ge \bigwedge \K [ \seta ]$.  The desired conclusion therefore follows from Lemma \ref{lem_closure_comp_eq}.
\end{proof}

A saecular homomorphism is \emph{s-natural (or simply natural) at $\seta$} if 
	\begin{align}
	\dia |\seta|_\aa = |\seta \wedge \nu(\aa)|_\a = |\seta|_\a \wedge \KK(\aa)_\a
	\label{eq_dicap}
	\\
	\iia |\seta|_\a = |\seta \vee \nu(\a) |_\aa = |\seta|_\aa \vee \K(\a)_\aa
	\label{eq_iicup}
	\end{align}
for any $\aa \dople \a$.   A saecular homomorphism is \emph{s-natural (or simply natural)} if it is natural at every set in its domain.  Conversely, a set $\seta$ is  \emph{natural with respect to $\dga$} if there exists a saecular homomorphism of $\dga$ which is natural at $\seta$.  We defer the proof of Theorem \ref{thm_cdf_is_proper} to \S\ref{sec_prove_factors_are_intervals}.

\begin{theorem}
\label{thm_cdf_is_proper}
For any chain diagram $\dga$,
	\begin{enumerate*}
	\item $\SCD(\mora)$ is natural, if it exists.
	\item $\SBD(\mora)$ is natural if $\SCD(\mora)$ exists.
	\item $\SBD(\mora)$ is natural iff it is natural at every set in  $\seta \in \Im(\dop \cpph) \cup \Im(\cpph)$.
	\end{enumerate*}
\end{theorem}

\begin{remark}
\label{rmk_canonical_realizations_of_sfun}
Subobjects, quotient objects, and subquotients in a p-exact category are defined only up to unique isomorphism.  However, canonical realizations of a given subquotient are sometimes  available, e.g.\ via the construction of a quotient group via cosets.  

When one wishes to work with such canonical realizations, one can define $\sfun$ on an object $X$ in $\lclc(\axl{\dop \tob \sqcup \tob})$ as $|\seta| /|\setb|$, where $(\setb, \seta)$ is the smallest nested pair in $\axlh^2( \dop \tob \sqcup \tob)$ such that $X = \seta - \dopiso(\setb)$ (concretely, $\seta = \dsh X$ and $\setb = \dopiso^{-1}(\seta - X)$).    We refer to the resulting functor  the \emph{canonical realization} of $\sfun$, and it is uniquely defined.  The canonical realization of $\sfuni$ is defined similarly. In particular, the canonical realization of $\sfuni$ satisfies $\sfuni \{ \itva \} = \xfrac{ \dsh \itva }{ \pdsh \itva }$ for each nonempty $\itva \in \allitv{ \tob}$.
\end{remark}

\subsection{Exact functors and interval factors}
\label{sec_spd_pexact}

The \emph{saecular CD functor} is the locally closed functor of $\SCD(\dga)$, denoted $\SCDF(\dga)$.  Locally closed functors $\SCDIF(\dga)$, $\SBDF(\dga)$, and $\SBDIF(\dga)$ are defined similarly.  For economy of notation, we will often write 
	$
	|\seta|
	$ 
for the value taken by one of this maps on a locally closed set $\seta$.

\begin{theorem}
\label{thm_saecular_factors_are_interval}
If $\SBDI(\mora)$ is natural at $\dsh \itva$ and $\pdsh \itva$, then 
	\begin{align*}
	| \{\itva \}| \in \allifun(\itva).
	\end{align*}
\end{theorem}
\begin{proof}
The assertion equates to Theorem \ref{thm_single2interval}, which we will prove in \S\ref{sec_prove_factors_are_intervals}.
\end{proof}

We call $| \{\itva \} |$ the  \emph{$\itva^{th}$ saecular factor} of $\dga$.  Technically, the $\itva^{th}$ saecular factor is defined only up to unique isomorphism.  However, if $\ecata$ has canonical realizations of subquotients then we can define it uniquely, as per Remark \ref{rmk_canonical_realizations_of_sfun}.  In this formulation, $| \{\itva \} |= \xfrac{ \dsh \itva }{ \pdsh \itva }$ for each $\itva$.  We refer to this unique construction as the \emph{canonical realization}.

\begin{counterexample}
\label{cx_s_unnatural_sets}
An interval factor $|\{\itva\}|$ may not be a type-$\itva$ interval functor if $\SBDI(\dga)$ is not natural at $\dsh \itva$.  

Let $X^m$ denote the set of finitely supported functions $\Z \to \fielda$  that vanish on $\Z_{> m}$.  Let $\dgb: \Z_{< 0} \to \fielda\catvect$ be the functor that carries $m$ to $X^m$ and  $m \le n$  to the inclusion $X^m \su X^n$.   Extend $\dgb$ to a functor $\dga: \Z_{\le 0} \to \fielda\catvect$ such that $\mora_0 = \fielda$ and $\mora( m \le 0) : (x_i)_{i \in \Z} \mapsto \sum_i x_i$.  Then $\SBD(\mora)$ exists, and $\SBD(\mora) \{\Z_{\le 0}\}_m = 0$ iff $m \neq 0$.  In particular, $\SBD(\mora) \{\Z_{\le 0}\} = |\{\Z_{\le 0} \}|$ is not an interval functor support type $\Z_{\le 0}$.
\end{counterexample}

\subsection{Classification of point-wise finite dimensional chain diagrams, redux}

The existence of Krull-Schmidt decomposition for arbitrary chain functors $\tob \to \fielda \catvectf$ was first shown by Botnan and Crawley-Boevey \cite{botnan2020decomposition}.  The saecular framework provides a short new proof -- in addition to a structural connection with saecular functors.

\begin{theorem}
\label{thm_findimdecomp}
If $\tob$ is a totally ordered set and $\ecata = \fielda \catvectf$, then every functor $\dga: \tob \to \ecata$ admits a saecular CDI functor $\sfun$.  Moreover,
	\begin{align*}
	\dga \cong \bigoplus_{\itva \in \allitv{\tob}} \sfun \{\itva \}.
	\end{align*}
\end{theorem}
\begin{proof}
In this case $\dga$ satisfies the well ordered criterion of Theorem \ref{thm_saecularexistencecriteria}, so the saecular CDF functor exists.  For each $\itva \in \allitv{\tob}$ fix an indexed family of sets $(B_\tela^\itva)_{\tela \in \tob}$ such that $B_\tela^\itva \su \sfun( \downarrow \itva)$ and  $q(B_\tela^\itva)$ is a basis for  $\sfun (\{\itva\})_\tela$ for each $\tela$, where $q$ is the quotient map 
	$
	\sfun( \downarrow \itva \onto \{\itva \})
	$.

We call each $B_\tela^\itva$ an \emph{$\itva$-prebasis} in $\f_\tela$.  We claim that one can choose $(B_\tela^\itva)$ such that $\dga(\tela \le \telb)(B_\tela^\itva) = B_\telb^\itva$ for all $\tela \le \telb \in \tob$.  Such a family will be called a \emph{consistent family of prebases}.  Existence of a consistent family is clear when $\itva$ contains a minimum element.  When it does not, first consider the class $\pseta$ of consistent prebases which are defined, not on all of $\itva$ as desired, but on up-sets of $\itva$.  Equip $\pseta$ with the partial order under inclusion.  Each chain in $\pseta$ has an upper bound in $\pseta$, so  $\pseta$ contains a maximal element $(C^\itva_\tela)_{\tela \in \itvb}$ by Zorn's lemma, where $\itvb$ is an up-set of $\itva$.  If $\itvb \neq \itva$, choose an element $\tela \in \itva - \itvb$.  Since $\f_\tela$ is finite-dimensional, the indexed family of inverse images $(\seta_\telb)_{\telb \in \itvb}$ defined by
	\begin{align*}
	\seta_\telb: = \f(\tela \le \telb) \inv(C^\itva_\telb)
	\end{align*}
takes at most finitely many distinct values.  Moreover, $\mora(\tela \le \telb)(\seta_\telb) = C^\itva_\telb$ for all $j \in \itvb$ since, $\sfun(\tela \le \telb)$ carries $\sfun(\dsh \itva)_\tela$ surjectively onto $\sfun(\dsh \itva)_\telb$.  Fix the minimum value $\setb = \min_\telb( \seta_\telb) \su \f_\tela$,  select an $\itva$-prebasis $D^\itva_\tela$ contained in $\setb$, and define $D_\telb^\itva = \f(\tela \le \telb)(D^\itva_\tela)$ for all $\telb$ greater than $\tela$ in $\itva$.  Then $C \lneq D$ in $\pseta$, a contradiction.  Thus $C$ is defined on all of $\itva$.  This establishes the claim.   

To complete the proof, it remains only to check that $B_\tela: = \bigcup_{\itva \in \allitv{\tob}} B_\tela^\itva$ is a basis for $f_\tela$, for each $\tela$.  This can be verified via Lemma \ref{lem_complete_map_to_fdim_vspace_basis}, by fixing a linearization $\lina$ of poset $\allitv{\tob}$, and verifying that $B_i^\itva$ is a prebasis for $\sfun ( \dsh_\lina \itva)_\tela / \sfun( \pdshp{\lina} \itva)_\tela \cong \sfun (\{ \itva \})$.
\end{proof}

As in the proof of Theorem \ref{thm_findimdecomp}, given a vector space $V$ with nested subspaces $U \su W$, let us define a \emph{prebasis} of $W/U$ to be a linearly independent subset $B \su W$ such that $q(B)$ is a basis of $W/U$, where $q$ is the projection map $\seta \mapsto W/U$.

\begin{lemma}
\label{lem_complete_map_to_fdim_vspace_basis}
Let $\lina$ be a totally ordered set, $V$ be a vector space, and $\cdf: \axl{ \lina } \to \Sub_V$ be a complete lattice homomorphism.  Suppose that $\Im(\cdf)$ is a well-ordered chain in $\Sub_V$, and let $B^\a$ be a pre-basis of $\morc( \dsh \a) / \morc( \pdsh \a)$ for each $\a \in \lina$.  Then $\bigcup_{ \a \in \lina} B^\a$ is a basis for $V$.
\end{lemma}
\begin{proof}
Let $\chaina = \Im(\cdf)$, and fix $v \in V$.  There exists a unique covering relation $\setc \lessdot \setc'$ in $\chaina$ such that $ v \in \setc' - \setc$.  Likewise, there exists a unique covering $\seta <: \seta'$ in $\axl{\lina}$ such that $\setc' = \cdf \seta'$ and $\setc = \cdf \seta$;
 concretely, $\seta' = \dsh \a_v$, where $\a_v = \min \{ x \in \lina : v \in \cdf \dsh x \}$, and  $\seta = \{ \a \in \lina : v \notin \cdf \dsh \a \} = \pdsh \a_v$ (note that every covering relation in $\axl{\lina}$ has form $\pdsh \a \lessdot \dsh \a$ for some $\a$).  
Writing $\Suq(\chaina)$ for the set of successors in $\chaina$, we may therefore define a unique function $\xi: \Suq(\chaina) \to \lina$ such that (i) $\setc = \cdf \dsh (\xi_\setc)$  for each $\setc \in \chaina$, and (ii) $V$ partitions as a disjoint union 
	\begin{align}
	V 
	= 
	\bigsqcup_{\setc \in \Suq(\chaina)} (\cdf \dsh \xi_\setc) - (\cdf \pdsh \xi_\setc)
	= 
	\bigsqcup_{\a \in \Im(\xi)} (\cdf \dsh \a) - (\cdf \pdsh \a)
	=
	\bigsqcup_{\a \in \lina} (\cdf \dsh \a) - (\cdf \pdsh \a)	
	\label{eq_wo_prebasis_disjoint_unions}
	\end{align}
	
Now let $B = (B^\a)_{\a \in \lina}$ be given, and let $\lina' = \Im(\xi)$.  A routine exercise shows that $\bar B: = \bigcup_{\a \in \lina} B^\a = \bigcup_{\a' \in \lina'} B^{\a'}$ is linearly independent. Let $\holdsfor(\a)$ denote the statement that $  \bigcup_{ \a' \le \a} B^{\a'}$ is a basis for $\cdf \dsh \a$.  By \eqref{eq_wo_prebasis_disjoint_unions} one has $\bigcup_{\a \in \lina'} \cdf \a = V$.  Thus if $\holdsfor(\a)$ for all $\a \in \lina'$, then $\bar B$ contains a basis for $V = \bigcup_{\a \in \lina'} \cdf \dsh \a$.  Since $\lina'$ is well ordered, it therefore suffices, by transfinite induction, to show that ($\holdsfor(\a')$ for all $\a' < \a$) $\implies \holdsfor(\a)$ for each $\a \in \lina'$.  Indeed, ($\holdsfor(\a')$ for all $\a' < \a$) implies that $\bar B$ contains a basis for $\cdf \pdsh \a$, whence $\bar B$ contains a basis for $\cdf \dsh \a$, by definition of $B^\a$.  The desired conclusion follows.
\end{proof}

\begin{corollary}
\label{cor_wo_bifiltration_prebasis}
Let $\tob$ be a well-ordered set, $V$ be a vector space, and $\cdf: \axlh^2(\dop \tob \sqcup \tob) \to V$ be a complete lattice homomorphism.  For each $\seta \in \axl{\dop \tob \sqcup \tob}$, let $B^\seta$ be a prebasis for $\cdf \dsh \seta / \cdf \pdsh \seta$.  Then $\bigcup_{\seta \in \axl{\dop \tob \sqcup \tob}} B^\seta$ is a basis for $V$.
\end{corollary}
\begin{proof}
Let $\lina$ be any linearization of $\axl{\dop \tob \sqcup \tob} \cong \axl{\dop \tob} \times \axl{\tob} \cong \axl{\tob}^2$.  Then $\lina$ is well-ordered, by Corollary \ref{cor_wolinearization}.  For each $\seta \in \axl{\dop \tob \sqcup \tob}$ there exists a regularly induced isomorphism 
	$
	\cdf \dsh_{\axl{\dop \tob \sqcup \tob}} \seta / \cdf \pdsh_{\axl{\dop \tob \sqcup \tob}} \seta 
	\cong 
	\cdf \dsh_{\lina} \seta / \cdf \pdsh_{\lina} \seta \cong \cdf \{ \seta \}
	$.  
Therefore $B^\seta$ is a prebasis for $\cdf \dsh_{\lina} \seta / \cdf \pdsh_{\lina} \seta$.  The desired conclusion therefore follows from Lemma \ref{lem_complete_map_to_fdim_vspace_basis}.
\end{proof}

\subsection{New decomposition results}

Theorem \ref{thm_wo_linear_decomp} is a natural companion to Theorem \ref{thm_findimdecomp}.  It indicates, among other things, that we can remove restrictions on dimension when $\tob$ is well-ordered.

\begin{theorem}
\label{thm_wo_linear_decomp}
The following are equivalent for any totally ordered set $\tob$: 
	\begin{enumerate*}
	\item  $\tob$ is well ordered
	\item every $\tob$-shaped diagram in $\fielda \catvect$ splits as a direct sum of interval functors, for some fixed choice of  $\fielda$
	\item every $\tob$-shaped diagram in $\fielda \catvect$ splits as a direct sum of interval functors, for every choice of $\fielda$
	\end{enumerate*}
\end{theorem}

\begin{proof}
It suffices to show that $(1)$ implies $(3)$, and that not (1) implies not (2).

First assume (1), and let $\fielda$ and $\dga: \tob \to \fielda \catvect$ be given.  Then $\dga$ admits a saecular CD functor $\sfun$, by the well-ordered criterion (Theorem \ref{thm_saecularexistencecriteria}).  For each $\itva \in \allitv{\tob}$, let  $B^\itva_{\min(\itva)}$ be a pre-basis for $\left( \xfrac{\dsh \itva}{\pdsh \itva}  \right)_{\min(\itva)}$. Then define an indexed family of sets $(B^\itva_\a)_{\a \in \tob}$ by
	\begin{align*}
	B^\itva_\a
	=
		\begin{cases}
		\dga( \min(\itva) \le \a) (B^\itva_{\min(\itva)}) & \a \in \itva \\
		\emptyset 	& else\\		
		\end{cases}
	\end{align*}
Since $\sfun \{\itva\} =\xfrac{\dsh \itva}{\pdsh \itva}$ is an $\itva$-interval module, it follows that $B^\itva_\a$ is a pre-basis for $\xfrac{\dsh \itva}{\pdsh \itva}$, for all $\a$.  Thus $\bigcup_{\itva \in \allitv{\tob}} B^\itva_\a$ is a basis for $\dga_\a$, by Corollary \ref{cor_wo_bifiltration_prebasis}.  It follows that $\dga$ is an internal direct some of 1-dimensional interval functors, with one summand for each basis vector in $B^\itva_{\min(\itva)}$.  Thus (1) implies (3).

Now assume that (1) fails, i.e.\ that $\tob$ is not well-ordered.  Then there exist poset homomorphisms $p: \tob \leftrightarrows \Z_{\le 0}: i$ such that $p\circ i$ is identity on $\Z_{\le 0}$.  If $\dga: \Z_{\le 0} \to \fielda\catvect$ is any functor, then $\dga = \dga (pi) = (\dga p)i$.  Thus every functor $\Z_{\le 0} \to \fielda\catvect$ can be expressed in form $\dgb i$ for some $\dgb: \tob \to \cata$.  Note that the composite $\dgb i$ splits as a direct sum of interval functors whenever $\dgb$ does.  In Counterexample \ref{cx_s_unnatural_sets}, we construct a functor $\dgb: \Z_{\le 0} \to \fielda\catvect$ where $\SBD(\dga)$ is not natural.  Functor $\dgb$  admits  no saecular CD homomorphism (by Theorem \ref{thm_cdf_is_proper}) and therefore  no decomposition as a direct sum of interval functors (by Theorem \ref{thm_saecular_capture_gabriel}).  Consequently, $\dgb i$ admits no decomposition as a direct sum of interval functors.  Thus not (1) implies not (2).  The desired conclusion follows.
\end{proof}

Up to canonical isomorphism, the saecular factors capture each summand of a direct sum of interval diagrams in a category of modules, in the following sense:

\begin{theorem} 
\label{thm_saecular_capture_gabriel}
Suppose that $\ecata$ is a category of modules over a ring and 
	$
	\dga = \bigoplus_{\itva \in \allitv{\tob}} \dgb(\itva)
	$
for some indexed family of interval modules $\dgb$ such that $\dgb(\itva) \in \allifun(\itva)$ for each $\itva$.  Then $\mora$ admits a saecular CDI functor $\sfun$, and  $\sfun(\seta) = \bigoplus_{\itva \in \seta} \dgb(\itva)$ for each $\seta \in \axlh \allitv{\tob}$.  

Consequently, there exists a canonical isomorphism
	\begin{align*}
	\sfun\{ I \} \cong \dgb(\itva)
	\end{align*}
for each interval $\itva \in \allitv{\tob}$.
\end{theorem}
\begin{proof}
Existence follows from the proof of Theorem \ref{thm_saecularexistencecriteria}, as does the formula for $\sfuni(\seta)$.  The claim for $\sfuni\{\itva\}$ follows.
\end{proof}

\subsection{Persistence diagrams: old and new}

The \emph{saecular persistence diagram} of  $\dga$ is defined as the function 
	\begin{align*}
	\spd : \allitv{\tob} \to [\tob, \ecata]
	&&
	\itva \mapsto  
	| \{\itva \}|
	\end{align*}
Equivalently, $\spd(\itva) = \frac{\cdf(\dsh \itva)}{\cdf(\pdsh \itva)}$, where $\cdf = \SBDI(\mora)$. Where context leaves room for confusion, we write $\spd_\dga$ to emphasize dependence of $\dga$. The saecular persistence diagram is well defined iff $\SBDI(\mora)$ exists.
	
Note that $\spd$ carries intervals to \emph{chain diagrams}, not to integers.  In particular, Theorem \ref{thm_saecular_factors_are_interval} states that the saecular persistence diagram carries  $\itva$ to an interval functor of support type $\itva$, so long as $\SBDI(\mora)$ is natural at $\dsh \itva$.  By contrast, the classical persistence diagram, $\cpd$,  assigns an \emph{integer} to each interval.  More precisely, if $\tob \to \fielda\catvectf$ is a chain functor and $\dga \cong \bigoplus_{\itva \in \allitv{\tob}} \dgb(\itva)$, then $\cpd$ can be defined as the unique function such that  $\cpd(\itva) = \dim( \dgb_\itva)$ for all $\itva$, where $\dim(\dgb_\itva)$ is the dimension of the object type of $\dgb_\itva$.

\begin{corollary}
\label{thm_saecular_pd_generalizes_linear_pd}
If $\ecata$ is the category of finite-dimensional vector spaces over a field $\fielda$, then the saecular persistence diagram exists and
	\begin{align*}
	\cpd(\itva) = \dim (\spd(\itva))
	\end{align*}
for each interval $\itva$.
\end{corollary}
\begin{proof}
Follows from Theorems \ref{thm_findimdecomp} and \ref{thm_saecular_capture_gabriel}. 
\end{proof}

\subsection{Subsaecular series}
\label{sec_saecularcalendars}

There is a useful analogy to be drawn between the interval factors of a chain functor, on the one hand, and the composition factors of finite-height ring module, on the other.    Theorem \ref{thm_cseries} formalizes this analogy.

\begin{theorem} 
\label{thm_cseries}
Suppose that $\tob$ is well ordered.  If $\lina$ is a linearization of   $\allitv{\tob}$, and $\Sub_\f$ is an upper continuous lattice, then there exists a unique $\forall$-complete lattice homomorphism $\cdflin: \axl{\lina} \to \Sub_\f$ such that
	\begin{align}
		\frac
		{ \cdflin \dsh_\lina   (\itva) }
		{ \cdflin \pdshp\lina (\itva) }		
	\in
	\allifun(\itva )
	\label{eq_intervalconstraints}
	\end{align}
for each $\itva \in \allitv{\tob}$.  In particular, $\dga$ admits a SCDI homomorphism $\sfun$, and $\cdflin = \sfun|_{\axlh^2(\lina)}$.
\end{theorem}

\begin{definition}
We call $\cdflin$ the \emph{subsaecular series} of $\dga$ subordinate to $\lina$.  We call \eqref{eq_intervalconstraints} the $\itva^{th}$ subsaecular factor of $\lina$. 
\end{definition}

\begin{remark}
If $\tob$ is well-ordered then so, too, are $\lina$ and $\axlh(\lina)$, by Corollary \ref{cor_wolinearization}.
\end{remark}

\begin{corollary}[Jordan-H\"older for subsaecular series]
\label{cor_saecular_factors_noether_isomorphic}
The $\itva^{th}$ saecular factors of any two linearizations, $\lina$ and $\lina'$, are Noether isomorphic. 
\end{corollary}

\begin{proof}
Let $\itva \in \allitv{\tob}$ be given.  
Since $ (\dsh_\lina \itva) \backslash (\pdsh_\lina \itva) = \{\itva\} = (\dsh_{\lina'}\itva) \backslash (\pdsh_{\lina'} \itva )$, Lemma \ref{lem_phi_M_relation_for_semitopologies}, implies that $\frac{ \sfun (\dsh_\lina \itva)}{ \sfun (\pdsh_\lina \itva)} \mdom  \frac{  \sfun (\dsh_{\lina'} \itva)} { \sfun  (\pdsh_{\lina'}\itva) }$.  The desired conclusion therefore follows from  Theorem \ref{thm_mutualdomination}. 
\end{proof}

\begin{proof}[Proof of Theorem \ref{thm_cseries}]
The saecular CDI functor $\sfun$ exists by the well-ordered criterion (Theorem \ref{thm_saecularexistencecriteria}), and the restriction $\sfun|_{\axl{\lina}}$  satisfies \eqref{eq_intervalconstraints} by Theorem \ref{thm_saecular_factors_are_interval}.  This establishes existence.

It remains only to prove uniqueness.  For this purpose we identify each interval $\itva$ with the disjoint union of down sets $\dop \seta \sqcup \seta$ such that $\itva = \dopiso(\dop \seta) - \seta$, as per custom.  Fix any $\forall$-complete lattice homomorphism $\cdflin$ satisfying criterion \eqref{eq_intervalconstraints}.  Given any subset $\pseta \su \axl{\lina}$, let $\holdsfor(\pseta)$ denote the statement that
	\begin{align}
	\cdflin \setc = \sfun \setc
	\label{eq_transfinitegoal}
	\end{align}
for each $\setc \in \pseta$. The chain $\axlh(\lina)$ is well-ordered, by Corollary \ref{cor_wolinearization}.  By transfinite induction, therefore, it suffices to show that $\holdsfor(\axl{\lina}_{\lneq \setb})$ implies $\holdsfor(\{\setb\})$	for each $\setb \in \axl{\lina}$.  By completeness, we may restrict to the successor case, so $\setb$ covers a unique element $\mathring \setb$, and    there exist $\dop \setd, \setd$ such that 
$\mathring \setb = \pdshp{\lina} \dop \setd \sqcup \setd$ and $\setb = \dsh_{\lina} \dop \setd \sqcup \setd$.

If $\setd - \dopiso(\dop \setd) = \emptyset$, then 
	$
	\cdflin \setb / \cdflin \mathring \setb
	$
and  $\sfun  \setb / \sfun  \mathring\setb$ are interval chain diagrams with empty support, hence 0 diagrams. Thus $\cdflin \setb = \cdflin \mathring \setb$ and $\sfun  \setb = \sfun \mathring \setb$.  By induction, the righthand sides of both these equations agree, as do the lefthand sides; this establishes the inductive step. 
%
Therefore assume, without loss of generality, that $\setd - \dopiso(\dop \setd)$ is nonempty, hence there exists a minimum element
	\begin{align*}
	\eld = \min (\setd -  \dopiso(\dop \setd)).
	\end{align*}
	
Since $\cdflin \setb / \cdflin \mathring \setb$ has support type $\setd - \dopiso(\dop \setd)$, one has
	\begin{align*}
	\cdflin \setb / \cdflin \mathring \setb
	=
	 \cdflin \setb / \sfun \mathring \setb
	\le 
	\KKs_{\dga/ \sfun \mathring \setb}(\dsh \dop \setd) \wedge \Ks_{\dga/ \sfun \mathring \setb}(\dsh \setd)
	=
	\sfun \setb / \sfun \mathring \setb.
	\end{align*}
Therefore $\sfun \mathring \setb \le \cdflin \setb \le  \sfun \setb $. 
\\

\noindent \ul{Claim 1.} $\cdflin \setb < \sfun \setb \iff \cdflin \setb_\eld < \sfun \setb_\eld$.  \emph{Proof.}  If $\eld \le \a$ then $(\sfun \setb / \sfun \mathring \setb)(\eld \le \a) $ is epi, and therefore $\sfun \mathring \setb_\a  \vee \dga( \eld \le \a) \di \sfun \setb_\eld = \sfun \setb_\a$.  Thus 
	$
		\cdflin \setb_\eld = \sfun_\eld 
	$
implies
	$$
		\cdflin_\a
		\ge
		 \cdflin \mathring \setb_\a \vee \dga( \eld \le \a) \di \cdflin \setb_\eld 
		=  
		 \sfun \mathring \setb_\a \vee \dga( \eld \le \a) \di \sfun \setb_\eld
		= 
		\sfun \setb_\a
		\ge 
		\cdflin_\a.
	$$
This establishes one direction.  The converse is trivial.
\\

By Claim 1, we have  $\cdflin \setb \neq \sfun \setb \iff \cdflin \setb_\eld \lneq  \sfun \setb_\eld $, in which case 
	\begin{align}
	\{ \seta \in \axl{\lina} 
	 	: 
	 	\cdflin \setb_\eld 
		\lneq
		\cdflin \seta_\eld \wedge \sfun \setb_\eld   
	\}	
	 \label{eq_obstructions}
	\end{align}
is nonempty.  Assume, for a contradiction, that such is the case.  Then  \eqref{eq_obstructions} contains a minimum element $\setc$. 
\\

\noindent \underline{Claim 2.}  $\setc$ is a successor. \emph{Proof.} If not, then $\setc$ is the supremum in $\axl{\lina}$ of $\sete = \{\seta \in \axl{\lina} : \seta \lneq \setc\} $, and by upper continuity 
	\begin{align*}
	\textstyle
	\cdflin \setb 
	\lneq
	\cdflin \setc \wedge \sfun \setb 
	= 
	\bigvee (\cdflin \sete) \wedge \sfun \setb
	=
	\bigvee \{ \cdflin \seta \wedge \sfun \setb : \seta \in \sete \}.
	\end{align*}
By Claim 1, it follows that 
	\begin{align}
	\textstyle
	\cdflin \setb_\eld
	\lneq
	\cdflin \setc_\eld \wedge \sfun \setb _\eld
	= 
	\bigvee (\cdflin \sete)_\eld \wedge \sfun \setb_\eld
	=
	\bigvee \{ \cdflin \seta_\eld \wedge \sfun \setb_\eld : \seta \in \sete \}.
	\label{eq_atomiccontradiction}
	\end{align}	
However, if $\seta \in \sete$ then $\seta$ is a strict lower bound of $\setc$, the minimum of \eqref{eq_obstructions}, hence $\seta$ lies outside \eqref{eq_obstructions} and therefore   
$
	 	\cdflin \setb_\eld 
		\ge 
		\cdflin \seta_\eld \wedge \sfun \setb_\eld  
$. Therefore   the righthand side of \eqref{eq_atomiccontradiction} bounds $\cdflin \setb$ from below, a contradiction. Thus $\setc$ is a successor.  
\\

Concretely, Claim 2 implies that there exists $\dop \setx \sqcup \setx \in \lina$ such that  $\setc = \dsh_\lina \dop \setx \sqcup \setx $ is the successor to $\mathring \setc :=  \pdshp{\lina} \dop \setx \sqcup \setx$.
\\

\noindent \underline{Claim 3.}
$\setb \le \mathring \setc$.  \emph{Proof.} By hypothesis $\cdflin \setb <  \sfun \setb$, hence $\setb$ lies outside \eqref{eq_obstructions}.  Therefore $\setb < \setc$, hence $\setb \le \mathring \setc$.
\\

The first of the following relations holds because minimality of $\setc$ in \eqref{eq_obstructions} implies $\cdflin \mathring \setc \wedge \sfun \setb \lneq \cdflin \setc \wedge \sfun \setb$.  The second combines the Noether isomorphism $\soa / \soa\wedge \sob \cong \soa \vee \sob/ \sob$, where $\soa = \cdflin \setc \wedge \sfun \setb$ and $\sob = \cdflin \mathring \setc$, with an application of the modular law to rearrange parentheses in the numerator.
The third  holds for each  $\a \notin \setd$, since, by Claim 3, $\mora(\eld \le \a) \di( \cdflin \setc \wedge \sfun \setb )_\eld \le \mora (\eld \le a) \di \sfun \setb_\eld \le \sfun \mathring \setb_\a  = \cdflin \mathring \setb_\a \le \cdflin \mathring \setc_\a$.
	\begin{align}
	0 
	\lneq
		\left(
		\frac
		{\cdflin \setc \wedge \sfun \setb}
		{\cdflin \mathring \setc \wedge \sfun \setb}
		\right)_\eld	
	\cong
		\left(
		\frac
		{	\cdflin \setc			
			\wedge 
			(\sfun \setb \vee \cdflin \mathring \setc)}
		{	\cdflin \mathring \setc }
		\right)_\eld
	\le
		\left(
		\frac
		{\cdflin \setc}
		{\cdflin \mathring \setc}
		\right)
		(\eld \le \lela) \ii 0
	\label{eq_sss_modular}
	\end{align}
In particular, map	$(\cdflin \setc / \cdflin \mathring \setc)(\eld \le \lela)$ has nontrivial kernel for each $\lela \notin \setd$.  Since $(\cdflin \setc / \cdflin \mathring \setc) \in \allifun ( \setx -  \dopiso(\dop \setx))$, it follows that $\a \notin \setd \implies \a \notin \setx$, hence $\setx \su \setd$.  Likewise, since $(\cdflin \setc / \cdflin \mathring \setc)$ is nonzero at $v = \min(\setd - \dopiso(\dop \setd))$, we have   $\dop \setx \su \dop \setd$.  Therefore $\setc \le \setb$.  But this contradicts Claim 3.  Thus $\cdflin \setb = \sfun \setb$, which was to be shown.
\end{proof}

\subsection{Uniqueness and characterization of saecular objects}
\label{sec_uniqueness}

Since saecular homomorphisms are formally defined as free homomorphisms, they are unique when they exist.  However, it is natural to ask whether the same objects have alternate characterizations, e.g.\  defining properties that do not rely on  saecular filtrations explicitly.  If $\tob$ is well ordered and $\Sub_\dga$ is upper-continuous, then the answer is yes:

\begin{theorem}
Suppose that $\tob$ is well-ordered and the subobject lattice of $\dga: \tob \to \ecata$ is upper-continuous.  Then $\dga$ admits a saecular CDI homomorphism, and  $\sfun:=\SCDI(\dga)$ is the unique $\forall$-complete lattice homomorphism $\axlh{\allitv{\tob}} \to \Sub_\dga$ such that 
	\begin{align}
	\sfun \{\itva \} \in \allifun(\itva)
	\label{eq_cx_interval_criterion}
	\end{align}
for all $\itva \in \allitv{ \tob }$.
\end{theorem}
\begin{proof}
Follows from uniqueness of subsaecular series (Theorem \ref{thm_cseries}), since each $\seta \in \axlh \allitv{\tob}$ lies in at least one lattice of form $\axl{\lina}$, for some linearization $\lina$ of $\allitv{\tob}$.
\end{proof}

When $\tob$ is not well ordered, the answer may be no:

\begin{counterexample}[Non-uniqueness of complete homomorphisms with interval factors]  
\label{cx_wo_needed_for_uniqueness}
Several distinct $\forall$-complete lattice homomorphisms may satisfy \eqref{eq_cx_interval_criterion},  when $\tob$ is not well ordered.  Uniqueness may even fail  when $\Im(\KK)$ and $\Im(\K)$ are \emph{finite} subsets of $\Sub_\mora$.

To illustrate, let $\tob$ be the unit interval $[0,1] \su \R$, and let $\dga: \tob \to \cmc$ be the functor such that
	\begin{align*}
	\dga_\lela = \tob && \dga(\aa \dople \a) = \id
	\end{align*}
for all $\a \in \tob$ and $\aa \dople \a$.  Then $\Im(\KK) \cup \Im(\K) = \{ 0, \dga\}$, hence $\dga$ admits a saecular CDI homomorphism $\sfun$.

Let $\Psi: \axlh{\allitv{\tob}} \to \tob$ be the complete lattice homomorphism  described in Example \ref{ex_unit_interval_map}, and let $X: \tob \to \Sub_\dga$ be the complete lattice homomorphism such that $X(t)_\a = [0,t]$ for all $\a \in \tob$.  Then $\sfun$ and $\Xi: = (X \circ \Psi): \axlh{\allitv{\tob}}\to \Sub_\dga$ are both complete lattice homomorphisms.  The saecular CDI functor is natural, and therefore carries singletons to interval diagrams by Theorem \ref{thm_saecular_factors_are_interval}.  We claim that $\Xi$ satisfies the same condition.  In fact,   $ \frac{ \Xi(\dsh_{\axl{\dop \tob \sqcup \tob}} \itva) } { \Xi(\pdsh_{\axl{\dop \tob \sqcup \tob}} \itva ) }  = 0 \in \allifun(\itva)$ by Equation \eqref{eq_psi_dsh_equa_pdsh}.

Thus $\sfun$ and $\Xi$ are both complete homomorphisms; their locally closed functors both carry  singletons to interval diagrams, but the interval diagrams (\emph{a fortiori}, the homomorphisms) are distinct.
\end{counterexample}

\begin{counterexample}[Non-uniqueness of series with interval factors]
Counterexample \ref{cx_wo_needed_for_uniqueness} also  illustrates what can go wrong in Theorem \ref{thm_cseries} if one drops the condition that $\tob$ be well ordered.

Let $\lina$ be any linearization of $\axlh{\allitv{\tob}}$.  We claim that $\frac{ \Xi(\dsh_\lina \itva) } { \Xi(\pdsh_\lina \itva) } = 0$ for every nonempty interval $\itva$, since $ \frac{ \Xi(\dsh_{\axl{\dop \tob \sqcup \tob}} \itva) } { \Xi(\pdsh_{\axl{\dop \tob \sqcup \tob}} \itva ) }  = 0$ by Equation \eqref{eq_psi_dsh_equa_pdsh}, and  $\frac{ \Xi(\dsh_\lina \itva) } { \Xi(\pdsh_\lina \itva) } \mdom  \frac{ \Xi(\dsh_{\axl{\dop \tob \sqcup \tob}} \itva) } { \Xi(\pdsh_{\axl{\dop \tob \sqcup \tob}} \itva ) }$.  However $\frac{ \sfun(\dsh_\lina \itva) } { \sfun(\pdsh_\lina \itva) } = \dga \neq 0$ when $\itva = \tob$.  Therefore $\Xi|_{\axl{\lina}}$ and $\sfun|_{\axl{\lina}}$ are two complete lattice homomorphisms $\axl{\lina} \to \Sub_\dga$ satisfying \eqref{eq_intervalconstraints} for each $\itva \in \allitv{\tob}$, but each engenders a distinct set of interval factors.
\end{counterexample}

\subsection{Existence of CD homomorphisms}

One of the first questions  one must ask, concerning saecular homomorphisms, is \emph{when do they exist}?  
Existence of a CD homomorphism is nontrivial to check, in general.  However, Theorem \ref{thm_simplifiedcdfextensioncriterion} gives a simple criterion for the special case of bifiltrations on a modular algebraic lattice.  Theorem \ref{thm_woextensioninmodularuc} yields an even simpler criterion, conditioned on well-ordering.

\begin{theorem}[\cite{hp_cdcp}]
\label{thm_simplifiedcdfextensioncriterion}
If $ \tob$ and $\dop \tob$ are $\forall$-complete chains in a modular algebraic lattice $\cdla$, then $\dop \tob \cup \tob$  extends to a CD sublattice of $\lata$ iff 
	\begin{align}
	\textstyle
	\bigwedge_{\lelb \in B} (\lela \vee \lelb) = \lela \vee \bigwedge_{\lelb \in B} \lelb
	\label{eq_algebraic_condition}
	\end{align}
for each  $\lela \in \dop \tob \cup \tob$ and each set $B$ contained either in  $\dop \tob$ or in $\tob$.
\end{theorem}

\begin{counterexample}[Bifiltrations of a vector space that do not extend to a CD sublattice]
\label{cx_infinite_dimensional_bifiltrations_cd}
An important class of modular algebraic lattices are those of form $\lata = \Sub(V)$, where $V$ is a vector space.  Subspace lattices have a number of special properties, but even within this restricted class, one can find bifiltrations that do not extend to CD sublattices.  To illustrate, let $V$ be the space of finitely-support funections $\Z \to \fielda$, and let  $\lata = \Sub(v)$.  Let $X^m$ be the subspace of $V$ consisting of functions supported on $(-\infty, m]$. Let $\tob = \{X^m : m \in \Z \}$, and let  $\dop \tob = \{ 0, Y, V \}$, where $Y$ is the subspace generated by $ \{ \delta_m - \delta_n : m, n \in \Z, \; m \neq n \}$.  Here $\delta_m$ denotes the Kronecker delta.  If we take $B: = \tob$ and $\a = Y$, then the lefthand side of \eqref{eq_algebraic_condition} equals $V$, but the righthand side equals $\a = Y < V$.
\end{counterexample}

\begin{counterexample}
\label{cx_infinite_dimensional_bifiltrations_cd_generalized}
One can modify Counterexample \ref{cx_infinite_dimensional_bifiltrations_cd} to obtain a new pathology by defining $V = X^0$, $\tob = \{ X^m \cap X^0 : m \in \Z \}$, and $\dop \tob = \{ 0, Y \cap X^0, V \}$.  If we take $B: = \tob$ and $\a = Y \cap X^0$, then the lefthand side of \eqref{eq_algebraic_condition} equals $V$, but the righthand side equals $\a = Y \cap X^0 < V$.
\end{counterexample}

\begin{corollary}
\label{cor_inf_dim_vec_space_bifilt_dont_extend}
The following are equivalent for any vector space $V$:
	\begin{enumerate*}
	\item Every pair of chains $\dop \tob, \tob \su \Sub(V)$ extends to a CD sublattice $\dop \tob \sqcup \tob \su \latb$.
	\item The space $V$ is finite-dimensional.
	\end{enumerate*}
\end{corollary}
\begin{proof}
If $V$ is finite-dimensional then every bifiltration extends to a CD sublattice by Theorem \ref{thm_simplifiedcdfextensioncriterion}.  If  $V$ is not finite-dimensional, then it contains a subspace $U$ isomorphic to the space of finitely-supported functions $\Z \to \fielda$. In this case we may apply the construction of Counterexample \ref{cx_infinite_dimensional_bifiltrations_cd} to find a bifiltration of $U$ that does not extend to a CD sublattice; since $\Sub(U)$ is a $\mathring \forall$-complete sublattice of $\Sub(V)$, this bifiltration cannot extend to a CD sublattice of $\Sub(V)$.  The desired conclusion follows.
\end{proof}

\begin{reptheorem}{thm_woextensioninmodularuc}
Each pair of well ordered chains in a complete, modular, upper continuous lattice $\lata$ extends to a complete, completely distributive sublattice of $\lata$.
\end{reptheorem}

\begin{proof}
We defer the proof to \S\ref{sec_wo_cd_extension}.
\end{proof}

Theorem \ref{thm_saecularexistencecriteria} gives several existence criteria specific to the saecular CDF homomorphism.  
Here, for each $\tela \in \tob$, each lattice homomorphism $\cdf: \dla \to \Sub_\f$, and each subset $\seta \su \dla$, we set
	\begin{align*}
	\cdf[\seta]_\tela : = \{ \cdf(\ela)_\tela : \ela \in \seta \}
	\end{align*}

\begin{theorem}[Existence of the saecular CDF homomorphism]  
\label{thm_saecularexistencecriteria}
The saecular CD homomorphism exists if any one of the following criteria is satisfied.
	\begin{enumerate*}
	\item[] \emph{(Height criterion)} Lattice $\Sub(\f_\tela)$ has finite height, for each $\tela \in \tob$.			
	\item[] \emph{(Well-ordered criterion)} Lattice $\Sub (\f_\tela)$ is complete and upper continuous, and chains $\KK[\tob]_\tela$ and $\K[\tob]_\tela$ are well ordered for each $\tela \in \tob$.	
	\item[] \emph{(Algebraic criterion)} For each $\tela \in \tob$, lattice $\Sub (\f_\tela)$ is complete and algebraic, and 
	\begin{align*}
	\textstyle
	\bigwedge_{\lelb \in B} (\lela \vee \lelb) = \lela \vee \bigwedge_{\lelb \in B} \lelb
	\end{align*}
whenever  
	$
	\setx, \sety \in \{ \KK[ \tobc]_\tela , \K[\tobc]_\tela\}, 
	$
	$
	\lela \in \setx 
	$
	, and
	$
	B \su \sety.
	$
	\item[] \emph{(Split criterion)} $\ecata \in \{ R \catmod, \catmod R\}$, and $\dga$ splits as a direct sum of interval functors.	
	\end{enumerate*}
\end{theorem}

\begin{proof}
Suppose that $\dga$ is a direct sum $\bigoplus_\itva \dgb(\itva)$, where  $\itva$ runs over all  intervals in $\tob$  and  $\dgb(\itva) \in \allifun(\itva)$ for each  $\itva$.  Then there exists a  complete lattice homomorphism 
	$
	\cdf : \pslh \axlh(\dop \tob \sqcup \tob) \to \Sub_\dga
	$ 
such that $\cdf \{\dop \seta \sqcup \seta\} = \dgb(\seta - \dop \seta)$ for every down set of $\dop \tob \sqcup \tob$.  In fact  $\KK = \cdf \dop \cpph$ and $\K = \cdf \dop \cpph$ so, by uniqueness, $\FCD(\KK, \K)$ is the restriction of $\cdf$ to $\axlh^2(\dop \tob \sqcup \tob)$.  This establishes the split criterion.

For the algebraic and well-ordered criteria, recall that a subset  $\seta \su \Sub_\f$ extends to a complete CD sublattice iff $\{\ela_\tela : \ela \in \seta \}$ extends to a complete CD sublatice of $\Sub(\f_\tela)$ for all $\tela$.  The desired conclusion thus follows from Theorems  \ref{thm_simplifiedcdfextensioncriterion} and \ref{thm_woextensioninmodularuc}.  The height criterion implies the well ordered one.
\end{proof}

\begin{corollary} Every functor from a totally ordered set to the (full) subcategory of Artinian modules over a commutative ring $\ringa$ has a saecular CDF homomorphism.
\end{corollary}
\begin{proof}
This holds by the well ordered criterion.
\end{proof}

\subsection{Proof of Theorems \ref{thm_cdf_is_proper} and \ref{thm_saecular_factors_are_interval}}
\label{sec_prove_factors_are_intervals}

Recall that a saecular homomorphism is \emph{natural} at $\seta$ if 
	\begin{align}
	\dia |\seta|_\aa = |\seta \wedge \nu(\aa)|_\a = |\seta|_\a \wedge \KK(\aa)_\a
	\tag{\ref{eq_dicap}}
	\\
	\iia |\seta|_\a = |\seta \vee \nu(\a) |_\aa = |\seta|_\aa \vee \K(\a)_\aa
	\tag{\ref{eq_iicup}}
	\end{align}
A saecular homomorphism is natural if it is natural at every element of its domain.

\begin{lemma}
\label{lem_pushpull_basic}
If $\SBD(\dga)$ exists then it is natural on $\Im(\nu)$.  

Concretely, if (i) $L = \K$ and $\lelc \in \tob$, or (ii)  $L = \KK$ and $\lelc \in \dop\tob$,  then for each $\aa \dople \a$  
	\begin{align}
	\dia \L(\lelc)_\aa =  \L(\lelc)_\a \wedge  \KK(\aa)_\a 
	\label{eq_dicap_reduced}
	\\
	\iia \L(\lelc)_\a = \L(\lelc)_\aa \vee \K(\a)_\aa
	\label{eq_iicup_reduced}
	\end{align}
\end{lemma}
\begin{proof}
Equation  \eqref{eq_dicap_reduced} is essentially trivial when $\L = \KK$: in particular, $\aa \le \lelc$ implies $\L(\lelc)_\aa = \mora_\aa$, hence both sides of \eqref{eq_dicap_reduced} equal $\KK(\aa)_\a$, while $ \lelc < \a$ implies that both sides equal $\KK(\lelc)_\a$.  Equation \eqref{eq_dicap_reduced} is likewise trivial when $\L = \K$ and $\c \le \aa$, for in this case both sides equal 0.  On the other hand, when $\aa \le \c$ one has
	\begin{align*}
	\dia \K(\c)_\aa 
	= 
	\dia \iia \f(\a \dople \c)\ii  0 
	= 
	\K(\c)_\a \wedge \KK(\aa)_\a.
	\end{align*}
This establishes \eqref{eq_dicap_reduced} in full generality. The proof of \eqref{eq_iicup_reduced} is dual.
\end{proof}

For ease of reference, let us recall Theorem \ref{thm_cdf_is_proper}:

\begin{reptheorem}{thm_cdf_is_proper}
Suppose that $\mora$ admits a saecular BDH.
	\begin{enumerate*}
	\item $\SCD(\mora)$ is natural, if it exists.
	\item $\SBD(\mora)$ is natural if $\SCD(\mora)$ exists.
	\item $\SBD(\mora)$ is natural iff it is natural at each $\seta \in \Im(\dop \cpph) \cup \Im(\cpph)$.
	\end{enumerate*}
\end{reptheorem}

\begin{proof}[Proof of Theorem \ref{thm_cdf_is_proper}]
Suppose that $\SCD(\mora)$ exists.  We can regard the formulae on the righthand sides of \eqref{eq_dicap} and \eqref{eq_iicup} as functions of $\seta \in \axlh^2(\dop \tob \sqcup \tob)$; in fact, by complete distributivity, these formulae determine $\forall$-complete lattice homomorphisms from $\axlh^2(\dop \tob \sqcup \tob)$ to $[0,\KK(\aa)_\a] \su \Sub(\mora_\a)$ and $[\K(\a)_\aa, 1] \su \Sub(\mora_\aa)$, respectively.  Since $\axlh^2(\dop \tob \sqcup \tob)$ is the free CD lattice on $\dop \tob \sqcup \tob$, it therefore  suffices to verify \eqref{eq_dicap} and \eqref{eq_iicup} in the special case where $\seta \in \Im(\nu)$.  Lemma \ref{lem_pushpull_basic} supplies this fact.  This establishes assertion 1.

If $\SCD(\mora)$ exists, then it restricts to $\SBD(\mora) = \SCD(\mora)|_{\dxlh \axlh (\dop \tob \sqcup \tob)}$, so $\SBD(\mora)$ is s-natural as a special case.  This establishes assertion 2.

If \eqref{eq_dicap} and \eqref{eq_iicup} hold for each $\seta \in \Im(\dop \cpph) \cup \Im(\cpph)$ then they hold for every $\seta \in  \dxlh \axlh(\dop \tob \sqcup \tob)$, since every element of $\dxlh \axlh(\dop \tob \sqcup \tob)$  can be obtained from $\Im(\dop \cpph) \cup \Im(\cpph)$ via a finite number of meet and join operations, and both $\dia$ and $\iia$ preserve binary meets and joins, by Lemma \ref{lem_direct_image_complete_on_CD_complete_sublattice}.  This establishes assertion 3, and completes the proof.
\end{proof}

\begin{corollary}
\label{cor_kerofsub}
Let  $\setb \le \seta \in \axlh^2(\dop \tob \sqcup \tob)$ be given.  If $\SCD(\mora)$ exists then $\SCD\left (\xfrac{\seta}{\setb} \right)$ exists, and\footnote{The numerator of \eqref{eq_inducedkfiltrations} requires no parentheses, by the modular law. }

	\begin{align}
	\Ks_{ |\seta| /| \setb| }(\setc )
	=
		\frac
		{|\seta \wedge (\cpph \setc) \vee \setb| }
		{|\setb|}
	&&
	\KKs_{ |\seta| /| \setb| }( \dop \setc )
	=
		\frac
		{|\seta \wedge (\dop \cpph \dop \setc) \vee \setb| }
		{|\setb|}	
	\label{eq_inducedkfiltrations}
	\end{align}	
for each $\dop \setc \in \axlh^2(\dop \tob)$ and  $\setc \in \axlh^2(\tob)$.  More generally, 
	\begin{align}
	\sfun \setx = \xfrac{\seta \wedge \setx \vee \setb}{\setb}
	\label{eq_subquo_mapping_formula}
	\end{align}
for $\setx \in \axlh^2(\dop \tob \sqcup \tob)$, where $\sfun: = \SCD (\xfrac{\seta}{\setb})$ 
\end{corollary}
\begin{proof}
Theorem \ref{thm_cdf_is_proper} implies the last equality in each of the following two lines, for any $\aa \le \a \in \tob$.
	\begin{align}
	\Ks_{|\seta| /| \setb| }(\a)_\aa
	=
	\frac{|\seta|}{|\setb|} (\aa \le \a) \ii 0
	&
	=
	\frac{|\seta|(\aa \le \a)\ii |\setb|_\a}{|\setb|_\aa}
	=
	\frac{|\seta|_\aa \wedge \iia |\setb|_\a}{|\setb|_\aa}
	=
	\frac{|\seta|_\aa \wedge \K(\a)_\aa \vee |\setb|_\aa }{|\setb|_\aa}	
	\label{eq_subquo_formula_k}
	\\
	\KKs_{|\seta| /| \setb| }(\aa)_\a
	=	
	\frac{|\seta|}{|\setb|} (\aa \le \a) \di 1
	&
	=
	\frac{|\seta|(\aa \le \a)\di 1 \vee |\setb|_\a}{|\setb|_\a}
	=
	\frac{|\seta|_\a \wedge \KK(\aa)_\a  \vee |\setb|_\a}{|\setb|_\a}.	
	\label{eq_subquo_formula_kk}
	\end{align}
Therefore  \eqref{eq_inducedkfiltrations} holds when $\setc = \fem(\a)$ and $\dop \setc = \dop \fem(\aa)$.  If we regard the righthand side of \eqref{eq_subquo_mapping_formula} as a function  $h(\setx)$, then this function can be regarded as a composite $h = q_\bullet \circ p \circ \SCD( \mora)$, where $p(\soa)  = \soa \wedge |\seta|$ and $q$ is the quotient map $|\seta| \to |\seta| /|\setb|$.  Function $p$ restricts to a $\mathring \forall$-complete lattice endomorphism on $\Im(\SCD(\mora))$ by complete distributivity, and $q \di $ restricts to a $\forall$-complete lattice homomorphism on $\Im(\SCD(\mora))$  by Lemma \ref{lem_direct_image_complete_on_CD_complete_sublattice},  Thus $h$ is a $\forall$-complete lattice homomorphism.  Moreover, $h \circ \nu = \KKs_{|\seta| /| \setb| } \sqcup \Ks_{|\seta| /| \setb| }$, by \eqref{eq_subquo_formula_k} and \eqref{eq_subquo_formula_kk}.  Thus $h = \SCD\left (\xfrac{\seta}{\setb} \right)(\setx)$, as desired.
\end{proof}

\begin{corollary}
\label{cor_kerofsub_bd}
Suppose that $\sfun = \SBD(\dga)$ exists, and that $\sfun$ is natural at $\seta$ and $\setb$, where  $\setb \le \seta \in \dxlh \axlh(\dop \tob \sqcup \tob)$.
	\begin{enumerate*}
	\item  Diagram $|\seta| /|\setb|$ admits a saecular BD functor.
	\item  Equations \eqref{eq_inducedkfiltrations} hold for each $\dop \setc \in \axlh^2(\dop \tob)$ and  $\setc \in \axlh^2(\tob)$.
	\item Equation \eqref{eq_subquo_mapping_formula} holds for each $\setx \in \dxlh \axlh(\dop \tob \sqcup \tob)$.
	\end{enumerate*}
\end{corollary}
\begin{proof}
The proof is essentially that of Corollary \ref{cor_kerofsub}.  One shows that $h \circ \dop \cpph $ and $h \circ  \cpph $ are the free CD homomorphisms of $\KKs_{|\seta| /| \setb| }$ and $\Ks_{|\seta| /| \setb| }$, respectively, where $h = q\di \circ p \circ \SBD(\mora)$.  One likewise shows that $h$ is a lattice homomorphism, by invoking Lemma \ref{lem_direct_image_complete_on_CD_complete_sublattice}.  The desired conclusion then follows by uniqueness of free BD homomorphisms.
\end{proof}

\begin{theorem}
\label{thm_single2interval}
Every natural SBD functor $\sfun$ carries singletons to interval functors.  More precisely, 
	\begin{align*}
	\sfun \{\dop \setd \sqcup \setd \} \in \allifun( \setd - \dop \setd)
	\end{align*}
for each $\dop \setd \sqcup \setd  \in \axlh(\dop \tob \sqcup \tob)$.  The same conclusion holds if $\sfun$ is any SCD functor.
\end{theorem}
\begin{proof}
Fix $\seta, \setb \in \dxlh \axlh(\dop \tob \sqcup \tob)$ such that $\frac{\seta}{\setb} = \seta - \setb = \{\dop \setd \sqcup \setd \}$.  Then
	\begin{align*}
	\seta \cup \nu(\aa) \cap \setb  =  \begin{cases} \seta & \aa \notin \dop \setd \\ \setb & else \end{cases}
	&&
	\seta \cup \nu(\a) \cap \setb  =  \begin{cases} \seta & \a \notin  \setd \\ \setb & else \end{cases}	
	\end{align*}
Thus, by Corollary \ref{cor_kerofsub_bd},
	\begin{align*}
	\KKs_{|\seta|/|\setb|}(\aa)  =  \begin{cases} 1 & \aa \notin \dop \setd \\ 0 & else \end{cases}
	&&
	\Ks_{|\seta|/|\setb|}(\a)  =  \begin{cases} 1 & \a \notin  \setd \\ 0 & else \end{cases}
	\end{align*}
It follows that $|\seta|/|\setb|$ vanishes outside $\setd - \pi(\dop \setd)$.  Moreover, for each $\aa \dople \a \in \tob$, the map $\dga( \aa \dople \a)$ has kernel 0 and image 1, hence is an isomorphism.  The desired conclusion follows.
\end{proof}

\begin{proof}[Proof of Theorem \ref{thm_saecular_factors_are_interval}]
This is naturally equivalent to Theorem \ref{thm_single2interval}.
\end{proof}

\section{Saecular persistence in $\catgroup$}
\label{sec_saecular_group}

Most of the results presented for p-exact categories thus far have natural analogs for chain diagrams in $\catgroup$.  In general, however, these results require independent proofs.  A fundamental point of disconnect between the two domains is the following structural distinction:  in a p-exact category, every subobject of form $\K(\a)$ or $\KK(\a)$ is \emph{normal}, meaning that it is the kernel of some arrow $\dga \to \dgb$ in $[\tob, \ecata]$.  Such is not the case for  chain diagrams in $\catgroup$.  In fact, $\dgb$ is normal in $\dga$ iff $\dgb_\a \unlhd \dga_\a$ for all $\a \in \tob$.  This becomes relevant for the following reasons:
 	\begin{itemize}
	\item The lattice of normal subgroups, like the lattice of normal objects in any p-exact category, is modular; however, subgroup lattices are not normal in general.  Since we relied heavily on modularity to prove existence of the saecular BD homomorphism for p-exact categories,  we therefore need a new proof of existence of SBD homomorphisms in $\catgroup$.
	\item While $\catgroup$ has kernels and cokernels, a monomorphism $\dgb: H \into L$ satisfies $\dgb = \ker(\coker(\dgb))$ iff $\dgb$ is normal.  In particular, even if $H, L$ lie in the image of a lattice homomorphism, $\ker(\coker(\dgb))$ may not.
	\end{itemize}

In light of these and other hurdles, it is surprising that the machinery built for p-exact generally categories carries over to $\catgroup$ in any way.  The successful transfer of ideas from one domain to the other therefore speaks both to the power of p-exact categories as a source of inspiration, and to the power of order lattices to bridge disciplines.

\subsection{Definitions}

Let $\dga: \tob \to \catgroup$ be a chain diagram.  The saecular filtrations of $\mora$, as well as the  saecular CDH, CDIH, BDH, BDIH, and corresponding saecular functors, are defined exactly as in \S\ref{sec_saecular}.  As per convention, we  write $|\seta|$ for the value taken by one of these maps on a locally closed set $\seta$.

Let us say that a subdiagram $\dgb$  is \emph{normal}, written $\dgb \unlhd \dga$, if $\dgb_\a \unlhd \dga_\a$ for all $\a \in \tob$.

\begin{theorem}
\label{thm_saecular_functor_exiests_for_group_iff_images_normal}
The diagram $\dga$
	\begin{enumerate*}
	\item admits a saecular BDH and BDIH, though these maps may not be natural.
	\item admits a saecular CDH iff it admits a saecular CDIH.	
	\item admits a saecular BDF iff it admits a saecular BDH and $\KK(\aa) \unlhd \dga$ for all $\aa \in \dop \tob$. 
	\item admits a saecular CDF iff it admits a saecular CDH and $\KK(\aa) \unlhd \dga$ for all $\aa \in \dop \tob$. 			
	\end{enumerate*}
Moreover, the condition that $\KK(\aa) \unlhd \dga$ for all $\aa \in \dop \tob$ is equivalent to the condition that $\Im(\dga( a \le b))$ be a normal subgroup of $\dga_b$, for all $a \le b \in \tob$.	
\end{theorem}
\begin{proof}
Existence of a saecular BDH follows from Theorem \ref{thm_half_normal_admits_bd}.  The remainder of assertions 1 and 2 may be deduced by repeating the proof of Theorem \ref{thm_sacular_existence_iff}, verbatim.

Now consider assertion 4.  If $\dga$ admits a saecular CDF functor $\sfun$ then it admits a saecular CD homomorphism \emph{a fortiori}.  Moreover, every subobject $|\sob| \su |\soa|$ is the kernel of some arrow in $[\tob, \catgroup]$, by exactness of $\sfun$.   In particular $\K(\aa) = |\nu(\aa)| \unlhd \dga$ for all $\aa \in \dop \tob$, hence $
\Im(\dga( a \le b) ) \unlhd \dga_b$ for all $a \le b \in \tob$.  This establishes the first direction.  

For the converse, suppose that the saecular CD homomorphism $\cdf$ exists, and $\Im(\dga( a \le b) ) \unlhd \dga_b$ for all $a \le b \in \tob$.  Then $\KK(\aa) \unlhd \dga$ and $\K(\a) \unlhd \dga$ for all $\aa \in \dop \tob$ and $\a \in \tob$.   Thus the lattice  of normal subdiagrams $\latb \su \Sub_\dga$ contains the image of $\cdf$, since $\latb$ is a $\forall$-complete sublattice of $\dga$ (by Lemma \ref{lem_normal_subdiagrams_make_complete_subllattice_in_group}), and every subgroup $\cdf \seta$ can be expressed in terms of meets and joins of subgroups of form $\KK(\aa)$ and $\K(\a)$, by Theorem \ref{thm_tunconstruction}.  Thus $\cdf$ induces an exact functor, by Remark \ref{rmk_induced_exact_fuctor_for_groups}.  This completes the proof of assertion 4.  The proof of assertion 3 is similar.
\end{proof}

\begin{theorem}
\label{thm_cdf_is_proper_group}
Theorem \ref{thm_cdf_is_proper} holds for chain functors valued in $\catgroup$.  In particular, 
	\begin{enumerate*}
	\item $\SCD(\mora)$ is natural, if it exists.
	\item $\SBD(\mora)$ is natural if $\SCD(\mora)$ exists.
	\item $\SBD(\mora)$ is natural iff it is s-natural at each $\seta \in \Im(\dop \cpph) \cup \Im(\cpph)$.
	\end{enumerate*}
\end{theorem}
\begin{proof}
The proof of Theorem \ref{thm_cdf_is_proper} carries over without change.
\end{proof}

The \emph{saecular persistence diagram} of $\dga$ is defined exactly as in \S\ref{sec_spd_pexact}, i.e.
	\begin{align*}
	\spd : \allitv{\tob} \to [\tob, \ecata]
	&&
	\itva \mapsto  
	\sfun \{\itva \}
	\end{align*}
where $\sfun = \SBDIF(\mora)$.   Note that $\spd$ carries intervals to \emph{chain diagrams}, not to integers.  In particular, Theorem \ref{thm_saecular_factors_are_interval_group} states that the saecular persistence diagram carries each $\itva$ to an interval functor supported on $\itva$, \emph{assuming  s-naturality}.

\begin{theorem}
\label{thm_saecular_factors_are_interval_group}
Fix $\itva \in \allitv{\tob}$.  If $\SBD(\dga)$ is natural at $\dsh \itva$ and $\pdsh \itva$, then 
	\begin{align*}
	\frac
		{ | \dsh \itva | }
		{ | \pdsh \itva | }		
	\in 
	\allifun(\itva)
	\end{align*}
\end{theorem}
\begin{proof}
We will prove a stronger result below, namely  Theorem \ref{thm_single2interval_group}.
\end{proof}

We call $| \{\itva \}|$ the  \emph{$\itva^{th}$ saecular factor} of $\dga$.  Technically, the $\itva^{th}$ saecular factor is defined only up to unique isomorphism, as is $\SBDI(\dga)$ itself.  However, if desired one can define both $\SBDI(\dga)$ and $\spd(\itva)$ uniquely, as per Remark \ref{rmk_canonical_realizations_of_sfun}.    In this formulation, $\spd(\itva) = \xfrac{ \dsh \itva }{ \pdsh \itva }$ is the canonical realization of a quotient group via cosets, for each $\itva$.  We refer to this unique construction as the \emph{canonical realization} of $\spd$.

\subsection{The coset persistence diagram}

The saecular \emph{coset} persistence diagram of $\dga$ is the function 
	\begin{align*}
	\spdSET : \allitv{\tob} \to [\tob, \catsetp]
	&&
	\itva 
	\mapsto
	\Cst 
	\left (
	\xfrac{ \dsh \itva} {\pdsh \itva}
	\right )
	\end{align*}
where $\Cst(\soa / \sob)$ is the family of left cosets of $\sob$ in $\soa$.
Since every functor $\dga: \tob \to \catgroup$ admits a saecular BDIH, every $\mora$ admits a coset persistence diagram. 

We call $\spdSET(\itva)$ the \emph{$\itva^{th}$ saecular coset factor} of $\dga$.  When $\SBD(\mora)$ is natural, Theorem \ref{thm_saecular_coset_functor_carries_singletons_to_itvs_group} states that the saecular coset factors are interval diagrams in the category $\catsetp$ of pointed sets, if we regard  $\Cst(\soa/\sob)$ to be a pointed set with distinguished element $H$.

\begin{theorem}
\label{thm_saecular_coset_functor_carries_singletons_to_itvs_group}
If $\dsh \itva$ and $\pdsh \itva$ are natural, then 
$\spdSET( \itva ) \in \allifun(\itva)$.  
\end{theorem}
\begin{proof}
This is a special case of Theorem \ref{thm_single2interval_group}, which we prove below.
\end{proof}

\subsection{The normalized persistence diagram}

The \emph{normal closure} of a subgroup $\soa$ in a group $\oba$ is the smallest subgroup of $\oba$ containing $\soa$, denoted $\ncl_\oba(\soa)$.  The \emph{normal closure} of a subdiagram $\soa \le \sob$, denoted $\ncl_\sob (\soa)$ is the smallest normal subdiagram of $\sob$ containing $\soa$.  One can check that $\ncl_\sob (\soa)_\a = \ncl_{\sob_\a} (\soa_\a)$ for each $\a \in \tob$.

The  \emph{normalized} saecular persistence diagram of $\dga$ is the function 
	\begin{align*}
	\spdNCL : \allitv{\tob} \to [\tob, \catgroup]
	&&
	\itva 
	\mapsto
		\frac
			{| \dsh \itva  |}
			{\ncl_{|\dsh \itva |} ( |\pdsh \itva | ) } 
		=
		\Cok \left( 	\xfrac{ \dsh \itva} {\pdsh \itva} \right)
	\end{align*}
Since every functor $\dga: \tob \to \catgroup$ admits a saecular BDH, every $\mora$ admits a normalized persistence diagram. 

We call $\spdNCL(\itva)$ the \emph{$\itva^{th}$ normalized saecular factor} of $\dga$.  Theorem \ref{thm_saecular_coset_functor_carries_singletons_to_itvs_group_ncl} states that the saecular normalized factors are interval diagrams in the category $\catsetp$, when  $\SBD(\mora)$ is natural.

\begin{theorem}
\label{thm_saecular_coset_functor_carries_singletons_to_itvs_group_ncl}
If $\dsh \itva$ and $\pdsh \itva$ are natural, then 
$\spdNCL( \itva ) \in \allifun(\itva)$.  
\end{theorem}
\begin{proof}
This is a special case of Theorem \ref{thm_single2interval_group}, which we prove below.
\end{proof}

\begin{remark}
The coset PD, normalized PD, and saecular PD all coincide when the saecular BD functor exists -- however it need not exist, in general.
\end{remark}

\subsection{Non-normal subdiagrams}

In light of Theorem \ref{thm_saecular_functor_exiests_for_group_iff_images_normal}, it is natural to ask whether some version of the saecular functor can still be constructed, even when there exists a map $\dga(a \le b)$ whose image is not normal  in $\dga_b$.  Remarkably, this can be done; the only essential change is that we must work with cosets, rather than quotient groups.

The overall strategy can be summarized as follows.  We will define a relation $\mdomG$ on coset families, which is similar to the relation $\mdom$ which exists for subquotients in a p-exact category.  Like $\mdom$, the relation $\mdomG$ is an equivalence relation.  Also like $\mdom$, the statement $(\soa/\sob) \mdomG (\soa'/\sob')$ will turn out to imply  that $\soa/\sob$ is Noether isomorphic to $\soa'/\sob'$ -- not in the category of groups or even the category of sets equipped with group actions, but in the category $\catsetp$ of pointed sets.  

The relation $\mdomG$ will allow us to define quasi-regular induction for $\catgroup$ in much the same way that we defined quasi-regular induction for p-exact categories, \emph{cf.}\ diagram \eqref{eq_define_quasiregularlyin_induced}.  We will then complete the construction using a cohesive family of isomorphisms related to $\mdomG$; the result is  a close analogue of the saecular functor.

\subsection{Elements of group theory}

\newcommand{\ga}{\oba}
\newcommand{\gb}{\obb}
\newcommand{\sga}{H}
\newcommand{\sgb}{K}
\newcommand{\sgc}{L}
\newcommand{\nsga}{N}
\newcommand{\nsgb}{M}
\newcommand{\Set}{\mathtt{Set}}

Let us briefly recall some core facts in group theory.

\begin{lemma}
If $\sga, \sgb$ are subgroups of a group $\ga$ and $\sga \sgb = \sgb \sga$, then $\sga \sgb = \sga \vee \sgb$ in the subgroup lattice $\Sub_\oba$.
\end{lemma}

\begin{lemma}[{Modular law for subgroups that commute in gross }]
\label{lem_modularity_of_commuting_subgroups}
If $\sga, \sgb$ are subgroups such that $\sga \sgb = \sgb \sga$, then for each $\sga \le \sgc$ one has $(\sga \vee \sgb) \wedge \sgc = \sga \vee (\sgb \wedge \sgc)$.
\end{lemma}
\begin{proof}
See \cite[p13, remark following Theorem 11]{birkhoff1973lattice}.
\end{proof}

\begin{lemma}
\label{lem_normal_subgroup_lattice_is_complete}
The lattice $\latb$ of normal subgroups of $\ga$ is a $\forall$-complete sublattice of $\Sub(\ga)$.
\end{lemma}
\begin{proof}
It suffices to show that $\latb$ is closed under arbitrary meets and joins.  If $\seta$ is an arbitrary family of normal subgroups then $\bigvee \seta = \bigcup  \left \{ \bigvee \setb : \setb \su \seta, \# \setb < \infty  \right \}$.  It is therefore simple to verify that $g (\bigvee \seta) g^{-1} \su \bigvee \seta$ for all $g \in \ga$.  Thus $\bigvee \seta$ is normal, hence $\latb$ is closed under join.  On the other hand,  $x \in \bigwedge \seta \implies x \in \sga$ for each $ \sga \in \seta \implies g x g^{-1} \in \sga$ for each $\sga \in \seta \implies g x g^{-1} \in \bigwedge \seta$.  Thus $\bigwedge \seta$ is normal, hence $\latb$ is closed under meet.
\end{proof}

\begin{lemma}
\label{lem_normal_subdiagrams_make_complete_subllattice_in_group}
If $\dga: \tob \to \catgroup$ is a chain diagram, then the family of normal subdiagrams of $\dga$ is a $\forall$-complete sublattice $\latb \su \Sub_\dga$.
\end{lemma}
\begin{proof}
Follows from Lemma \ref{lem_normal_subgroup_lattice_is_complete}.
\end{proof}

\subsection{Regular and canonical induction}
\label{sec_induction_for_groups}

Let $\catps(\oba)$ denote the poset of nested pairs of subgroups in a group $\oba$, \emph{cf.}\ \S\ref{sec_induced_maps}.  For each doubly nested pair  $(\sob, \soa) \le (\sob', \soa')$ we refer to the map of cosets  $\soa/\sob \to \soa'/\sob', \; x \sob \mapsto x \sob'$ as a \emph{regularly induced map on left-cosets}.  Regularly induced maps on right cosets are defined similarly.  Lemma \ref{lem_functor_exists_from_nested_pair_poset_coset} follows directly from this definition.

\begin{lemma}[Coherence of regularly induced maps of cosets]
\label{lem_functor_exists_from_nested_pair_poset_coset}
Regularly induced maps of left cosets are closed under composition.  More precisely, 
there exists a functor $F: \catps(\oba) \to \catsetp$ sending each pair ${(\sob, \soa)}$ to the family of left cosets $\soa / \sob$ (equipped with distinguished point $e \sob = \sob$), and each relation $(\sob, \soa) \le (\sob', \soa')$ to the regularly induced map of left cosets $\soa/\sob \to \soa'/\sob'$.
\end{lemma}

The relation $\mdom$ was introduced for subquotients in a p-exact category in \S\ref{sec_induced_maps}.  We adapt it to the setting of groups by adding a constraint (membership in a lattice of ``commuting subgroups'') which ensures that Theorem \ref{thm_mutualdomination}  holds true in $\catgroup$ -- at least, true in the formal sense of Theorem \ref{thm_mutualdomination_groups}.  Concretely, given nested\footnote{We do not require that  $(\sob, \soa) \le (\sob', \soa')$.} pairs $(\sob, \soa)$ and  $(\sob', \soa')$, define
	\begin{align}
	\frac{\soa}{\sob} \; \mdomG  \; \frac{\soa'}{\sob'}
	&&
	\iff
	&&
	\begin{cases}
	M \vee N' = N \vee M'
	\quad
	\text{and}
	\quad
	M \wedge N' = M' \wedge N.
	\\
	(\soa \wedge \soa') \sob = \sob (\soa \wedge \soa') 
	\quad
	\text{and}
	\quad
	(\soa \wedge \soa') \sob' = \sob' (\soa \wedge \soa') 
	\\
	(\sob \vee \sob') \soa = \soa (\sob \vee \sob') 
	\quad
	\text{and}
	\quad
	(\sob \vee \sob') \soa' = \soa' (\sob \vee \sob') 			
	\end{cases}
	\label{eq_defofmutualdomination_groups}
	\end{align}

\begin{counterexample}
\label{cex_normal_mdom_implication}
It is generally false that $\frac \sob \soa \mdomG \frac {\sob'}{\soa'}$ implies $\soa \unlhd \sob \iff \soa' \unlhd \sob'$.    To illustrate, consider the case where $\ga = N \rtimes H$ is the semidirect product of a normal subgroup $N$ with non-normal subgroup $H$.  Then $\ga = NH = HN$ and $N \cap H = 0$, so it is easy to verify that $(\frac N 0) \mdomG ( \frac \ga  H)$.  However, $0$ is normal in $N$, while $H$ is not normal in $\oba$.
\end{counterexample}

\begin{theorem}[Induced bijections of cosets]  
\label{thm_mutualdomination_groups}
If $\soa, \soa'$ are subgroups (not necessarily normal) of a group $\ga$ and $\soa \soa' = \soa' \soa$, then the commutative square of subgroup inclusions on the lefthand side of \eqref{eq_noetheriso_for_groups} extends to a unique commutative diagram of regularly induced maps on left cosets \eqref{eq_noetheriso_for_groups}.  Moreover,  $\phi$ is bijective, and $\soa' \unlhd \soa \vee \soa' \implies \soa \wedge \soa' \unlhd \soa$.
\begin{equation}
    \begin{tikzcd}
    	\soa  \arrow[r, tail] & \soa \vee \soa'  \arrow[r, two heads] & (\soa \vee \soa') / \soa \\
		\soa \wedge \soa' \arrow[u, tail] \arrow[r, tail] & \soa' \arrow[u, tail] \arrow[r, two heads] & \soa' / (\soa \wedge \soa')  \arrow[u, "\phi" ', dashed]
    \end{tikzcd}
    \label{eq_noetheriso_for_groups}
\end{equation} 

If, in addition,  $\sob \le \soa$ and $\sob' \le \soa'$ are nested pairs of subgroups such that  $\soa / \sob \; \mdomG \; \soa' / \sob'$, then  there  exists a  commutative diagram 
\begin{equation}
    \begin{tikzcd}
    	& (\soa \vee \soa') / (\sob \vee \sob')    \\
		\soa/\sob   \arrow[ur] \arrow[rr, "\psi", dashed] && \soa' / \sob' \arrow[ul] \\
		& (\soa \wedge \soa') / (\sob \wedge \sob') \arrow[ul] \arrow[ur]
    \end{tikzcd}
    \label{eq_canonicallyinducediso_for_groups}
\end{equation} 
where each solid arrow is a regularly induced bijection of left cosets.  The bijection $\psi$ is a homomorphism of quotient groups iff $\sob \unlhd \soa$ and $\sob' \unlhd \soa'$.
\end{theorem}

\begin{proof}
The map $\phi$ is surjective because every element of $\soa \vee \soa' = \soa' \soa$ can be expressed in form $ab$ for some $a \in \soa'$ and $b \in \soa$.    If $x,y \in \soa'$ and $x \soa = y \soa$ then $y^{-1} x \in \soa$, hence $y^{-1} x \in \soa \wedge \soa'$ and therefore $y(\soa \wedge \soa') = y y^{-1} x(\soa \wedge \soa') = x (\soa \wedge \soa')$.  Thus $\phi$ is injective.  This establishes that $\phi$ is a bijection.  Uniqueness and commutativity of \eqref{eq_noetheriso_for_groups} is clear.   The implication $\soa' \unlhd \soa \vee \soa' \implies \soa \wedge \soa' \unlhd \soa$ is a standard property of normal subgroups.

Uniqueness and commutativity of the solid arrow diagram \eqref{eq_canonicallyinducediso_for_groups} are also clear.  Thus it suffices to check that each solid arrow in \eqref{eq_canonicallyinducediso_for_groups} is  bijective.  By the result we have just proved, it suffices to show that each solid arrow has form $\sga / (\sga \wedge \sgb) \to (\sga \vee \sgb) / \sgb$ for some pair of subgroups $\sga, \sgb$ such that $\sga \sgb = \sgb \sga$.

Consider the lower lefthand arrow in \eqref{eq_canonicallyinducediso_for_groups}.  In this case, let us set $\sga = \soa \wedge \soa'$ and $\sgb = \sob$.  We have  $\sga \sgb = \sgb \sga$  by the second of the three lines of conditions in \eqref{eq_defofmutualdomination_groups}.  The adapted modular law (Lemma \ref{lem_modularity_of_commuting_subgroups}) implies the second of the following equalities, and \eqref{eq_defofmutualdomination_groups} provides the third
	\begin{align*}
	\sga \vee \sgb
	=
	(\soa \wedge \soa') \vee \sob 
	= 
	\soa \wedge (\soa' \vee \sob)
	= 
	\soa \wedge (\soa \vee \sob')
	=
	\soa
	\end{align*}
Condition \eqref{eq_defofmutualdomination_groups} likewise implies the third equality in 
	$
	\sga \wedge \sgb 
	=
	(\soa \wedge \soa') \wedge \sob
	= 
	\sob \wedge \soa'
	=
	\sob \wedge \sob'.
	$
This shows that the lower left arrow in \eqref{eq_canonicallyinducediso_for_groups} has form $\sga / (\sga \wedge \sgb) \to (\sga \vee \sgb) / \sgb$, and is therefore bijective.  To verify that the upper left arrow is bijective, take $\sga = \soa$ and $\sgb = \sob \vee \sob'$.  Then  $\sga \wedge \sgb = (\sob \vee \sob') \wedge \soa = \sob \vee (\sob' \wedge \soa) = \sob \vee (\sob \wedge \soa') = \sob$ and $\sga \vee \sgb = (\sob \vee \sob') \vee \soa = \soa \vee \sob' = \soa \vee \soa'$.  Thus the two lefthand arrows are bijective.  The two righthand arrows are bijective by symmetry.

If $\sob \unlhd \soa$ and $\sob' \unlhd \soa'$ then all solid arrows in \eqref{eq_canonicallyinducediso_for_groups} are isomorphisms of quotient groups, hence $\psi$ is also a group isomorphism.  On the other hand, if $\sob$ is not normal in $\soa$ or $\sob'$ is not normal in $\soa'$ then either the domain or the codomain of $\psi$ fails to inherit the structure of a quotient group.
\end{proof}

\begin{definition}
We call $\psi$ the bijection \emph{canonically induced on left cosets}.  The bijection canonically induced on right cosets is defined similarly.
\end{definition}

\begin{remark}
Diagram \eqref{eq_noetheriso_for_groups} is a special case of the so-called \emph{product formula} for groups, which states that for \emph{any} pair of subgroups $\sga, \sgb$, the induced map on cosets $\sga / (\sga \cap \sgb) \to \sga \sgb / \sgb$ is bijective.
\end{remark}

\begin{remark}
\label{rmk_canonically_induced_bijection_being_group_hom_depends}
Counterexample \ref{cex_normal_mdom_implication} shows that the bijection $\psi$ defined by diagram \eqref{eq_canonicallyinducediso_for_groups} may not be a group homomorphism,  even when $\soa \unlhd \sob$.  However, it is  a group homomorphism when  $\sob \unlhd \soa$ and $\sob' \unlhd \soa'$.
\end{remark}

\begin{remark}
It is tempting to try to interpret canonical isomorphism in terms of pointed, transitive $G$-sets, however problems arise almost immediately.  To illustrate, suppose that $\morc: \ga \to \ga'$ is a group homomorphism and $\soa \le \sob \le \ga$ and $\soa' \le \sob' \le \ga'$ satisfy $\morc(\soa) \su \soa'$ and $\morc(\sob) \su \sob'$.  The regularly induced map of cosets $\sob/\soa \to \sob'/\soa', \; x \sob \mapsto \morc(x) \soa'$  can be interpreted as an arrow in the category $\catgsetpt$ where (i) objects are pairs $(\seta, \gb)$ such that $\gb$ is a group and $\seta$ is a pointed left $\gb$-set, and (ii) morphisms are pairs $( p: \seta \to \seta', q: \gb \to \gb')$ such that $q$ is a group homomorphism and $p( g \ela) = q(g) p(\ela)$ for all $\ela \in \seta$ and $g \in \gb$.  However, no such interpretation is available for the canonical isomorphism $\psi$ defined in \eqref{eq_canonicallyinducediso_for_groups}, in general.
\end{remark}

\begin{theorem}[Distributive cohesion for canonically induced bijections]
\label{thm_grandis_distributive_coherence_for_groups}
Let $\ga$ be a group and $\lata$ be a distributive sublattice of $\Sub(\oba)$ such that $\sga \sgb = \sgb \sga$ for all $\sga, \sgb \in \lata$.  Let $\setx$ be the family of  canonically induced bijections on left cosets $\soa/\sob \to \soa'/\sob'$, where $\soa, \soa', \sob, \sob' \in \lata$.  Then $\setx$ is closed under composition.
\end{theorem}
\begin{proof}
The proof provided by Grandis \cite[p 23-24, Theorem 1.2.7]{grandis12} carries over without alteration, if we use Theorem \ref{thm_mutualdomination_groups} in place of the diagram defined in \cite[p 20, equation (1.16) ]{grandis12}.
\end{proof}

\begin{theorem}[Compatibility of canonical bijections with regular induction on cosets ]
\label{thm_regular_cohesion_for_cosets}
Let $(\sob, \soa) \le (\sob', \soa')$ and $(\sod, \soc) \le (\sod', \soc')$ be doubly nested pairs.  Let $\oba_{{\soa}/{\sob}}$ and $\oba_{{\soa'}/{\sob'}}$ be objects in $\catsetp$ that realize the coset families $\soa/\sob$ and $\soa'/\sob'$, respectively, and define $\obb_{\soc/\sod}$ and $\obb_{\soc'/\sod'}$ similarly.

If $\soa/ \sob \; \mdomG  \; \soc / \sod$ and $\soa'/ \sob' \; \mdomG  \; \soc' / \sod'$, then the following diagram commutes, where horizontal arrows are regularly induced maps and vertical arrows are canonically induced isomorphisms.
\begin{equation}
    \begin{tikzcd}
    	\oba_{{\soa}/{\sob}} \ar[rr] 
		&& 
		\oba_{{\soa'}/{\sob'}} 
		\\
    	\obb_{{\soc}/{\sod}} \ar[rr, "" ]  \ar[u, dashed, "\cong"]
		&& 
		\obb_{{\soc'}/{\sod'}} \ar[u, dashed, "\cong"']
    \end{tikzcd}
    \label{eq_regular_cohesion}
\end{equation} 
\end{theorem}

\begin{proof}
The proof given in \cite[p 180, Theorem 4.3.7]{grandis12} adapts to this setting in a straightforward manner, substituting Lemma \ref{lem_functor_exists_from_nested_pair_poset} for \cite[p 112-113, Theorem 2.6.9]{grandis12} to guarantee coherence of regularly induced maps of cosets.
\end{proof}

A quasi-regular induction \emph{of left cosets} on a  sublattice $\lata \su \Sub_\oba$ is a commutative square of form
\begin{equation}
    \begin{tikzcd}
    	\frac{\soa}{\sob} \ar[rr, "\mora"] 
		&& 
		\frac{\soc}{\sod} 
		\\
    	\frac{\soa'}{\sob'} \ar[rr, "" ]  \ar[u, dashed, "\cong"]
		&& 
		\frac{\soc'}{\sod'} \ar[u, dashed, "\cong"']
    \end{tikzcd}
    \label{eq_define_quasiregularlyin_induced}
\end{equation} 
where the numerator and denominator of each fraction lie in $\lata$, vertical arrows are canonically induced bijections of cosets (in particular. we require that  $\soa/\sob \mdomG \soa'/ \sob'$ and  $\soc/\sod \mdomG \soc'/ \sod'$) and the lower horizontal arrow is regularly induced.  In this case we say that $\mora$ is quasi-regularly induced by $\lata$, and write $f : (\soa, \sob, \soc, \sod) \sim (\soa', \sob', \soc', \sod')$.

\begin{lemma}
\label{lem_quasiregularly_induced_maps_unique_for_cosets}
If  $\lata$ is a distributive sublattice of $\Sub(\oba)$ and $\soa \le \sob$ and $\soc \le \sod$ are nested pairs in $\lata$, then there exists at most one $\lata$-quasi-regularly induced map of cosets $\soa / \sob \to \soc/\sod$.
\end{lemma}
\begin{proof}
The proof is identical to that of Lemma \ref{lem_quasiregularly_induced_maps_unique}.
\end{proof}

Essentially all results in the present section, \S\ref{sec_induction_for_groups}, have natural generalizations for the functor category $[\cata, \catgroup]$, where $\cata$ is any category.  Joins and  meets in this category obtain object-wise, meaning that for each diagram $\oba: \cata \to \catgroup$ and each pair of  subdiagrams $\soa, \sob \in \Sub(\oba)$, one has  $(\soa \vee \sob)_\a = \soa_\a \vee \sob_\a$ and $(\soa \wedge \sob)_\a = \soa_\a \wedge \sob_\a$.  We write $\soa   \sob = \sob  \soa$ if $\soa_\a \sob_\a = \sob_\a \soa_a$ for each $\a$, and in this case we have a well-defined functor $\soa \cdot \sob = \soa \vee \sob : \cata \to \catgroup$ such that $(\soa \sob)_\a = \soa_\a \sob_\a$.  Regular induction on cosets has a natural (object-wise) definition in this setting, as does the relation $\mdomG$ and, therefore, the notion of $\lata$-quasi-regular induction.  Thus Lemmas \ref{lem_modularity_of_commuting_subgroups}, \ref{lem_functor_exists_from_nested_pair_poset_coset}, and \ref{lem_quasiregularly_induced_maps_unique_for_cosets} carry over essentially unaltered, as do Theorems \ref{thm_mutualdomination_groups} - \ref{thm_regular_cohesion_for_cosets}.  The details are left as an exercise to the reader.  To distinguish results regarding diagrams $\cata \to \catgroup$ from results regarding groups, we refer to the former as as \emph{$\cata$-shaped groups}.

\subsection{Coset functors}
\label{sec_coset_functors}

Let $\lata$ be the lattice of closed sets in a semitopological space $\setx$.  Let $\oba$ be either a group or a diagram $\cata \to \catgroup$, for some index category $\cata$.  Let $\cdf: \lata \to \Sub (\oba)$ be a bound-preserving lattice homomorphism from $\lata$ to the subobject lattice of $\oba$ and, as usual, write $|\seta|$ for $\cdf \seta$.  

We say that a functor $\sfun: \lclc( \lata) \to \catset$ is \emph{$\cdf$-quasi-regular} if it satisfies the following two conditions.

	\begin{enumerate}
	\item The functor $\sfun$ agrees with the homomorphism $\cdf$ on  $\lata \su \lclc(\lata)$.  More precisely, the restricted functor  $ \sfun|_{\lata}$ equals the functor $\funcforget \circ \cdf$ obtained by composing $\cdf$ either with the forgetful functor $\catgroup \to \catset$ (if $\oba$ is a group), or with the forgetful functor $[\cata, \catgroup] \to [\cata, \catset]$ (if $\oba$ is a $\cata$-shaped group).
	\item  For each regularly induced map $\mora: \soa / \sob \to \soc/\sod$ in the category $\lclc(\lata)$, there exists a quasi-regular induction of cosets on     $\Im(\cdf)$
\begin{equation}
    \begin{tikzcd}
    	\sfun \frac{\soa}{\sob} \ar[rr, "\sfun \mora"] 
		&& 
		\sfun \frac{\soc}{\sod} 
		\\
    	\xfrac{\soa}{\sob} \ar[rr, "" ]  \ar[u, dashed, "\cong"]
		&& 
		\xfrac{\soc}{\sod} \ar[u, dashed, "\cong"']
    \end{tikzcd}
    \label{eq_define_quasiregular_functor}
\end{equation} 
In particular, we require that $\xfrac{\soa}{\sob} \mdomG  \left (\sfun \frac{\soa}{\sob} \right)$ and $\xfrac{\soc}{\sod} \mdomG \left( \sfun \frac{\soc}{\sod} \right)$. Vertical arrows in \eqref{eq_define_quasiregular_functor} are the corresponding canonically induced bijections of cosets.

	\end{enumerate}

\begin{theorem}
\label{thm_existence_uniqueness_lclc_functors_into_group}
There exists a $\cdf$-quasi-regular functor $\sfun$ if $(\cdf \seta)(\cdf \setb) = (\cdf \setb) (\cdf \seta)$ for all $\seta, \setb \in \lata$

This functor is unique, in the following sense.  If $\sfun'$ is any other $\cdf$-quasi-regular functor, then there exists a unique natural isomorphism $\eta: \sfun \cong \sfun'$ such that $\eta|_{\lata}$ is identity.  Transformation $\eta$ assigns the canonically induced isomorphism $\eta_{\oba/\obb} : \sfun_{\oba/\obb} \cong \sfun'_{\oba/\obb}$ to each nested pair $\obb \le \oba \in \lata$.
\end{theorem}
\begin{proof}
The proof of Theorem \ref{thm_lclcextension} carries over with minimal alteration, using Lemmas \ref{lem_functor_exists_from_nested_pair_poset_coset} and \ref{lem_quasiregularly_induced_maps_unique_for_cosets} and Theorems \ref{thm_mutualdomination_groups} - \ref{thm_regular_cohesion_for_cosets} in place of Lemmas \ref{lem_functor_exists_from_nested_pair_poset} and \eqref{lem_quasiregularly_induced_maps_unique} and Theorems \ref{thm_mutualdomination} - \ref{thm_regular_cohesion}.
\end{proof}

As stated in Theorem \ref{thm_existence_uniqueness_lclc_functors_into_group},  the functor $\sfun$  is defined only up to unique isomorphism.  However, if $\lata$ is the family of closed sets in an Alexandrov topology, then for each locally closed set $X$ there exists a unique minimum nested pair $\setb_X \le \seta_X \in \lata$ such that $X = \seta_X - \setb_X$.

\begin{definition}
If $\lata$ is the lattice of closed sets in an Alexandrov topology, then the \emph{minimum quasi-regular functor induced by $\cdf$}  is the unique $\cdf$-quasi-regular functor such that $\sfun X = |\seta_X|/|\setb_X|$ for all locally closed $X$.  
\end{definition}

The minimum quasi-regular functor is important for the following reason:  given a map $\mora: X \to Y$ in $\lclc(\lata)$, there may exist two distinct realizations of the induced functor, $\sfun$ and $\sfun'$, and $\sfun \mora$ may be a group homomorphism while $\sfun'$ is not, \emph{cf.}\ Remark \ref{rmk_canonically_induced_bijection_being_group_hom_depends}.  The status of $\sfun \mora$ as a bona fide group homomorphism depends entirely on the status of $\sfun X$ and $\sfun Y$ as bona fide quotient groups, by Theorem \ref{thm_mutualdomination_groups}. Also  by Theorem \ref{thm_mutualdomination_groups}, the minimal functor $\sfun''$ carries $X$ to a bona fide quotient group iff there exists \emph{any} induced functor that does so.  Therefore the minimal realization carries $\mora$ to a bona fide group homomorphism iff there exists any induced functor that does so.

\subsection{Subsaecular series}

In yet another surprise, Theorem \ref{thm_cseries} carries over unchanged for groups:

\begin{theorem} 
\label{thm_cseries_group}
Suppose that $\tob$ is well ordered.  If $\lina$ is a linearization of   $\allitv{\tob}$, then  there exists a unique complete lattice homomorphism $\cdflin: \axl{\lina} \to \Sub_\f$ such that
	\begin{align}
		\frac
		{ \cdflin \dsh_\lina   (\itva) }
		{ \cdflin \pdsh_\lina (\itva) }		
	\in
	\allifun(\itva )
	\label{eq_subsaecular_coset}
	\end{align}
for each $\itva \in \allitv{\tob}$, where \eqref{eq_subsaecular_coset} is interpreted as diagram of coset families in $\catsetp$.  In particular, $\cdflin  = \SCD(\dga)|_{\axl{\lina}}$.
\end{theorem}
\begin{proof}
The proof of Theorem \ref{thm_cseries} carries over essentially unchanged.  Existence holds by Theorem \ref{thm_woextensioninmodularuc}.  One must invoke Lemma \ref{lem_modularity_of_commuting_subgroups} to justify use of the modular law in \eqref{eq_sss_modular}.
\end{proof}

\newcommand{\downa}{\axl{\dop \tob \sqcup \tob}}

\subsection{Proof of Theorems \ref{thm_saecular_coset_functor_carries_singletons_to_itvs_group} and \ref{thm_saecular_coset_functor_carries_singletons_to_itvs_group_ncl}}

\begin{theorem}
\label{thm_interval_coset_normality_conditions}
Let $\itva \in \allitv{\tob}$ be given.  Let $\seta:= \dsh \itva$ and $\seta', \setb, \setb'$ be down-closed sets such that $\seta/\setb = \seta'/\setb' = \{\itva\}$.
	\begin{enumerate*}
	\item $|\setb'| \unlhd |\seta'| \implies |\setb| \unlhd |\seta|$
	\item $|\setb| \unlhd |\seta| \iff |\setb|_\a \unlhd |\seta|_\a$ for some $\a \in \itva$
	\item If $\seta$ and $\setb$ are natural with respect to $\dga$, then $|\seta|_\a = |\setb|_\a$ for all $\a \notin \itva$; moreover, for each $\aa \dople \a$ such that $\pi(\aa), \a \in \itva$ and each $X \in \{ |\seta|, |\setb|, \ncl_{|\seta|}(|\setb|) \}$, one has
		\begin{align}
		|\seta|( \aa \dople \a) \di  X_\aa  = X_\a
		&&
		|\seta|( \aa \dople \a) \ii  X_\a   = X_\aa 
		\label{eq_pushpull_group}
		\end{align}
	\end{enumerate*}
\end{theorem}
\begin{proof}
Fix $\dop \setd \sqcup \setd$ such that $\itva = \setd - \dopiso(\dop \setd)$.  If $\a \in \tob$ and $\aa \in \dop \tob$ then 
	\begin{align*}
		\seta \wedge \nu(\aa)  = \dsh_{\downa} ((\pdsh \aa) \sqcup   \setd ) \su \setb    \iff  \aa \in \dop \setd	
		&&
		\seta \wedge \nu(\a)  = \dsh_{\downa} ( \dop \setd \sqcup (\pdsh \a)) \su \setb   \iff  \a \in \setd	
	\end{align*}
Suppose that $\seta$ is natural. Then for any  $\aa \dople \a$ and any $\setc \in \{ \seta, \setb\}$, one has 
	\begin{align}
	|\seta| (\aa \dople \a) \ii |\setc|_\a 
		= 
		|\seta|_\aa \wedge ( \iia |\setc|_\a )
		=
		|\seta|_\aa \wedge ( |\setc|_\aa \vee \K(\a)_\aa )		
		=
		| \seta \wedge ( \setc \vee \nu_\a )|_\aa
		\quad
		&
		\begin{cases}
		  =  |\setc|_\aa &  \a \in \setd \\
		  =  |\seta|_\aa &  \a \notin \setd .
		\end{cases}
	\label{eq_group_pullback_relations}
	\\	
	|\seta| (\aa \dople \a) \di |\setc|_\aa 
		= 
		\dia |\setc|_\aa 
		=
		|\setc|_\a \wedge \KK(\aa)_\a		
		=
		| \setc \wedge \nu_\aa |_\a
		\quad
		&
		\begin{cases}
		 \le    |\setb|_\a &  \aa \in \dop \setd \\
		  =  |\setc|_\a &  \aa \notin \dop \setd .
		\end{cases}	
	\label{eq_group_pushforward_relations}
	\end{align}
Fix $\aa \dople \a$ such that $\pi(\aa), \a \in \itva$; equivalently, such that  $\aa \notin \dop \setd$ and $\a \in \setd$.  Then \eqref{eq_pushpull_group} holds for each $X \in \{ |\seta|, |\setb|\}$, by \eqref{eq_group_pullback_relations} and \eqref{eq_group_pushforward_relations}.  It follows, by Lemma \ref{lem_reg_induction_of_cosets_plus_normality}, that $|\seta| (\aa \dople \a) \di$ and $|\seta| (\aa \dople \a) \ii$ induce mutually inverse bijections between the family of normal subgroups of $|\seta|_\aa$ containing $|\setb|_\aa$ and the family of normal subgroups of $|\seta|_\a$ containing $|\setb|_\a$.  Thus  \eqref{eq_pushpull_group} holds for  $X = \ncl_{|\seta|}(|\setb|)$.
	
Now suppose that $\a = \dopiso(\aa)$, so that $|\seta|(\aa \dople \a)$ is the identity map.  If $\a \notin \setd$, then setting $\setc = \setb$ in \eqref{eq_group_pullback_relations} shows that $|\seta|_\a = |\setb|_\a$.
On the other hand, if $\aa \in \dop \setd$, then setting $\setc = \seta$ in \eqref{eq_group_pushforward_relations} shows $|\seta|_\aa \le |\setb|_\aa \le |\seta|_\aa$.  This completes the proof of assertion 3.

Assertion 2 follows from assertion 3, since $|\setb|_\a \unlhd |\seta|_\a \iff |\setb|_\a = \ncl_{|\seta|_\a}(|\setb|_\a)$.  It can be shown via basic set operations that $(\setb, \seta) \le (\setb', \seta')$ is a doubly nested pair, and $\setb = \setb' \cap \seta$.  Thus $|\setb'| \unlhd |\seta'|$ implies $|\setb| = |\setb'| \cap |\seta| \unlhd |\seta|$.  This establishes assertion 1, and completes the proof.
\end{proof}

\begin{theorem}
\label{thm_single2interval_group}
If $\itva \in \allitv{\tob}$ and both $\dsh \itva$ and $\pdsh \itva$ are natural with respect to $\dga$, then 
	\begin{align*}
	|\{ \itva\}| \in \allifun( \itva )
	\end{align*}
and
	\begin{align*}
		\frac
			{| \dsh \itva  |}
			{\ncl_{|\dsh \itva |} ( |\pdsh \itva | ) } \in \allifun( \itva )
	\end{align*}
\end{theorem}
\begin{proof}
Follows from assertion 3 of Theorem \ref{thm_interval_coset_normality_conditions}.
\end{proof}

\begin{lemma}
\label{lem_reg_induction_of_cosets_plus_normality}
Let $\phi: \ga \to \ga'$ be a group homomorphism with kernel $K$, and let $\sob, \soa$ be subgroups of $\ga$ such that $K \cap \soa \le \sob \le \soa$.  We do not require $\soa$ or $\sob$ to be normal.  Then there exist well-defined bijections of left cosets  
\begin{equation}
    \begin{tikzcd}
    	{ \frac {\soa} {\sob} } \ar[rr, "\cong"]
		&&
		\frac {\soa  K } {\sob K }   \ar[ rr, "\phi" ] 
		&&
		\frac{\phi(\soa)}{\phi(\sob)}
		&&
		x \sob \ar[r, mapsto ] 
		&  
		{x(\sob K)} \ar[r, mapsto] 
		&
		{ [\phi(x) \phi(\sob)] }
    \end{tikzcd}
    \label{eq_coset_iso_induced_by_phi}
\end{equation} 
Moreover,  $\sob \unlhd \soa \iff \sob K \vartriangleleft \soa K \iff \phi(\sob) \vartriangleleft \phi(\soa) $.
\end{lemma}
\begin{proof}
We claim that that the lefthand map in \eqref{eq_coset_iso_induced_by_phi}  is a regularly induced bijection, \emph{cf.}\ diagram \eqref{eq_noetheriso_for_groups}.  To prove this, it suffices to establish commutation, i.e.\ $\soa (\sob K) = (\sob K) \soa$, as well as  $\soa \wedge (\sob K) = \sob$ and $\soa \vee (\sob K) = \soa K$.  The first and third of these statements are simple exercises.  The second follows from the observation that $m \in \sob$, $k \in K$, $mk \in \soa \implies k \in \soa \implies k \in \soa \wedge K \le \sob \implies k \in \sob \implies mk \in \sob$.

Since $\phi^{-1}( \phi(\sob K)) = \sob K$, we have $\phi(a) \phi(\sob K) = \phi(b) \phi(\sob K) \implies \phi(a^{-1} b)  \in \phi(\sob K) \implies a^{-1}b \in \sob K \implies a (\sob K) = b (\sob K)$.  Thus the map $x (\sob K)  \mapsto \phi(x) \phi(\sob K)$ is well-defined and 1-1.  It is clearly surjective.  If $a (\sob K) a^{-1} \su \sob K$ then $ \phi(a)  \phi(\sob K) \phi( a^{-1}) \su  \phi(\sob K)$, hence $\sob K \unlhd \soa K$ implies $\phi(\sob K) \unlhd \phi(\soa K)$.  The converse may be deduced from the observation that $\phi^{-1}( \phi(\sob K)) = \sob K$.

Finally, we have seen that $\sob = \soa \wedge (\sob K)$, hence $\sob K \unlhd \soa K \implies \sob \unlhd \soa$.  Conversely, if $\sob \unlhd \soa$ then for each $k \in K$ and $n \in \soa$ we have $kn ( \sob K) = k(K \sob) n = K \soa n = \soa K n = \soa K kn = K \soa k n$, hence $\sob K \unlhd \soa K$.
\end{proof}

\subsection{Half-normal bifiltrations as genera of distributive sublattices}

\newcommand{\chainh}{ H}
\newcommand{\chainn}{ N}
\newcommand{\subgrouph}{H}
\newcommand{\subgroupn}{N}

Let $\lata$ be the lattice of subgroups of a group $G$.  Let $\chainh: \tob \to \lata$ and $\chainn: \toc \to \lata$ be lattice homomorphisms indexed by  totally ordered sets, where each  $\chainn_i$ is a normal for all $i$.  We will prove Theorems \ref{thm_halfnormalbifiltration} and \ref{thm_half_normal_admits_bd} after establishing some supporting results.

\begin{theorem}	
\label{thm_halfnormalbifiltration}
The set $\Im(\chainh) \cup \Im(\chainn)$ extends to a distributive sublattice of $\lata$.
\end{theorem}

\begin{theorem}
\label{thm_half_normal_admits_bd}
The copairing $\chainh \sqcup \chainn$ (in the category of posets $\catps$) admits a free BD homomorphism.
\end{theorem}

\begin{remark}
Theorem \ref{thm_half_normal_admits_bd} relates strongly to the free modular lattice on two chains.  This lattice was  identified by Birkhoff \cite{birkhoff1973lattice} in his exploration of order lattices and their applications in algebra.  We do not know whether this result was known to him.  However several results in \cite{birkhoff1973lattice} seem perfectly formulated to enable the proof of Theorem \ref{thm_half_normal_admits_bd}.
\end{remark}

Our proof of Theorem \ref{thm_halfnormalbifiltration} closely follows that of \cite{birkhoff1973lattice} for modular lattices.  The lattice of subgroups $\lata$ is not modular in general, so we must work without several key assumptions, e.g.\ duality.  To begin, let $u^i_j = \chainh_i \wedge \chainn_j$ and $v^i_j = \chainh_i \vee \chainn_j$. Let $X$ be the set of all subgroups expressible as as joins of finitely many $u^i_j$.  

\begin{lemma}
Each element of  $X$ can be written in form $u^{i_1}_{j_1} \vee \cdots \vee u^{i_m}_{j_m}$ for some  $i_1 > \cdots > i_m  $ and $j_1 < \cdots < j_m$.
\end{lemma}
\begin{proof}
Given a finite subset $\seta \su \lata$, the join $\bigvee \seta$ equals $\bigvee \max(\seta)$.  Note, in particular, that $\max(\seta)$ is an antichain.
\end{proof}

\begin{lemma}
\label{lem:eltsofXcommute}
For any $x,y \in X$, one has $x \vee y = xy = yx$.  
\end{lemma}
\begin{proof}
It suffices to show that $u^i_j  u^k_l = u^k_l u^i_j$ for any $i,j,k,l$.  Without loss of generality, $i \le k$, so $u^i_j, u^k_l \le \chainh_k$.  The desired conclusion follows, since $u^k_l$ is normal in $\chainh_k$.
\end{proof}

\begin{lemma}
\label{lem:Xobeysmodularlaw}
The elements of $X$ obey the modular law, in the sense that $x,y,z \in X$ and $x \le z$ implies $x \vee (y\wedge z) = (x \vee y) \wedge z$.
\end{lemma}
\begin{proof}
This result follows from Lemma \ref{lem:eltsofXcommute} .  
\end{proof}

\begin{lemma}
\label{lem:transmuteuandv}
If $i_1 \ge \cdots \ge i_r$ and $j_1 \le \cdots \le j_r$, then
	\begin{align}
	\label{eq_wedgejoinrewrite}
	(\chainh_{i_1} \wedge \chainn_{j_1}) \vee \cdots \vee (\chainh_{i_r} \wedge \chainn_{j_r})
	& =
	\chainh_{i_1} \wedge (\chainn_{j_1} \vee \chainh_{i_2}) \wedge \cdots \wedge (\chainn_{j_{r-1}} \vee \chainh_{i_r}) \wedge \chainn_{j_r}
	\\
	\label{eq_joinwedgerewrite}	
	(\chainn_{j_1} \vee \chainh_{i_1}) \wedge \cdots \wedge (\chainn_{j_r} \vee \chainh_{i_r})
	& =
	\chainn_{j_1} \vee ( \chainh_{i_1} \wedge \chainn_{j_2}) \vee \cdots \vee (\chainh_{i_{r-1}} \wedge \chainn_{j_r}) \vee \chainh_{i_r}	
	\end{align}
\end{lemma}
\begin{proof}
We proceed by induction on $r$.  Two applications of the modular law (Lemma \ref{lem:Xobeysmodularlaw}) show that the lefthand sides of \eqref{eq_wedgejoinrewrite} and \eqref{eq_joinwedgerewrite} can be rewritten as 
	\begin{align}
		\label{eq_wedgejoinrewrite_halfway}
	\chainh_{i_1} \wedge [ \chainn_{j_1} \vee (\chainh_{i_2} \wedge \chainn_{j_2}) \vee \cdots \vee (\chainh_{i_{r-1}} \wedge \chainn_{j_{r-1}}) \vee \chainh_{i_r} ] \wedge \chainn_{j_r}
	\\
		\label{eq_joinwedgerewrite_halfway}	
	\chainn_{j_1} \vee [ \chainh_{i_1} \wedge (\chainn_{j_2} \vee \chainh_{i_2}) \wedge \cdots \wedge (\chainn_{j_{r-1}} \vee \chainh_{i_{r-1}}) \wedge \chainn_{j_r}  ] \vee \chainh_{i_r}	
	\end{align}
respectively.  We may apply the inductive hypothesis for \eqref{eq_joinwedgerewrite} to rewrite the expression in the square brackets of \eqref{eq_wedgejoinrewrite_halfway}, such that \eqref{eq_wedgejoinrewrite_halfway} becomes the righthand side of \eqref{eq_wedgejoinrewrite}.  A similar substitution transforms \eqref{eq_joinwedgerewrite_halfway} into the righthand side of \eqref{eq_joinwedgerewrite}.  The desired conclusion follows.
\end{proof}

\begin{lemma}
\label{lem:Xisadistributivesublattice}
Set $X$ is a distributive sublattice.
\end{lemma}
\begin{proof}
Lemma \ref{lem:transmuteuandv} shows that every join of elements of form $u^i_j$ is a meet of elements of form $v^i_j$, and vice versa.  Thus $X$ is closed under both meet and join, hence a sublattice.  It is modular by Lemma \ref{lem:Xobeysmodularlaw}.  Since $X$ is generated by two chains, it is therefore distributive, by \cite[p66, Theorem 9 and exercise 12]{birkhoff1973lattice}.
\end{proof}

\begin{proof}[Proof of Theorem \ref{thm_halfnormalbifiltration}]
Clearly $X$ contains $\chainh(\toa) \cup \chainn(\tob)$ and is contained in the sublattice generated by $\chainh(\toa) \cup \chainn(\tob)$.  Thus $X$ is the sublattice generated by $\chainh(\toa) \cup \chainn(\tob)$.  The desired conclusion follows from Lemma \ref{lem:Xisadistributivesublattice}.
\end{proof}

\begin{proof}[Proof of Theorem \ref{thm_half_normal_admits_bd}]
Let $\bar \chainh: \axlh^2( \tob) \to \lata$ and $\bar \chainn: \axlh^2(\toc) \to \lata$ be the free CD homomorphisms of $\chainh$ and $\chainn$, respectively.  By Lemma \ref{lem_normal_subgroup_lattice_is_complete}, every element of $\Im(\bar \chainn)$ is normal.  Thus $\bar \chainh \sqcup \bar \chainn$ factors through a distributive sublattice of $\lata$, by Theorem \ref{thm_halfnormalbifiltration}.
\end{proof}

\section{Persistent Homology}
\label{chp:hompersistence}

Let $\ecata$ be a p-exact category and $\catcx(\ecata)$ be the p-exact category of chain complexes in $\ecata$. Let $\tob$ be a totally ordered set with top element $1$, and $\fcxa: \tob \to \ecatb$ be a nested family of subcomplexes, i.e., a chain functor that carries each arrow to an inclusion.  Fix a dimension $\inta$, and let $\dga = \funh_\inta \circ \fcxa : \tob \to \ecata$ be the functor obtained from $\fcxa$ by post composition with the $\inta$th homology functor.  Then we have filtrations $\Ks: \tob \to \Sub( \funz_\inta \fcxa_1 )$ and $\KKs: \dop \tob \to  \Sub( \funz_\inta \fcxa_1 )$ such that 
	$
	\KKs(\a) = \funz_\inta \fcxa_\a
	$
	and
	$
	\Ks(\a) =  \funbd_\inta \fcxa_\a 
	$
where $\funbd_\inta \fcxa_\a$ and $\funz_\inta \fcxa_\a$ denote the $\inta$-dimensional boundary and cycle subobjects, respectively, of $\fcxa_\a$.  Then
	\begin{align}
	\KK = q \circ \KKs
	&&
	\K  = q \circ \Ks
	\label{eq_ph_qcomposition}
	\end{align}
where $q$ is the poset map $\Sub(\funz_\inta \fcxa_1) \to \Sub_\f$ such that 
	\begin{align*}
	q(X)_\a & = 
		\frac
		{\funz_\inta \fcxa_\a \wedge X  \vee \funbd_\inta \fcxa_\a }
		{\funbd_\inta \fcxa_\a}
	\end{align*}
for each $\a \in \tob$ and $\setx \in \axlh^2(\dop \tob \sqcup \tob)$.
The numerator requires no parentheses, by the modular law (Proposition \ref{prop_pexact_modular_principles}). 

\begin{theorem} 
\label{thm_homologicalconnection}
If $\KKs$ and $\Ks$ admit a joint CD homomorphism $\sfun = \FCD(\KKs \sqcup \Ks)$, then  $\KK$ and $\K$ admit a joint CD homomorphism $\sfun_\dga : = \FCD(\KK \sqcup \K)$, and 
	\begin{align}
	\sfun_\dga = (q \circ \sfun).  
	\label{eq_ph_cdf_formula}
	\end{align}
In this case there exists a second-type Noether isomorphism
	\begin{align*}
	\sfun_\dga \{\dop \seta \sqcup \seta \}_\a \cong \sfun \{\dop \seta \sqcup \seta \}
	\end{align*}
for each $\dop \seta \sqcup \seta \in \axlh^2(\dop \tob \sqcup \tob)$ and each $\a \in \seta - \dopiso( \dop \seta) $.  
\end{theorem}
\begin{proof}
Suppose that $\KKs$ and $\Ks$ admit a joint CD homomorphism $\sfun$.  Then $\Im(\sfun)$ is a $\forall$-complete sublattice of $\Sub(\fcxa_1)$ containing $\funz_\inta \fcxa_\a $ and $\funbd_\inta \fcxa_\a$ for all $\a$.  Consequently, the rule $X \mapsto \funz_\inta \fcxa_\a \wedge X  \vee \funbd_\inta \fcxa_\a$ determines a $\mathring \forall$-complete lattice homomorphism $\Im(\sfun) \to \Im(\sfun)$. Thus $q$ is a $\forall$-complete lattice homomorphism.  So, too, is $(q \circ \sfun)$. The righthand identity in \eqref{eq_ph_cdf_formula} follows, by \eqref{eq_ph_qcomposition} and definition of free CD homomorphisms.

Now fix a nonempty interval $\itva = \seta - \dopiso(\dop \seta)$ and choose $\a \in \itva$. Since $\cpph(\dsh \dop \seta) \le \nu( \dopiso^{-1}\a)$, 
one has 
	\begin{align*}
	\{\dop \seta \sqcup \seta \}
	=
	\frac
		{\dsh( \dop \seta \sqcup \seta)}
		{\pdsh( \dop \seta \sqcup \seta)}		
	=
	\frac
		{\nu_\a \vee \dsh( \dop \seta \sqcup \seta)}
		{\nu_\a \vee \pdsh( \dop \seta \sqcup \seta)}			
	=
	\frac
		{(\nu_\a \vee \dsh( \dop \seta \sqcup \seta))  / \nu_\a}
		{(\nu_\a \vee \pdsh( \dop \seta \sqcup \seta)) / \nu_\a}					
	\end{align*}
It follows, by exactness of $\sfun$, that there exist (coherent) Noether isomorphisms 
	\begin{align*}
	\sfun \{\dop \seta \sqcup \seta \}
	\cong
	\frac
		{\sfun \dsh( \dop \seta \sqcup \seta)}
		{\sfun \pdsh( \dop \seta \sqcup \seta)}		
	\cong
	\frac
		{\sfun \nu_\a \vee \dsh( \dop \seta \sqcup \seta)}
		{\sfun \nu_\a \vee \pdsh( \dop \seta \sqcup \seta)}			
	\cong
	\frac
		{\sfun (
			\frac
				{	\nu_\a \vee \dsh( \dop \seta \sqcup \seta) 	}
    			{ 	\nu_\a	}
			)
		}
		{\sfun (
			\frac
				{	\nu_\a \vee \pdsh( \dop \seta \sqcup \seta) 	}
    			{ 	\nu_\a	}
			)
		}
	=
	\frac
		{\sfun_\dga \dsh( \dop \seta \sqcup \seta)_\a}
		{\sfun_\dga \pdsh( \dop \seta \sqcup \seta)_\a}
	=
	\sfun_\dga \{\dop \seta \sqcup \seta \}_\a								
	\end{align*}
The desired conclusion follows.
\end{proof}

\section{Generalized persistence diagrams}
\label{chp_pepersistence}

For much of its history, work in persistent homology has focused on homology with field coefficients. Patel \cite{patel2018generalized} proposed several notions of persistence diagram suitable for any Krull-Schmidt category and, in particular any abelian one, subject to certain tameness requirements.

This approach characterizes the persistence diagram as the M\"obius inverse of a so-called \emph{rank function}.  This is consistent with the general framework in which classical $\fielda$-linear persistence diagrams were first defined \ \cite{ELZTopological02}.  The key insight of Patel, however, was a suitable definition of M\"obius inverse when the the rank function takes values not in integers, but in objects of a Krull-Schmidt category.  We briefly review this definition, then relate it with the saecular approach. 

\subsection{Background}
\label{sec_amitbackground}

Let $\ecata$ be an essentially small symmetric monoidal category with identity object $e$.  Write $[a]$ for the isomorphism class of an object $a$, and let $\amitj(\ecata)$ be the family of all isomorphism classes in $\ecata$, regarded as a commutative monoid under the the binary monoidal operation $[\lela] + [\lelb] = [\lela \square \lelb]$.  

Let $\amita(\ecata)$ be the group completion of $\amitj(\ecata)$.  If, in addition, $\ecata$ is abelian, then one may define a quotient monoid $\amitb(\ecata)= \amita(\ecata)/\sim$, where $[\lela] \sim [\lelb] + [\lelc]$ iff there exists a (not necessarily split) short exact sequence  $0 \to \lelb \to \lela \to \lelc \to 0$.  It is a fact that $\amita(\ecata)$ and $\amitb(\ecata)$ admit translation invariant partial orders compatible with subobject inclusion.  

Let $\pitvc$ be the poset whose elements are half-closed intervals $[\lela, \lelb) \su \R$, ordered such that $[a,b) \le [c,d)$ iff $[c,d) \su [a,b)$.  Let $\seta$ be a finite set of form
	\begin{align}
	\posts = \{\post_0 \lneq \cd \lneq \post_\intmaxa = \infty\} \su \R \cup \{\infty\},
	\label{eq_postset}
	\end{align}
If $\poga$ is an abelian group equipped with a translation invariant ordering of elements, then a  function
	$
	X: \pitvc \to \poga
	$
is defined to be \emph{$\posts$-finite} if 
	$
	\supp(X) \su  \{[\post_\inta, \post_\intb) : 0 \le \inta, \intb \le \intmaxa \}
	$. 
We call $X$ a \emph{generalized persistence diagram} if it is $S$-constructible for some $S$.	
	
A chain functor
	$
	\dga: \R \to \ecata
	$
is  \emph{$\posts$-constructible} if

\begin{itemize} 
\item   $\dga(p \le q)$ is the identity map on object $e$  whenever  $p \le q < \ela_1$, and 
\label{constructibleconda}

\item   $\dga(p \le q)$ is an isomorphism whenever $\ela_\inta \le p \le q < \ela_{\inta +1}$.
\label{constructiblecondb}
\end{itemize}
When call $\dga$  \emph{constructible} if it is $S$-constructible for some $S$.

For any $\seta$-constructible functor $\dga$ there exist maps
	\begin{align*}
	\prfa{\dga}:
	& \pitvc \to \amita(\ecata)
	&
	\prfa{\dga}(\itva) 
	&= 
	\lim_{t \to \post - } [\Im \dga(p \lneq \ela_\inta - \delta) ]	
	\\
	\prfb{\dga}:
	& \pitvc \to \amitb(\ecata)	
	&
	\prfb{\dga}(\itva) 
	&= q \com \prfb{\dga}(\itva)
	\end{align*}
where $q$ is the quotient map $\amita(\ecata) \to \amitb(\ecata)$.  The type-$\amita$ persistence diagram of $\dga$ is the M\"obius inversion of $\prfa{\dga}$, denoted $\pdta{\dga}: \pitvc \to \amita(\ecata)$.  If $\ecata$ is Abelian, then the type-$\amitb$ persistence diagram of $\dga$ is the M\"obius inversion of $\prfb{\dga}$, denoted $\pdtb{\dga}: \pitvc \to \amitb(\ecata)$.  

\begin{theorem}[Patel, \cite{patel2018generalized}] 
The type-$\amita$ and $\amitb$ persistence diagrams are $\posts$-finite when $\dga$ is $\posts$-constructible.
\end{theorem}

\noindent There exist several stability results for diagrams of this type, \emph{cf.}\ \cite{patel2018generalized}.

\subsection{Saecular enumeration theorem for generalized persistence}

Let $\ecata$ be an essentially small abelian category, regarded as a symmetric monoidal category where the monoidal operation is the biproduct, i.e. $\oba \square \obb = \oba \oplus \obb$.  The unit object is then $e = 0$.  Let  $\posts$ be a set of form \eqref{eq_postset}, and  
	\begin{align*}
	\f: \R \to \ecata
	\end{align*} 
be an $\posts$-constructible functor.  An object $\oba$ in category $\ecata$ has \emph{finite height} if $\Sub_\oba$ has finite height as an order lattice.  In this case each  maximal-with-respect-to-inclusion chain of subobjects $\lina = \{0 = \lina_0 <  \cdots < \lina_k  = \oba \} \su \Sub_\oba$ engenders a set function 
	\begin{align}
	\jhv{\lina} :  \smic \to \Z
	&&
	x \mapsto  \#\{ \inta \in (0, k] : [\lina_\inta/\lina_{\inta-1}] = x  \}
	\label{eq_def_jhvector}
	\end{align}
where $[\lina_\inta/\lina_{\inta-1}]$ is the isomorphism class of $\lina_\inta/\lina_{\inta-1}$ and $\smic$ is the family of isomorphism classes of simple objects in $\ecata$.  In fact, every maximal chain in $\Sub_\oba$ engenders the same function, \emph{cf.} \cite[\S 6.1.6, p. 250]{grandis12}. Thus we may refer to \eqref{eq_def_jhvector}  unambiguously as the \emph{Jordan-H\"older vector} of $\oba$, denoted $\jhv{\oba}$.  

\begin{lemma}  
If all objects in $\ecata$ have finite height then the map 
	\begin{align*}
	\amitb(\ecata) \to \Z^{\smic}
	&&
	[\oba] \mapsto \jhv{\oba}
	\end{align*}
is an embedding of monoids.
\end{lemma}
\begin{proof}
It suffices to show that $\jhv{ \oba \oplus \obb } = \jhv{ \oba } + \jhv{ \obb }$ for any pair of objects $\oba, \obb$.  If $\lina = \{0 = \lina_0 < \cdots < \lina_k = \oba \}$ and $\lina' = \{0 = \lina'_0 < \cdots < \lina'_l = \obb \}$ are maximal chains in $\oba$ and $\obb$ then, by a routine exercise, the combined set $\{\lina_0  < \cdots < \lina_k \oplus \lina'_0 < \cdots < \lina_k \oplus \lina'_l = \oba \oplus \obb \}$ is a maximal chain in $\oba \oplus \obb$.  The desired conclusion follows.
\end{proof}

\begin{proposition} 
\label{prop_finitecf}
Let $\oba$ be an object of finite length in $\ecata$.  If $\pseta$ is a partially ordered set and $\cdf: \axl{\pseta} \to \Sub(\oba)$ is a complete lattice homomorphism then
	\begin{align}
	\jhv{\oba} = \sum_{\psela \in \pseta}  
	\jhv{  \sfun \{ \psela \} }
	\label{eq_pseltocompfac}
	\end{align}
where $\sfun$ is the locally closed functor of $\cdf$.
\end{proposition}
\begin{proof}
Let $\lina$ be a linearization of $\pseta$, and let $\toc = \Im( \cdf|_{\axl{\lina}})$.  Since $\toc$ is finite and the restriction $\cdf|_{\axl{\lina}}$ is complete, the inverse image of each $\lela \in \toc$ under $\cdf|_{\axl{\lina}}$ is a closed interval, which we can denote $[\seta_\lela, \setb_\lela]_\lina$.  If $\seta \lessdot_\lina \setb$ is a covering relation in $\lina$, it then follows that $\cdf \seta < \cdf \setb$ iff $(\seta, \setb) = (\setb_\lela , \seta_\lelb)$ for some $\lela \lessdot_\toc \lelb$.    Thus	
	\begin{align}
	\sum_{\seta \lessdot_\lina \setb} \jhv{\cdf_\setb / \cdf_\seta}
	= 
	\sum_{\lela \lessdot_\toc \lelb} \jhv{\lelb / \lela}
	= 
	\jhv {\oba}.
	\label{eq_covertocompfac}
	\end{align}
Since every covering relation in $\axl{\lina}$ has form $\pdsh \psela \lessdot \dsh \psela$ for some $\psela \in \pseta$, the lefthand side of \eqref{eq_covertocompfac} coincides exactly with the righthand side of \eqref{eq_pseltocompfac}. The desired conclusion follows.
\end{proof}

\begin{theorem} 
Every constructible functor from $\R$ to $\ecata$ has a saecular CDF homomorphism.
\end{theorem}
\begin{proof}
By construction, the kernel and image filtrations on $\Sub(\mora)$ are finite.  Thus they generate a finite distributive (hence completely distributive) sublattice of $\Sub(\mora)$.
\end{proof}

Given an interval functor $\dgb: \R \to \ecata$ of support type $\itva$, we can define $\jhv{\dgb}$ as $\jhv{\dgb_\a}$ for any $\a \in \itva$, or, equivalently, as $\max_{\aa \in \itva} \jhv{\dgb_\a}$.

\begin{theorem}  
\label{thm_amit_pd_is_jh}
If $\dga$ is a constructible functor from $\R$ to an abelian category $\ecata$, then
	\begin{align}
	\pdtb{\dga} =  \jhvh \circ \spd_\dga
	\label{eq_mobiusinverse}
	\end{align}
where $\spd_\dga$ is the saecular persistence diagram of $\mora$.  That is, $\pdtb{\dga}(\itva)$ is the Jordan-H\"older vector of (the object type of) of the interval functor $\spd_\dga(\itva)$.
\end{theorem}
\begin{proof}
Let $[\lela, \lelb)$ be given, and put $\aa = \dopiso^{-1}(\a)$.  Let $\sfun$ be the saecular CD functor. Choose a real number $\lelc$  such that $\{ s \in \seta : s \le c\} = \{ s \in \seta : s < b\}$.  Then
	$
	\prfb{\dga} [a,b) = [\sfun(\aa)_\lelc]
	$.

The image of the saecular CD homomorphism, $\latb$, is a finite distributive sublattice of $\Sub_\dga$.  Thus the self map $q:  x \mapsto x \wedge \KK(\aa)$ is a $\mathring \forall$ complete lattice homomorphism $\latb \to \latb$.  It is therefore simple to verify that the composite map
	\begin{align}
	\cdf 
	= 
	(q \circ \sfun) : \axlh \allitv {\R} \to \Sub( \sfun( \aa ))
	\quad 
	\quad
	\setc \mapsto  \KK( \aa )  \wedge  \sfun \setc = \sfun(  \nu( \aa ) \wedge \setc )
	\label{eq_rank_fun_lattice_hom_def}
	\end{align}
is a $\forall$-complete lattice homomorphism, where, by abuse of notation, we write $\nu(\aa)$ for $(\iid \ii \nu)(\aa)$.  By a further abuse, we  write $\cdf$ for the locally closed functor of $\cdf$.    Since $\cdf$ is $\forall$-complete, Proposition \ref{prop_finitecf} implies the rightmost equality in 
	\begin{align}
	\prfb{\dga} [a,b) 
	= 
	[\sfun(\aa)_\lelc]
	= 
	\jhv{\sfun(\aa)_\lelc} 
	= 
	\sum_{\itva \in \allitv{\R}} 
		\jhv	{\ \cdf \{ \itva \}_\lelc
		}	
	\label{eq_jh_summands_basic}
	\end{align}
	
By inspection of \eqref{eq_rank_fun_lattice_hom_def}, $\itva \notin \nu(\aa) \implies \nu(\aa) \wedge \dsh \itva = \nu(\aa) \wedge \pdsh \itva$, hence 
	\begin{align}
	\cdf \{ \itva \}
	=
		\begin{cases}
		\sfun \{ \itva \} & \itva \in \nu( \aa) \\
		0 & else
		\end{cases}
	\label{eq_wedge_formulae}
	\end{align}
for each $\itva \in \allitv{ \R}$.	
\\

\noindent \ul{Claim 1.} $\sfun \{ \itva \}_\lelc = \cdf \{ \itva \}_\lelc = 0$ if interval $\itva$ cannot be expressed in form $[\ela_i, \ela_j)$ for some $\ela_i, \ela_j \in \seta$.  \emph{Proof.} If $X \in \{ \KK, \K\}$ then $X$ is constant on each interval $[\post_\inta, \post_{\inta+1})$, hence 
	$
	X(\pdsh \setb) \lneq  X(\dsh \setb) \implies \dsh \setb \in \{ \fem (x) : x \in \posts \}.
	$
Recalling that $\fem(x) = \axlh(\pdsh x)$, it follows that $\sfun$ vanishes on $\{\dop \setb \sqcup \setb\}$ unless $\dop \setb$ and $\setb$ are open intervals with suprema in $\seta$.
\\

\noindent \ul{Claim 2.} $\cdf \{ [\ela_i, \ela_j) \}_\lelc$ vanishes if $\ela_j < b$.  \emph{Proof.} 
In this case $\ela_j < c$.  Thus $\cdf  \dsh [\ela_i, \ela_j) \le \sfun  \dsh [\ela_i, \ela_j) \le \K( \ela_j)$.  Since $\K( \ela_j)_c = 0$ it follows that $\cdf  \dsh [\ela_i, \ela_j) = 0$, and the desired conclusion follows.
\\

\noindent \ul{Claim 3.} $\cdf \{ [\ela_i, \ela_j) \}_\lelc$ vanishes if $\a < \ela_i$.  \emph{Proof.} In this case $ \nu(\aa) \wedge \dsh [\ela_i, \ela_j) = \nu(\aa) \wedge \pdsh [\ela_i, \ela_j) = \nu(\aa) \wedge \nu(\ela_j)$, therefore $\cdf  \{ [\ela_i, \ela_j) \} = (\cdf \dsh   [\ela_i, \ela_j) /  (\cdf \pdsh   [\ela_i, \ela_j) = 0$.
\\

Combining \eqref{eq_wedge_formulae} with Claims 1-3, we find that 
	\begin{itemize*}
	\item $[a,b) \su \itva \implies \cdf \{\itva \} = \sfun \{\itva \}$
	\item $[a,b) \nsubseteq \itva \implies \cdf \{\itva \} = 0$
	\item $\itva \notin \pitvc \implies \cdf \{\itva\} = 0$
	\end{itemize*}
Therefore the righthand side of \eqref{eq_jh_summands_basic} equals
	\begin{align*}
	\sum_{\itva \in  Z }
		\jhv	{ \sfun \{ \itva \}_\lelc
		}	
	= 
	\sum_{\itva \in  Z }
		\jhv	{ \sfun \{ \itva \} }					
	= 
	\sum_{\itva \in  Z }
		(\jhvh \circ \spd )( \itva ) 
	.
	\end{align*}
where $Z = \{\itva \in \pitvc : [\lela, \lelb) \su \itva\}$.  Since $[a,b)$ was arbitrary, it follows that $(\jhvh \circ \spd)$ is  M\"obius inverse to $\prfb{\dga}$.  Equation \eqref{eq_mobiusinverse} follows by definition.
\end{proof}

\subsection{An elementary example}

Define $\qt: [0, \infty) \to \Z$ by
	\begin{align*}
	\qt_\tela
	=
		\begin{cases}
		0 		&		0 \le \tela < 1 \\
		4 		& 		1 \le \tela < 2 \\
		2 		& 		2 \le \tela
		\end{cases}
	\end{align*}
and let $\dga$ be the $\R$ shaped diagram in the category of abelian groups such that $\dga_\tela = 0$ for $\tela < 0$, and 
	\begin{align}
	\dga_\tela = \Z / \qt_\tela && \dga(\tela\le \telb)_\b ( x + \qt_\tela)  = x + \qt_\telb
	\label{eq_amit_rp2_example}
	\end{align}
for  $0 \le \tela \le \telb$.  The values of $\dga$ on objects are represented schematically in Figure \ref{fig_intervals}. 

Filtration $\KK$ then takes exactly two values, $0$ and $\dga$.  Filtration  $\K$ takes two nontrivial values, the  submodules with support  $[0, 1)$ and $[0,2)$.  Consequently $\supp(\prfa{\dga}) = \{ \linela^1, \linela^2, \linela^3 \}$, where
	\begin{align*}
	  \linela^1 = [0, 1) 
	  && 
	  \linela^2 =  [0, 2) 
	  &&
	  \linela^3=   [0, \infty).
	\end{align*}
This set admits a unique linear order $\linela^1 < \linela^2 < \linela^2$ compatible with inclusion.  The values of $\dgb^p$ and of the interval modules $\dgb^p / \dgb^{p-1}$ appear in Figure \ref{fig_intervals}.  

The lattice-theoretic approach allows us to associate \emph{generators} to the persistence diagram.   Suppose, for example, that $S$ is a unit circle, $\oba$ and $\obb$ are disks, and $\hk: \partial \oba \to S$ and $\kappa: \partial \obb \to S$ are continuous maps with winding numbers 4 and 2, respectively.  Let $X$ be the functor from $\R$ to the category of continuous maps and topological spaces such that 
	\begin{align*}
	X_\tela = 
		\begin{cases}
		\emptyset 		&  	\tela < 0 \\
		S 			& 	0 \le \tela < 1 \\
		S \sqcup_{\hk} \oba & 1 \le \tela < 2 \\
		S \sqcup_{\hk} \oba \sqcup_{\kappa} \obb & 2  \le \tela
		\end{cases}
	\end{align*}
and $X(\tela \le \telb)$ is inclusion for all $\tela \le \telb$.  Postcomposing $X$ with the degree 1 singular homology functor -- with integer coefficients -- yields a chain diagram isomorphic to $\dga$:
	\begin{align*}
	\dga \cong \funh_1 \com X .
	\end{align*}
Moreover, if $\generator$ is a fundamental class for $S$, then the nonzero objects in the cyclic quotient modules $\dgb^1/\dgb^0$, $\dgb^2/\dgb^1$, and $\dgb^3/\dgb^2$ are generated by $[4\generator]$, $[2\generator]$, and $[\generator]$, respectively.  In the parlance of topological data analysis, these generators represent three ``features'' of data born at time zero, that vanish at times 1, 2, and $\infty$, respectively.


\newcommand{\chaindgm}[8]{
	\def\xbase{#1}
	\def\ybase{#2}
	\def\xwidth{#3}
	\def\UL{#4}
	\def\tickmargey{#5}
	\def\ilabl{#6}
	\def\ilabc{#7}
	\def\ilabr{#8}	
	\draw (\xbase,\ybase) -- (\xbase+\xwidth, \ybase);

	\draw ({\xbase+1*\UL},\ybase+\tickmargey)-- ({\xbase+1*\UL},\ybase-\tickmargey);
	\draw ({\xbase+2*\UL},\ybase+\tickmargey)-- ({\xbase+2*\UL},\ybase-\tickmargey);
	\draw ({\xbase+3*\UL},\ybase+\tickmargey)-- ({\xbase+3*\UL},\ybase-\tickmargey);
			

	\node [below] at ({\xbase+1*\UL},\ybase-\tickmargey) {0};	
	\node [below] at ({\xbase+2*\UL},\ybase-\tickmargey) {1};
	\node [below] at ({\xbase+3*\UL},\ybase-\tickmargey) {2};	
	
	\node [above] at ({\xbase+0.5*\UL},\ybase) {$0$};
	\node [above] at ({\xbase+1.5*\UL},\ybase) {#6};
	\node [above] at ({\xbase+2.5*\UL},\ybase) {#7};	
	\node [above] at ({\xbase+3.5*\UL},\ybase) {#8};	
	
	\node  at (\xbase,\ybase) {(};	
	\node  at (\xbase+\xwidth,\ybase) {)};				
}

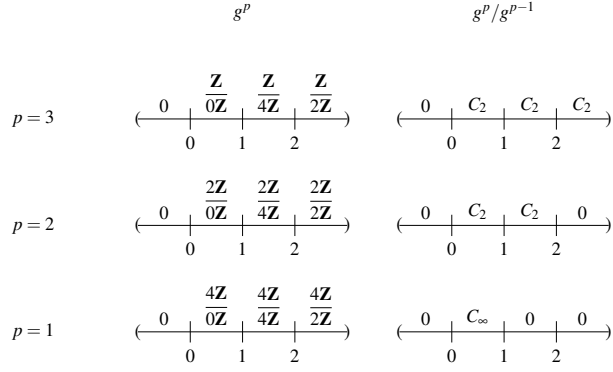
\begin{figure}
    \centering
  	\resizebox{0.5\textwidth}{0.3\textwidth}{      
        \begin{tikzpicture}
        
        \def\xbase{0}
        \def\xw{4}
        \def\xbasestep{5}
        
        \def\x{0}
        \chaindgm{\x}{4}{\xw}{1}{0.2}{$\dfrac{\Z}{0\Z}$}{$\dfrac{\Z}{4\Z}$}{$\dfrac{\Z}{2\Z}$}        
        \chaindgm{\x}{2}{\xw}{1}{0.2}{$\dfrac{2\Z}{0\Z}$}{$\dfrac{2\Z}{4\Z}$}{$\dfrac{2\Z}{2\Z}$}      
        \chaindgm{\x}{0}{\xw}{1}{0.2}{$\dfrac{4\Z}{0\Z}$}{$\dfrac{4\Z}{4\Z}$}{$\dfrac{4\Z}{2\Z}$}      

        \def\x{5}

        \chaindgm{\x}{4}{\xw}{1}{0.2}{$C_2$}{$C_2$}{$C_2$}        
        \chaindgm{\x}{2}{\xw}{1}{0.2}{$C_2$}{$C_2$}{$0$}      
        \chaindgm{\x}{0}{\xw}{1}{0.2}{$C_\infty$}{$0$}{$0$}   

        
        
        \def\h{6}        
	\node  at (0+\xw/2 ,\h) {$\dgb^p$};
	\node  at (5+\xw/2 ,\h) {$\dgb^p/\dgb^{p-1}$};

	\node  at (-2 ,0) {$p = 1$};
	\node  at (-2, 2) {$p = 2$};
	\node  at (-2, 4) {$p = 3$};			
	        
        \end{tikzpicture} 
    }
        
    \caption{Subquotients of chain diagram $\dga = \dgb^3$.} 
    \label{fig_intervals}
\end{figure}


\begin{figure}
    \centering
    \resizebox{0.7\textwidth}{0.29\textwidth}{  
        \begin{tikzpicture}
        
      \draw  (-1.5,-1) node[left]{$-1$} -- (1.5,-1) ;
      \draw (-1,-1.5) node[below] {$-1$} -- (-1,3.5) ;
    
      \foreach \x/\xtext in {0/0}
        \draw[shift={(\x,-1)}] (0pt,2pt) -- (0pt,-2pt) node[below] {$\xtext$};
    
      \foreach \y/\ytext in {0/0, 1/1, 2/2, 3/\infty}
        \draw[shift={(-1,\y)}] (2pt,0pt) -- (-2pt,0pt) node[left] {$\ytext$}; 	
        
	\node at (0,3) [circle,draw=black,fill=black,scale=0.3,label=right:{$\; \; \; \; \langle  \generator \rangle = \isqh(0,\infty) \cong C_2$}] {};  	          
	\node at (0,2) [circle,draw=black,fill=black,scale=0.3,label=right:{$\; \; \; \;  \langle 2 \generator \rangle \cong C_2$}] {};  
	\node at (0,1) [circle,draw=black,fill=black,scale=0.3,label=right:{$\; \; \; \; \langle 4 \generator \rangle \cong  C_\infty $}] {};

	\node at (4,1) [circle,fill=gray!45, scale=2]{$A$};		  
	\draw (4,1) circle (0.65cm);	

	\node at (7,1) [circle, scale=2]{$S$};		  
	\draw (7,1) circle (0.65cm);	
	
	\node at (10,1) [circle,fill=gray!45, scale=2]{$B$};		  
	\draw (10,1) circle (0.65cm);	

	\draw [->] (5,1) -- (6,1) node[midway,label=above:{$\hk$},label=below:{$\times4$}]{};
	\draw [->] (9,1) -- (8,1) node[midway,label=above:{$\kappa$},label=below:{$\times2$}]{};	
	
	\def\circledarrow#1#2#3{ 
	\draw[#1,->] (#2) +(80:#3) arc(80:-260:#3);
	}
	\circledarrow{dashed}{7,1}{0.8cm};
	
	\node at (6.4,2) []{$\generator$};	
	        
        \end{tikzpicture} 
    }
    \caption{
    \emph{Left.}  The type-$\amitb$ persistence diagram of chain functor $\dga$ defined in \eqref{eq_amit_rp2_example}.  
    \emph{Right.} Components of the identification space $S \coprod_{\hk} \oba \coprod_{\kappa} \obb $.  Element $\generator \in \funh_1(S ; \Z)$ is a fundamental class of $S$.  
    } 
    \label{fig_generators}
\end{figure}
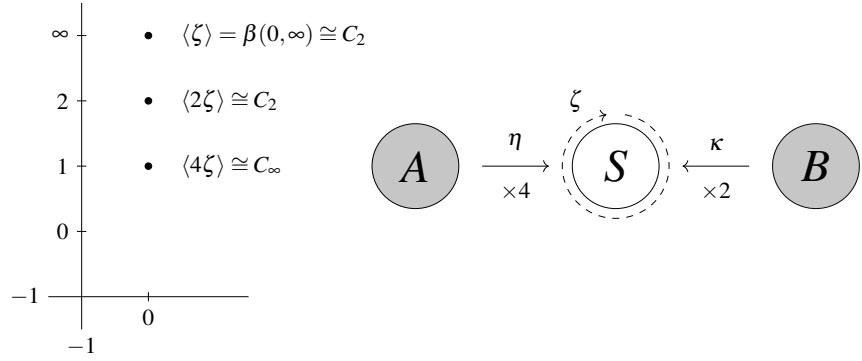

\section{The Leray-Serre Spectral Sequence}
\label{sec_lsss}

Let us resume the definitions and notation of \S\ref{chp:hompersistence}, with the additional assumption that  $\tob =  \{0, \ld, \intmaxa\}$.  For convenience,\footnote{This is a nonstandard assumption; it is easy to adapt to the case where $\fcxa(0) \neq 0$, but one must take care in reindexing.} suppose that $\fcxa$ preserves bounds as a poset map $\tob \to \Sub(\fcxa(\intmaxa))$.  

For each homology degree $\inta$, define  bound-preserving filtrations  $\filta^\inta: \tob \to \Sub_{\fcxa(\intmaxa)_\inta}$ and $\filtb^\inta: \toc = \{0, \ldots, 2\intmaxa + 1 \} \to \Sub_{\fcxa(\intmaxa)_\inta}$ by
	\begin{align*}
	\filta_\a^\inta
	=
	\fcxa(\a)_\inta
	&&
	\filtb_\a^\inta
	=
    	\begin{cases}
    	\partial_\bullet \filta^{\inta+1}_\a & \a \le \intmaxa \\
    	\partial^\bullet \filta^{\inta-1}_{\a -(\intmaxa+1)} & else 
    	\end{cases}
	\end{align*}
The superscripts on $\filta^\inta_\a$ and $\filtb^\inta_\a$ can often be inferred from context and thus suppressed.  For illustration, see Figure \ref{fig_lsss_sets}.

				\begin{figure}
                  \centering
                    \includegraphics[width=0.55\textwidth]{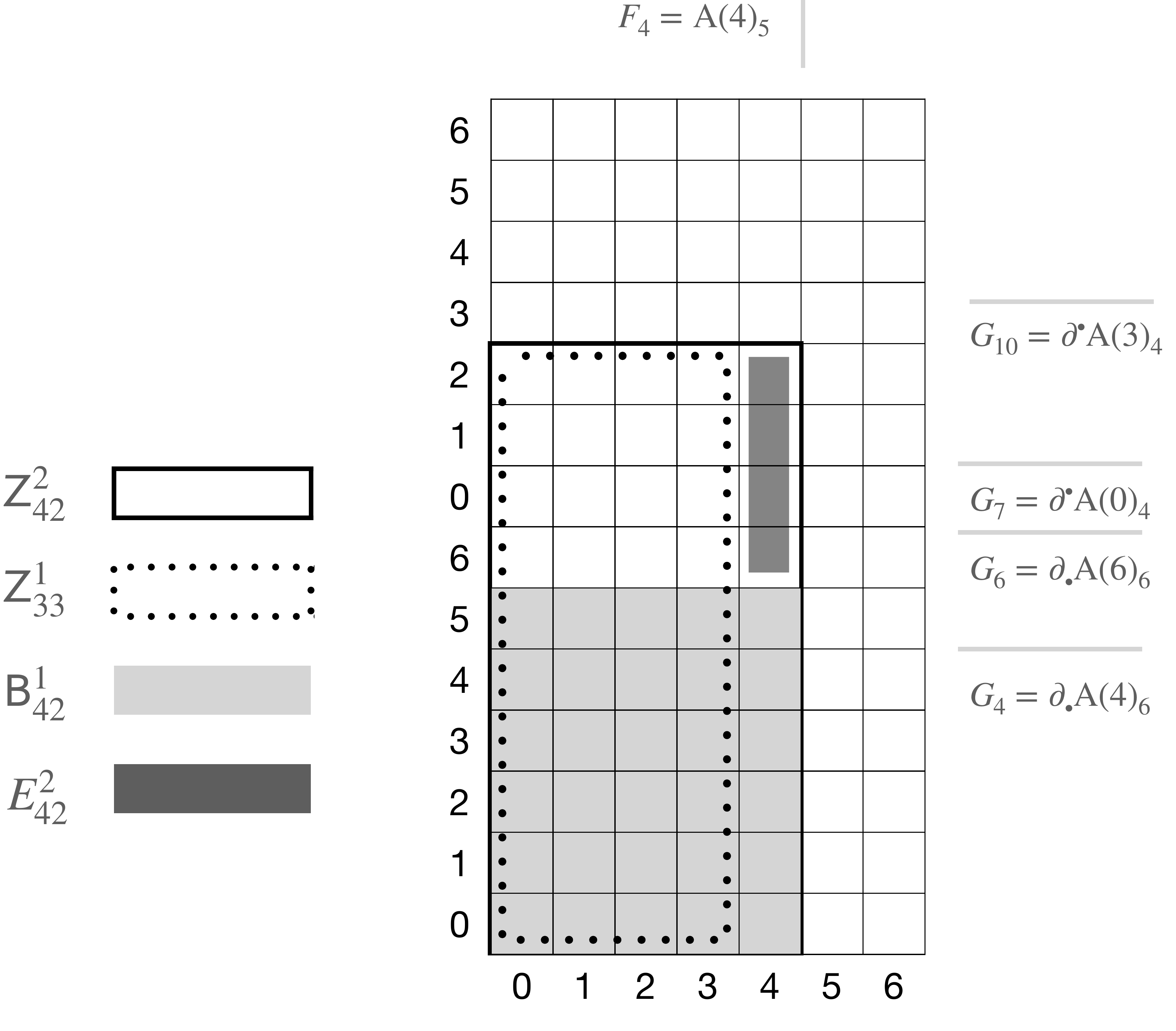}
                \caption{                	
					Filtrations of $\fcxa(\intmaxa)_6$. 
                }
                \label{fig_lsss_sets}
                \end{figure}

Write $\omega^\inta: \axl{ \tob \times \toc } \to \Sub_{\fcxa(\intmaxa)_\inta}$ for the (complete) lattice homomorphism such that $\omega^\inta( [0,p] \times \toc) = \filta^m_p$ and $\omega^\inta( \toc \times [0,q]) = \filtb^m_q$.  As usual,  write $\omega^\inta$ for the locally-closed functor of $\omega^\inta$.  

Define $\dga^\inta = \funh_\inta \circ \fcxa$ for each $\inta \in \Z$.  Since the index category of $\dga^\inta$ is finite, $\dga^\inta$ admits a saecular CDI functor $\sfun^\inta$. Let us say that a subquotient $X/Y$ is \emph{Noether sympathetic} with an interval functor $\sfun^\inta\{ \itva \}$ if $X/Y$ is Noether-isomorphic to $\sfun^\inta\{ \itva \}_\a$ for each $\a \in \itva$.

\begin{lemma}
\label{lem_proto_lsss_thm}
For each $0 \le x < y \le 2\intmaxa+1$, 
	\begin{align*}
	\omega^{\inta} \{(x,y)\}
		\begin{cases}
		\text{is Noether sympathetic to $\sfun \{ [x,y) \}$} & y \le \intmaxa + 1 \\
		\text{is Noether sympathetic to $\sfun \{ [y-\intmaxa-1, x) \}$} & else
		\end{cases}
	\end{align*}
\end{lemma}
\begin{proof}
The proof may be completed via two exercises.  Exercise 1: the case $y \le \intmaxa$ is essentially Theorem \eqref{eq_ph_cdf_formula}.  Exercise 2: when $y > \intmaxa$, the boundary operator $\partial$ carries $\omega^{\inta} \{(x,y)\}$ isomorphically onto $\omega^{\inta} \{(y-\intmaxa-1,x)\}$.
\end{proof}

Recall \cite[p 38]{grandis12} that the terms of the Leray-Serre spectral sequence satisfy
	\begin{align}
	\funz_{pq}^r (\fcxa) 
	& 
	= 
	\fcxa(p)_\inta \wedge \partial \ii \fcxa(p-r)_{\inta-1}  
	=
	\filta^\inta_p \wedge \filtb^\inta_{p-r + (\intmaxa + 1)}
	&&
	=
	\omega^\inta ([0,p] \times [0, p-r+(\intmaxa +1)])
	\label{eq_ls0}
	\\
	\funbd_{pq}^r(\fcxa) 
	&  
	= 
	\fcxa(p)_\inta \wedge \partial \di \fcxa(p+r)_{\inta+1} 
	=
	\filta^\inta_p \wedge \filtb^\inta_{p+r}
	&&
	=
	\omega^\inta ([0,p] \times [0, p-r])	
	\label{eq_ls1}
	\\
	E_{pq}^r(\fcxa) 
	& 
	= 
	\frac{\funz^r_{pq} }{ \funz^{r-1}_{p-1,q+1} \vee \funbd_{pq}^{r-1} }
	&&
	=
	\omega^\inta ( \{p \} \times [p+r, p-r + (\intmaxa + 1)])	
	\label{eq_ls2}
	\end{align}
where $\inta = p+q$.

				\begin{figure}
                  \centering
                    \includegraphics[width=0.6\textwidth]{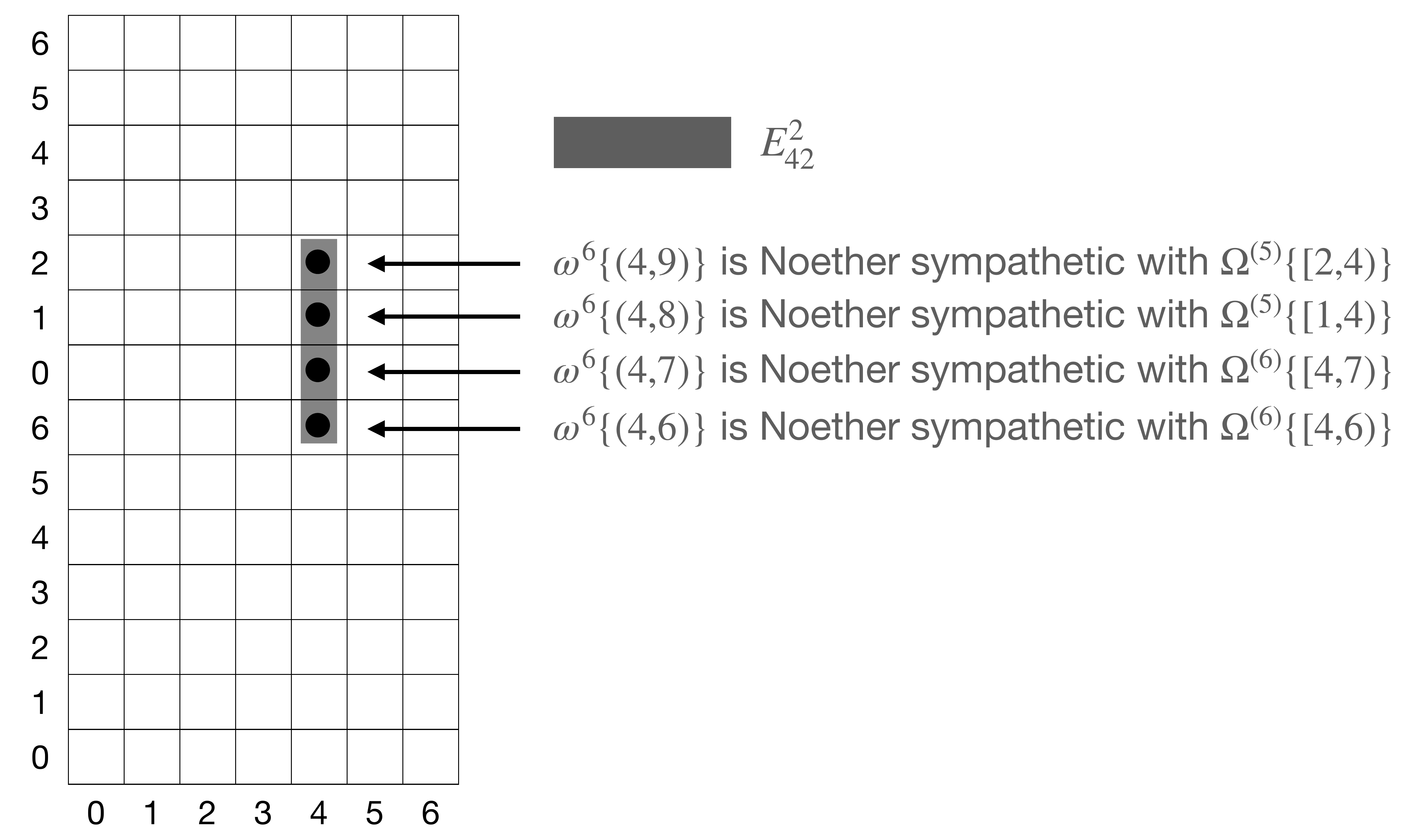}
                \caption{                	
					Values taken by the exact functor $\omega^6$ on singletons. 
                }
                \label{fig_lsss_points}
                \end{figure}

\begin{theorem}
Each term $E^r_{pq}$ of the Leray-Serre spectral sequence admits a unique filtration 
	\begin{align*}
	0 = \filtc_{p+r-1}  \le \cd \le \filtc_{p-r + (\intmaxa + 1)} = E_{pq}^r
	\end{align*}
such that 
	\begin{align*}
	\filtc_{x} / \filtc_{x-1}
		\begin{cases}
		\text{is Noether sympathetic with $\sfun^\inta \{ [x,y) \}$} & y \le \intmaxa + 1 \\
		\text{is Noether sympathetic with $\sfun^{\inta-1} \{ [y-(\intmaxa+1), x) \}$} & else
		\end{cases}
	\end{align*}
for each $ x\in [p + r , p-r + (\intmaxa +1)]$.
\end{theorem}
\begin{proof}
Follows  from \eqref{eq_ls2} and Lemma \ref{lem_proto_lsss_thm}.
\end{proof}

Now let
	\begin{align*}
	\tau^\inta(\itva)  : =  \jhvh ( \sfun ^\inta \{ \itva \})
	\end{align*}
be the Jordan-H\"older vector of (the object type of) $\sfun ^ \inta \{ \itva \}$.  In particular, $\tau^\inta(\itva)$ is a function $\amits(\ecata) \to \Z \cup \{\infty\}$. Elements of $(\Z \cup \{\infty\})^{\amits(\ecata)}$ can be added coordinatewise.  We say $\tau^\inta(\itva)$  is \emph{coordinatewise-finite} if every coordinate of $\tau^\inta(\itva)$  is an integer.

\begin{theorem}[Leray-Serre enumeration theorem for saecular persistence]  
\label{thm_leray_serre_count_saecular}
If  $\tau^\inta(\itva)$  is coordinatewise-finite, then 
	\begin{align}
	\jhvh (B^r_{pq}) & = \sum_{\max(\itva) \le  p-r} \tau^\inta (\itva) 	\label{eq_jhls0}\\
	\jhvh (Z^r_{pq}) & = \sum_{\min (\itva) \le p} \tau^\inta(\itva) + \sum_{ \min(\itvb) \le p-r, \; \max(\itvb) \le p} \tau^{\inta-1}(\itvb) \label{eq_jhls1}\\
	\jhvh (E^r_{pq}) & = \sum_{p \lneq s \le \intmaxa} \tau^\inta [p,s] + \sum_{s \le p - r} \tau^{\inta-1} [s,p] \label{eq_jhls2}
	\end{align}
for all $p$, $q$, and $r$.  Here intervals appearing as subscripts are compared via the partial order on $\allitv{\tob}$.
\end{theorem}
\begin{proof}
Apply Proposition \ref{prop_finitecf} to $\omega^\inta$.  Combine this result with  Equations \eqref{eq_ls0}-\eqref{eq_ls2} and Lemma \ref{lem_proto_lsss_thm}.
\end{proof}

There is a natural way to extend any chain diagram $\dga: \{0, \ldots, \intmaxa\} \to \ecata$ to an $\{0, \ldots, \intmaxa\}$-constructible functor $\dgb : \R \to \ecata$ such that $\dgc_\a = 0$ for $\a < 0$ and $\dgc_\a = \dgb_{\lfloor \a \rfloor}$ for $ a \in [0,\intmaxa+1]$.  For each $\dga^\inta$, apply this construction to obtain a constructible functor $\dgb^\inta$; denote the corresponding type-$\amitb$ generalized persistence diagram by $\gpdb_\inta$.

\begin{theorem}[Leray-Serre enumeration theorem for generalized persistence]
\label{thm_leray_serre_count_generalized}
Suppose that $\ecata$ is abelian and each vector $\tau_\inta(\itva)$ has finite support.  Then Theorem \ref{thm_leray_serre_count_saecular} remains true if we replace $\tau_\inta [a,b)$ with $\gpdb_\inta(a,b) $.
\end{theorem}
\begin{proof}
Theorem \ref{thm_amit_pd_is_jh} states that $\tau_\inta [a,b) = \gpdb_\inta(a,b)$ for all $a,b \in \R$.
\end{proof}

\section*{Acknowledgements}
%
This work would not have been possible without the assistance of Vidit Nanda. 
We also thank Hans Riess for useful comments.  

This work is principally supported by ONR N00014-16-1-2010 and NSF-1934960. G.\ H.-P. acknowledges   support from NSF grant DMS-1854748, the Centre for Topological Data Analysis of Oxford, EPSRC grant EP/R018472/1, and the Swartz Center for Theoretical Neuroscience at Princeton University.

\bibliography{bib}

\appendix

\section{Free modular and distributive lattices}	
\label{sec_modified_birkhoff}

\begin{theorem}[Birkhoff] 
\label{thm_birkhoff}
Let $\lata$ be a bounded lattice and $\pseta = \dop \tob \sqcup \tob$ be the coproduct in $\catps$ of two bounded totally ordered sets $\dop \tob$ and $\tob$.  Suppose that $\mora: \pseta \to \lata$ is an order-preserving function that preserves the bounds of both $\dop \tob$ and $\tob$.  Let $\stl{ \dop \tob_{>0} \times \tob_{>0}}$ be the sublattice of $\axl{ \dop \tob_{>0} \times \tob_{>0}}$ consisting of all finite unions of sets of form $(0, \aa] \times (0, \a]$ for some $\aa \in \dop \tob$  and $\a \in \tob$. Finally, define $\dop \cpph( \aa ) = [0, \aa] \times \tob$ and $\cpph(\a) = \dop \tob \times [0,\a]$.  Then the following are equivalent.
	\begin{enumerate*}
	\item There exists a modular sublattice of $\lata$ containing the image of $\psma$.
	\item There exists lattice homomorphism $\morc$ (which happens to preserve bounds) such that $\mora = \morc \fem$.  
	\end{enumerate*}	
\begin{equation} 
            \begin{tikzcd}
        				&&  \stl{\pseta} \arrow[dll , "\dop \cpph \sqcup \cpph" ' , <-] \arrow[drr, "\morc", dashed]	&& \\
        			\pseta \arrow[rrrr,  "\mora"] &&&& \cdla
        	\end{tikzcd}
    \label{eq_freediagram_modular}
\end{equation}  
\end{theorem}
\begin{proof}
The finite case is proved in \cite[p. 66, III, Theorem 9]{birkhoff1973lattice}; the infinite case is left as an exercise.  A complete treatment of the infinite case can be found in \cite[Theorem 1.3.7]{grandis12}.
\end{proof}

\begin{lemma}
\label{lem_fbd_alternate_form}
Lattice $\dxlh \axlh (\dop \tob \sqcup \tob)$ is isomorphic to the free bounded distributive lattice on $\axlh^2(\dop \tob) \sqcup \axlh^2(\tob)$.  Specifically, there exists a commutative diagram
\begin{equation} 
    \begin{tikzcd}
    & \axlh^2(\dop \tob) \sqcup \axlh^2(\tob )
    	\ar[dl, "" ]
    	\ar[dr, "\dop \cpph \sqcup \cpph"]		
	\\
	\stl{ \axlh^2(\dop \tob)_{>0} \times  \axlh^2(\dop \tob)_{>0} } 
		\ar[rr, "\rho" ]
	&&
	\dxlh \axlh( \dop \tob \sqcup \tob)
    \end{tikzcd}
    \label{eq_fbd_alternate_form}
\end{equation} 
where the lefthand diagonal arrow is defined as in  Theorem \label{thm_birkhoff}, the righthand diagonal arrow is defined as in \eqref{eq_free_diagram_xdi}, and $\rho: (0,\dop \seta] \times (0,\seta] \mapsto \dop \cpph(\dop \seta) \wedge \cpph(\seta)$.
\end{lemma}
\begin{proof}
The proof is a straightforward exercise in definition checking.
\end{proof}

\section{CD extension of well-ordered bifiltrations in a complete, modular, upper continuous lattice}
\label{sec_wo_cd_extension}

Here we prove Theorem \ref{thm_woextensioninmodularuc}.

Given a poset $\pseta$, let $\acl{\pseta}$ denote the lattice of all finite (including empty) antichains in $\pseta$, ordered such that $\seta \le \setb \iff \seta \su \dsh \setb$.  Let $\stl{\pseta}$ denote the sublattice of $\axl{\pseta}$ composed of all finite (including empty) unions of sets of form $\dsh \psela$.

\begin{lemma}
\label{lem_sviso}
For any poset $\pseta$ there exist mutually inverse isomorphisms
	\begin{align*}
	\isoa: \stl{\pseta} & \to \acl{\pseta}, \quad  \seta     \mapsto \max(\seta) \\
	\isob: \acl{\pseta} & \to \stl{\pseta}, \quad \achaina   \mapsto \dsh \achaina.
	\end{align*}
\end{lemma}
\begin{proof} 
Clearly $\isoa \isob =1$ and $\isob \isoa = 1$.
\end{proof}

\begin{lemma}
If $\tob$ is a well ordered set then so is every linearization of $\tob^2$.
\end{lemma}
\begin{proof}	
let $\lina$ be a linearization of $\tob^2$, and fix any nonempty subset $\seta \su \lina$.  The subset $\achaina = \min_{\tob^2}(\seta)$  is an antichain in $\tob^2$, so $\tela \lneq \tela'$ implies $\telb' \lneq \telb$ for any distinct elements $(\tela,\telb)$, $(\tela', \telb') \in \achaina$.  Since every nonempty subset $ \{\tela : (\tela, \telb) \in \achaina \}$ contains a  minimum, it follows that every subset of  $\{ \telb: (\tela, \telb) \in \achaina \}$ contains a maximum, and conversely.  Thus $\achaina$ is finite.  Every element of $\seta$ has an $\tob^2$-lower bound in $\achaina$, hence also an $\lina$-lower bound in $\achaina$.  Thus $\seta$ contains an $\lina$-minimum element in $\achaina$.
\end{proof}

\begin{corollary}  
\label{cor_wolinearization}
If $\tob$ is a well ordered set then so is every linearization of $\axlh(\tob)^2$.
\end{corollary}
\begin{proof}
It is elementary to show that $\tob$ is well ordered iff $\axl{\tob}$ is well ordered.  
\end{proof}

\begin{lemma}  
\label{lem_wofiniteantichains}
Let $\dop \tob$ and $\tob$ be totally ordered sets satisfying the ascending chain condition (ACC) or the descending chain condition (DCC), and set
	\begin{align}
	\pseta : = \dz{\axlh^2(\dop \tob)} \times \dz{\axlh^2(\tob)}
	\label{eq_monsterposet}
	\end{align}
Then there exists an isomorphism 
	\begin{align*}
	\isoa:
	\acl{ 
		\pseta
		}
	\to
	\axlh^2( \dop \tob \sqcup \tob)
	&&
	X
	\mapsto
	\bigcup_{(\dop \seta, \seta) \in X}  \{\dop \setb \sqcup \setb : \dop \setb \in \dop \seta, \; \setb \in \seta \}.		
	\end{align*}

\end{lemma}
\begin{proof}
Let $\setc \in \axlh^2(\dop \tob \sqcup \tob)$ be given.  Then the maps 
	\begin{align*}
	\polaa: & \axlh^2(\dop \tob) \to \axlh^2(\tob),  \quad \dop \seta \mapsto \{ \setb \in \axl{\tob} : \dop \setb \sqcup \setb \in \setc \; \forall  \; \dop \setb \in \dop \seta  \} \\
	\pola: & \axlh^2(\tob) \to \axlh^2( \dop \tob),  \quad \seta \mapsto \{ \dop \setb \in \axl{\dop \tob} :  \dop \setb \sqcup \setb \in \setc \; \forall  \;  \setb \in  \seta \} 
	\end{align*}
satisfy
	\begin{align*}	
	\polaa(\dop \seta ) \le \seta
	&&
	\iff
	&&
	\dop \setb \sqcup \setb \in \setc \; \; \forall \; (\dop \setb,  \setb) \in \dop \seta \times \seta
	&&
	\iff
	&&
	\pola(\seta) \le \dop \seta.
	\end{align*}  
In particular, $\polaa$ and $\pola$ are \emph{polarities}.  Consequently, these maps restrict to mutually inverse poset anti-isomorphisms between the images of the closure operators  $\clh := \polaa \pola$ and $\dop \clh : = \pola \polaa$.  If 
	\begin{align*}
	\achaina 
	&
	: = 
	\{ (\dop \seta, \polaa(\dop \seta)) : \dop \seta \in \Im(\clh^*) \} 
	= 
	\{ ( \pola(\seta), \seta) :  \seta \in \Im(\clh) \} 	 
	\end{align*}
and
	\begin{align*}
	\isob(\setc)
	&
	: = 
	\{ (\dop \seta, \seta) \in \achaina : \dop \seta \neq \emptyset, \; \seta  \neq \emptyset \}	
	\end{align*}
then by construction 
	\begin{align*}
	\setc 
	= 
	\isoa( \isob(\setc))
	\end{align*}
Moreover, since $\pseta$ is the product of two totally ordered sets satisfying ACC (respectively, DCC), and since $\achaina$ is an antichain in $\pseta$, only finitely many elements can lie in $\achaina$.  Thus
	\begin{align*}
	\isob(\setc) \in \acl{\pseta}.
	\end{align*}
Consequently $\isob$ yields a well-defined map $\acl{\pseta} \to \axlh^2(\dop \tob \sqcup \tob)$.  We have already seen that $\isoa \isob = 1$.  It is straightforward to check that $\isob \isoa = 1$, also.
\end{proof}

\begin{corollary} 
\label{cor_accAisS}
If $\dop \tob$ and $\tob$ are totally ordered sets satisfying ACC, then 
	\begin{align*}
	\axlh^2(\dop \tob \sqcup \tob)  = \stlh^2(\dop \tob \sqcup \tob).
	\end{align*}
\end{corollary}
\begin{proof}
Let $\isoa$ be defined as in Lemma \ref{lem_wofiniteantichains}.  
Since $\dz{\axlh^2(\dop \tob)}$ and $\dz{\axlh^2(\tob)}$ satisfy ACC, to each ordered pair $(\dop \seta, \seta)$ in the product \eqref{eq_monsterposet} corresponds a pair of  down-closed sets s $\dop \setb = \max \seta$ and $\setb = \max \seta$ such that
	\begin{align*}
	\phi( \{ \dop \seta, \seta \}) = \dsh ( \dop \setb \sqcup \setb).
	\end{align*}
It follows from surjectivity of $\isoa$ that every element of $\axlh^2(\dop \tob \sqcup \tob)$ can be expressed as a finite union of sets of form $\dsh (\dop \setb \sqcup \setb)$.
\end{proof}

\begin{corollary}  
\label{cor_accisalgebraic}
If $\tob$ and $\dop \tob$ satisfy ACC then $\axlh^2(\dop \tob \sqcup \tob)$ is algebraic.
\end{corollary}
\begin{proof}
Let $\setc \in \axlh^2(\dop \tob \sqcup \tob)$ be given.   One has $\setc = \dsh \max(\setc)$ by Corollary \ref{cor_accAisS}, so 
	\begin{align*}
	\textstyle
	\setc \le \bigvee \seta \iff \max(\setc) \su \bigcup \seta.
	\end{align*}
for any $\seta \su \axlh^2(\dop \tob \sqcup \tob)$. Since $\max (\setc)$ is finite, it follows that $\setc$ is compact.
\end{proof}

\begin{lemma} 
\label{lem_algebrasavesjoins}
Every lattice homomorphism $\mora$ from an algebraic lattice $\dla$ to an arbitrary lattice $\cdla$ preserves infinite joins.
\end{lemma}
\begin{proof}
Let $\seta$ be a nonempty subset of $\dla$.  By hypothesis there exists a finite subset $\setb \su \seta$ such that $\bigvee \seta = \bigvee \setb$.  One has
	\begin{align*}
	\textstyle
	\mora(\bigvee \seta)  = \mora(\bigvee \setb) = \bigvee \mora(\setb) \le \bigvee \mora(\seta) 
	\end{align*}
The reverse inequality holds by definition of least upper bound, so strict equality holds throughout.
\end{proof}

\begin{theorem}[Grandis] 
\label{thm_grandisbirkhoff}
Let $\KKs: \dop \tob \to \lata$ and $\Ks: \tob \to \lata$ be order-preserving maps from totally ordered sets to a modular lattice $\lata$.  Then there exists a unique lattice homomorphism
	\begin{align*}
	\cdf: \stl{\dz {\dop \tob} \times \dz\tob} \to \lata
	\end{align*}
such that
	\begin{align}
	\cdf( \dsh ( \dop \lela, \lela)) = \KKs(\dop \lela) \wedge \Ks(\lela)
	\label{eq_freedistlatticecondition}
	\end{align}
for each $(\dop \lela, \lela) \in \dz {\dop \tob} \times \tob$.
\end{theorem}

\begin{remark}
\label{rmk_anti-isomorphism}
To each poset $\psetc$ corresponds an anti-isomorphism
	\begin{align*}
	\isoa: \axl{\psetc} \to \axl{\psetc^*},
	&&
	\seta 
	\mapsto 
	\psetc - \seta.
	\end{align*}
Through repeated application, this construction can produce an anti-isomorphism $\isoa^{(\inta)}: \axlh^\inta(\pseta) \to \axlh^\inta(\pseta^*)$ for any positive $\inta$.
\end{remark}

\begin{theorem}  
\label{thm_woextensioninmodularuc}
Let $\dop \tob$ and $\tob$ be well-ordered sublattices of a complete upper continuous lattice $\lata$.  Let $\dop \toc$ and $ \toc$ denote the smallest $\forall$-complete sublattices containing $\dop \tob$ and $\tob$, respectively.  If there exists a modular sublattice $\latb \su \lata$ such that $\dop \toc \cup \toc \su \latb$, then there exists a $\forall$-complete CD sublattice $\latb' \su \lata$ such that $\dop \toc \sqcup \toc \su \latb'$.
\end{theorem}
\begin{proof}
Let $\KKs : \dop \tob \to \lata$ and $\Ks: \tob \to \lata$ be the inclusions of two well-ordered sets.  Since $\lata$ is complete, each of these maps factors through a complete totally ordered (hence completely distributive) sublattice.  Thus $\KKs$ and $\K$ have CDF extensions
	\begin{align*}
	\KKs : \axlh^2( \dop \tob) \to \lata &&
	\Ks  :  \axlh^2(      \tob) \to \lata.
	\end{align*}
By Theorem \ref{thm_grandisbirkhoff}, this pair of maps extends to a lattice homomorphism
	\begin{align*}
	\cdf: \stl{\dz{\axlh^2(\dop \tob)} \times \dz{\axlh^2(\tob)}} \to \lata
	\end{align*}
 satisfying \eqref{eq_freedistlatticecondition}.	Lemmas \ref{lem_sviso} and \ref{lem_wofiniteantichains} yield an isomorphism
	\begin{align*}
	\isoc: 
	\stl{\dz{\axlh^2(\dop \tob)} \times \dz{\axlh^2(\tob)}} 
	\xrightarrow{\cong }
	\axlh^2(\dop \tob \sqcup \tob)
	\end{align*}
where 
	\begin{align}
	\isoc(\dsh(\dop \seta, \seta)) : = \{ \dop \setb \sqcup \setb : \dop \setb \in \dop \seta, \; \setb \in \seta \}.
	\label{eq_finitebuildingblocks}
	\end{align}  
  Every element of the complete lattice $\axlh^2(\dop \tob \sqcup \tob)$ can  be expressed as a finite join of elements of form \eqref{eq_finitebuildingblocks}.  By abuse of notation, we will identify $\cdf$ with $\cdf \isoc^{-1}$.

  We claim that 
	\begin{align*}
	\cdf (\dsh(\dop \seta, \seta)) 
	= 
	\bigvee_{\dop \setb \in \dop \seta, \; \setb \in \seta} \cdf (\dsh( \dsh \dop \setb, \dsh \setb))
	=
	\bigvee_{\dop \setb \in \dop \seta, \; \setb \in \seta} \KKs(\dsh \dop \setb) \wedge \Ks(\dsh \setb)
	\end{align*}
for each $\dop \seta \in \dz{\axlh^2(\dop \tob)}$ and each $ \seta \in \dz{\axlh^2( \tob)}$.  This identity can be checked in three cases.  In the first case, suppose that $\dop \seta$ and $\seta$ each contain a maximum element.  Then the desired conclusion is vacuous.  In the second case suppose that at least one of $\dop \seta$ and $\seta$ contains a maximum element.  Here the desired conclusion holds by upper continuity.  In the third case, allow that neither $\dop \seta$ nor $\seta$ may contain a maximum.  Then a second application of upper continuity, applied to the already-established second case, yields desired conclusion, and the identity holds in full generality.

Now fix an arbitrary  subset $\sete \su \axlh^2(\dop \tob \sqcup \tob)$, and choose a finite sequence of pairs $(\dop \seta_\inta, \seta_\inta)$ such that $\bigvee \sete = \bigcup_\inta \isoc(\dsh (\dop \seta, \seta))$.  One then has
	\begin{align*}
	\cdf ( {\textstyle \bigvee \sete }) 
	= 
	\bigvee_\inta \cdf \dsh ( \dop \seta_\inta , \seta_\inta) 
	= 
	\bigvee_\inta
	\bigvee_{(\dop \setb, \setb) \in \dop \seta_\inta \times \seta_\inta} \cdf (\dsh( \dsh \dop \setb, \dsh \setb))
	\le
	\bigvee \cdf ( \sete ).
	\end{align*}
The reverse inequality holds by the definition of join, so $\cdf(\bigvee \sete) = \bigvee \cdf(\sete)$.  As $\sete$ was arbitrary, it follows that $\cdf$ preserves arbitrary joins.

Since $\axlh^2(\dop \tob \sqcup \tob)$ is anti-isomorphic to $\axlh^2(\dop \tob^* \sqcup \tob^*)$ by Remark \ref{rmk_anti-isomorphism}, and the latter is algebraic by Corollary \ref{cor_accisalgebraic}, it follows from Lemma \ref{lem_algebrasavesjoins} that $\cdf$ preserves arbitrary nonempty meets.  It preserves the empty meet by construction.  Therefore $\cdf$ is complete.  The image of $\cdf$ is thus a complete, completely distributive sublattice of $\lata$ containing $\dop \tob$ and $\tob$.  The desired conclusion follows.
\end{proof}

\section{Proof of Theorems \ref{thm_sample_hom} - \ref{thm_build_omega}}
\label{sec_intro_proofs}

\begin{proof}[Proof of Theorem  \ref{thm_sample_hom} ]
The lattice of subobjects in any category of $R$ modules is complete and upper-continuous.  Therefore existence of the CD homomorphism follows from the well-ordered criterion in Theorem \ref{thm_saecularexistencecriteria}.
\end{proof}

\begin{proof}[Proof of Theorem  \ref{thm_main_result_finite_fun} ]
Existence and uniqueness follow from Theorems \ref{thm_sample_hom} and \ref{thm_lclcextension}.
\end{proof}

\begin{proof}[Proof of Remark  \ref{rmk_sample_cd_finite_existence} ]
When $\tob$ is finite, each saecular filtration takes finitely many values.  The saecular BD homomorphism exists, and coincides with the saecular CD homomrophism.  If each $\dga_\a$ has finite height, then the saecular CD homomorphism exists for similar reasons.  Note, in particular, that we do not require $\Sub(\dga_\a)$ to be a complete lattice, for any $\a$.
\end{proof}

\begin{proof}[Proof of Theorem  \ref{thm_sample_hom_group} ]
To establish existence of the saecular CD homomorphism, first note that  $\dga$ admits a saecular BDH, by Theorem \ref{thm_saecular_functor_exiests_for_group_iff_images_normal}.  The image of this homomorphism is a distributive (a fortiori, modular) sublattice of $\Sub_\dga$ which contains the images of the FCD homomorphisms $\KK(\dop \tob): \axlh^2(\dop \tob) \to \Sub_\dga$ and $\K( \tob): \axlh^2( \tob) \to \Sub_\dga$.  Therefore the saecular filtrations factor through a $\forall$-complete CD sublattice of $\Sub(\dga)$, by Theorem \ref{thm_woextensioninmodularuc}.  Thus the SCDH exists.

The SCDH is natural, by \ref{thm_cdf_is_proper_group}.  The saecular coset factors and cokernel factors of $\dga$ are therefore interval functors, by Theorem \ref{thm_single2interval_group}.  
\end{proof}

\begin{proof}[Proof of Theorem \ref{thm_sample_cpd_agrees}]
Follows from Theorem \ref{thm_saecular_pd_generalizes_linear_pd}.
\end{proof}

\begin{proof}[Proof of Theorem \ref{thm_sample_gpd_agrees}]
This is identical to Theorem \ref{thm_amit_pd_is_jh}.
\end{proof}

\begin{proof}[Proof of Theorem \ref{thm_cseries_sample}]
Follows from Theorems \ref{thm_cseries} and \ref{thm_cseries_group}.
\end{proof}

\begin{proof}[Proof of Theorem \ref{thm_sample_saecular_factors_unique}]
If $\ecata$ is p-exact, then this result follows from exactness of the SCDI functor.  If $\ecata = \catgroup$, then the desired conclusion follows from property \eqref{eq_define_quasiregular_functor} of the quasi-regular functor; in particular, this property implies isomorphism between $\frac{|\dsh_{\allitv{\tob}} \itva|}{|\pdsh_{\allitv{\tob}} \itva|}$ and $\frac{|\dsh_{\lina} \itva|}{|\pdsh_{\lina} \itva|}$.
\end{proof}

\begin{proof}[Proof of Theorem \ref{thm_build_omega} and Remark \ref{rmk_generalization_of_interval_join_formula}]
%
If $\pseta$ is a poset, then each down-set $\seta \in \axl{ \pseta }$ is a join (i.e.\ union) of sets of form 
	$
	\seta = \bigvee_{\psela \in \seta}\dsh_\pseta (\psela) 
	$
As a special case, when $\pseta = \axl{ \dop \tob \sqcup \tob}$ or $\pseta =  \alljtv{\dop \tob \sqcup \tob}$, each element of $\axlh(\pseta)$ can be expressed  as a join of sets of form 
	$
	\seta = \bigvee_{\dop \setb \sqcup \setb \in \seta} \dsh_{\pseta} (\dop \setb \sqcup \setb)
	$.
Since $\sfuni: = \SCDI(\dga)$ and $\sfun: = \SCD(\dga)$ are $\forall$-complete lattice homomorphisms, it follows that
	\begin{align*}
	\textstyle
	\sfuni \seta 
	= 
	\sfuni \bigvee_{\dop \setb \sqcup \setb \in \seta} \dsh_{\pseta} (\dop \setb \sqcup \setb)
	= 
	\bigvee_{\dop \setb \sqcup \setb \in \seta} \sfuni \dsh_{\pseta} (\dop \setb \sqcup \setb)	
	= 
	\bigvee_{\dop \setb \sqcup \setb \in \seta} \sfun \dsh_{\pseta} (\dop \setb \sqcup \setb)	
	\end{align*}	
Since we identify $\allitv{\tob}$ with $\alljtv{\dop \tob \sqcup \tob}$, the desired conclusion follows from the observation that for any $\itva = \dop \setb - \dopiso(\setb)$, 
	\begin{align*}
	\sfun \dsh_{\pseta} (\dop \setb \sqcup \setb) 
	= 
	\KK(\dsh \dop \setb) \wedge \K( \dsh \setb) 
	= 
	\textstyle
	( \bigwedge_{\aa \notin \dop \setb}\KK(\aa)) \wedge ( \bigwedge_{\a \notin \setb}\K(\a))
	= 	
	\At_{\dopiso( \dop \setb), \setb}
	\end{align*}
\end{proof}


\end{document}